\documentclass[11pt]{article}               
\date{~}
\usepackage{dratex,calc,amsmath,amsthm,amscd,amssymb,amsfonts}
\usepackage[all]{xy}
\UseComputerModernTips
\theoremstyle{plain}
	\newtheorem{theorem}{Theorem}[section]

	\newtheorem{lemma}[theorem]{Lemma}
	\newtheorem{proposition}[theorem]{Proposition}
	\newtheorem{corollary}[theorem]{Corollary}

\theoremstyle{definition}
	\newtheorem{example}[theorem]{Example}
	\newtheorem{definition}[theorem]{Definition}
	\newtheorem{remark}[theorem]{Remark}
	\newtheorem{caveat}[theorem]{Caveat}
	\newtheorem{notation}[theorem]{Notation}
\def\cy#1{\mathcal{#1}}

\def\functor{\,\underline{\phantom{A}}\,} 
\def\Z{\mathbb{Z}}

\def\R{\mathbb{R}}
\def\Q{\mathbb{Q}}
\def\C{\mathbb{C}}

\def\ab{{\rm ab}}

\def\rat{{\rm rat}}
\def\proj{\text{{\rm -Proj}}}

\DeclareMathOperator{\Rank}{Rank}

\DeclareMathOperator{\Coker}{Coker}
\DeclareMathOperator{\Ker}{Ker}
\DeclareMathOperator{\Image}{Image}

\DeclareMathOperator{\End}{End}
\DeclareMathOperator{\Hom}{Hom}
\DeclareMathOperator{\Ext}{Ext}

\DeclareMathOperator{\id}{id}

\DeclareMathOperator{\Fix}{Fix}
\def\Bl{{\cal F}{\rm lk}}
\def\Sei{{\cal S}{\rm ei}}
\def\Prim{{\cal P}{\rm rim}}
\begin{document}
\begin{center}
{\Large Invariants of boundary link cobordism II. \\
The Blanchfield-Duval form \\[2ex]}
{\large Desmond Sheiham}
\end{center}
\vspace*{3ex}
\begin{abstract}
We use the Blanchfield-Duval form to define complete invariants for
the cobordism group $C_{2q-1}(F_\mu)$ of $(2q-1)$-dimensional
$\mu$-component boundary links (for $q\geq2$).

The author solved the same problem in earlier work via Seifert
forms. Although Seifert forms are convenient in explicit computations,
the Blanchfield-Duval form is more intrinsic and appears naturally in
homology surgery theory.

The free cover of the complement of a
link is constructed by pasting together infinitely many copies of the
complement of a $\mu$-component Seifert surface. 
We prove that the algebraic analogue of this construction, a functor
denoted $B$, identifies the author's earlier invariants with those defined
here. We show that $B$ is equivalent to a universal localization of
categories and describe the structure of the modules sent to
zero. Taking coefficients in a semi-simple Artinian ring, we deduce that
the Witt group of Seifert forms is isomorphic to the Witt group of
Blanchfield-Duval forms. 
\end{abstract}
\section{Introduction}{\def\thefootnote{}\footnote{Submitted September
23rd 2003.}}\addtocounter{footnote}{-1}
This paper is the second in a series on cobordism (=concordance) groups
of a natural class of high-dimensional links.
Chapter $1$ of the first work~\cite{She03mem} discusses background to
the problem at greater length but we summarize here some of the key
ideas.

\subsection{Background}
A knot is an embedding of spheres\footnote{Manifolds are assumed
oriented and embeddings are assumed locally flat. One may work in the
 category of smooth, $PL$ or topological manifolds according to taste,
 with the understanding that $S^n$ is permitted exotic structures if
 one selects the smooth category.}
 $S^n\subset S^{n+2}$. The following are generalizations:
\begin{itemize}
\item A $\mu$-component {\it link} is an embedding of $\mu$ disjoint
spheres
\begin{equation*}
L=\overbrace{S^n\sqcup\cdots \sqcup S^n}^{\mu} \subset S^{n+2}.
\end{equation*}
\item A {\it boundary link} is a link whose components bound disjoint
$(n+1)$-manifolds. The union of these $(n+1)$-manifolds
is called a {\it Seifert surface}.
\item An {\it $F_\mu$-link} is a pair $(L,\theta)$
where $L$ is a link and 
$\theta$ is a homomorphism from the fundamental group
$\pi_1(X)$ of the link complement $X=S^{n+2}\backslash L$ onto the free group
$F_\mu$ on $\mu$ (distinguished) generators such that some meridian of
the $i$th link component is sent to the $i$th generator. 
\end{itemize}
Not every link is a boundary link;
a link $L$ can be refined to an $F_\mu$-link $(L,\theta)$ if and only if
$L$ is a boundary link. 

Let us call a homomorphism $\theta:\pi_1(X)\to F_\mu$ permissible if it
sends some meridian of the $i$th link component to the $i$th generator.
There may be many permissible homomorphisms for a given boundary
link but if $\theta$ and $\theta'$ are permissible then
$\theta'=\alpha\theta$ where $\alpha$ is some ``generator conjugating''
automorphism of $F_\mu$ (Cappell and Shaneson~\cite{CapSha80}, Ko~\cite[p660-663]{Ko87}). Homomorphisms $\pi_1(X)\to F_\mu$
correspond to homotopy classes of maps from the link complement $X$ to
a wedge of $\mu$ circles and the permissible homotopy classes
correspond, by the Pontrjagin-Thom construction, 
to cobordism classes of Seifert surfaces (rel $L$). 

Every knot is a ($1$-component) boundary link and admits precisely one
permissible homomorphism, namely the abelianization
\begin{equation*}
\theta:\pi_1(X)\to\pi_1(X)^\ab\cong\Z.
\end{equation*}
Among the three
generalizations above it is the theory of $F_\mu$-links which seems to
bear the closest resemblance to knot theory.

Although one does not hope for a complete classification of knots or
$F_\mu$-links in higher dimensions much is known about their classification
up to the equivalence relation known as cobordism (or concordance).
Two links $L^0$ and $L^1$ are called {\it cobordant} if there is an embedding
\begin{equation*}
LI=(S^n\sqcup\cdots \sqcup S^n) \times [0,1] \subset S^{n+2}\times[0,1]
\end{equation*}
which joins $L^0\subset S^{n+2}\times\{0\}$ to $L^1\subset
S^{n+2}\times\{1\}$. One requires\footnote{$LI$ is also required to
meet $S^{n+2}\times\{0\}$ and $S^{n+2}\times\{1\}$ transversely.} that
$(S^n\sqcup\cdots \sqcup S^n) \times\{i\} \subset S^{n+2}\times\{i\}$
for $i=0$ and $i=1$ but no such requirement is made when $0<i<1$. Boundary links are said to be {\it boundary cobordant} if there is a
cobordism $LI$ whose components bound disjoint $(n+2)$-manifolds in
$S^{n+2}\times[0,1]$. Two $F_\mu$-links $(L^0,\theta^0)$ and
$(L^1,\theta^1)$ are called cobordant if there is a pair
\begin{equation*}
(LI\ ,\ \Theta:\pi_1(S^{n+2}\times [0,1] \backslash LI) \to F_\mu)
\end{equation*}
such that the restrictions of $\Theta$ to $\pi_1(X^0)$ and
$\pi_1(X^1)$ coincide with $\theta^0$ and $\theta^1$ (up to inner
automorphism).

The cobordism classes of knots form an abelian group $C_n(F_1)$ under (ambient)
connected sum but this operation does not extend to links in any
obvious way. If one attempts to add links $L^0$ and $L^1$
there are many inequivalent choices of connecting arc
from the $i$th component of $L^0$ to the $i$th component of $L^1$.

However when $n\geq2$ connected sum $[L_1,\theta_1]+[L_2,\theta_2]$ of
{\it cobordism classes} of $F_\mu$-links
is well-defined; one can remove the ambiguity in the choice of paths
by assuming, perhaps after some surgery, that $\theta_1$ and $\theta_2$ are
isomorphisms.
The set $C_n(F_\mu)$ of cobordism classes of $F_\mu$-links is therefore an
abelian group.

When $n$ is even, $C_n(F_\mu)$ is in fact the trivial
group~\cite{Ker65,CapSha80,Ko87,Mio87}; we sketch a proof
in~\cite[Ch1\S4.1]{She03mem}. On the other hand J.Levine obtained a
complete system of invariants for 
odd-dimensional knot cobordism groups $C_{2q-1}(F_1)$ for
$q\geq2$~\cite{Lev69B} and showed that each is isomorphic to a
countable direct sum
\begin{align}\label{iso_class_of_knot_cobordism_gp}
&C_{2q-1}(F_1)~\cong~\Z^{\oplus\infty}\oplus
\left(\frac{\Z}{2\Z}\right)^{\oplus\infty}
\oplus \left(\frac{\Z}{4\Z}\right)^{\oplus\infty}.
\intertext{The computation of $C_1(F_1)$ remains
open. In~\cite{She03mem} the author obtained a complete system 
of invariants for odd-dimensional $F_\mu$-link cobordism groups
$C_{2q-1}(F_\mu)$, $q\geq2$ (including some secondary invariants
defined only if certain primary invariants vanish) and found that}
\label{iso_class_of_F_link_cobordism_gp}
&C_{2q-1}(F_\mu)~\cong~\Z^{\oplus\infty}\oplus
\left(\frac{\Z}{2\Z}\right)^{\oplus\infty}\oplus\left(\frac{\Z}{4\Z}\right)^{\oplus\infty}
\oplus\left(\frac{\Z}{8\Z}\right)^{\oplus\infty}
\end{align}
for all $q\geq2$ and all $\mu\geq2$. 

Both~(\ref{iso_class_of_knot_cobordism_gp})
and~(\ref{iso_class_of_F_link_cobordism_gp}) were deduced from a purely
algebraic reformulation of $F_\mu$-link cobordism associated to
Seifert surfaces: It was proved by Levine~\cite{Lev69} in the knot
theory case $\mu=1$ and by Ko~\cite{Ko87} and Mio~\cite{Mio87} independently in the
general case that $C_{2q-1}(F_\mu)$ is isomorphic to the ``Witt group
of Seifert forms''. In the notation of the present paper, which we
explain more carefully in
Sections~\ref{section:define_hermitian_forms_and_witt_gp},~\ref{section:seifert_modules}
and~\ref{section:cobordism_invariants_seifert},
\begin{equation}\label{identify_F_link_cobordism_with_W(Seifert)}
C_{2q-1}(F_\mu)~\cong~W^{(-1)^q}(\Sei(\Z)) \quad\quad (q\geq3).
\end{equation}
The symbol $\Sei(\Z)$ denotes\footnote{The
category $\Sei(\Z)$ was  
denoted $(P_\mu\mbox{--}\Z)\proj$ in~\cite{She03mem}.} a category of
``Seifert modules'' designed to contain the homology modules of Seifert
surfaces among the objects (see Notation~\ref{Sei_notation}). In the
case $\mu=1$ an object in $\Sei(\Z)$ is a finitely
generated free $\Z$-module $V$ together with an endomorphism $V\to V$ which
carries information about how a Seifert surface is embedded. If
$\mu>1$ then the definition of Seifert module also includes a direct sum
decomposition $V=V_1\oplus\cdots\oplus V_\mu$ which reflects the
connected components of a Seifert surface.

The intersection form in a Seifert surface is an isomorphism
$\phi:V\to V^*$ in $\Sei(\Z)$ which satisfies
$\phi^*=(-1)^q\phi$. Such $(-1)^q$-hermitian forms are the 
generators of the Witt group $W^{(-1)^q}(\Sei(\Z))$. The relations
say that certain ``metabolic forms'' are identified with 
zero; see Definitions~\ref{define_metabolic} and~\ref{define_Witt_group} below.
 
Although Seifert surface methods are convenient in
explicit computations, it is preferable to define $F_\mu$-link invariants without
making a choice of Seifert surface. In the present paper we focus
instead on the covering space $\overline{X}\to X$ of a link
complement determined by the homomorphism
$\theta:\pi_1(X)\twoheadrightarrow F_\mu$. This approach sits more
naturally in homology surgery theory and is more
amenable to generalization from boundary links to arbitrary links
or other manifold embeddings.

We take as starting point the identification
\begin{equation}\label{identify_F_link_cobordism_with_W(Flink)}
C_{2q-1}(F_\mu)~\cong~W^{(-1)^{q+1}}(\Bl(\Z))\quad\quad (q\geq3)
\end{equation}
where $\Bl(\Z)$ is a category designed to contain homology modules
of the cover $\overline X$ (see
Definition~\ref{define_F_link_module} and
Notation~\ref{Flk_notation}). The objects in
$\Bl(\Z)$ are certain modules over the group ring
$\Z[F_\mu]$ of the free group; they are called
$F_\mu$-link modules in the present paper although they are more commonly
known as link modules. 

The $F_\mu$-equivariant Poincar\'e duality in ${\overline X}$
leads to a $(-1)^{q+1}$-hermitian form $\phi$ in
the category $\Bl(\Z)$. This is the Blanchfield-Duval form of the title,
originally introduced by Blanchfield~\cite{Bla57} in the knot theory
case $\mu=1$. The
identity~(\ref{identify_F_link_cobordism_with_W(Flink)}) was proved by
Kearton for $\mu=1$~\cite{Kea75,Kea75'} and by Duval~\cite{Duv86}
for $\mu\geq2$.
Cappell and Shaneson earlier identified the cobordism group
$C_n(F_\mu)$ with a $\Gamma$-group, an obstruction group in their homology
surgery theory~\cite{CapSha74,CapSha80}. The identification of this
$\Gamma$-group with the Witt group $W^{(-1)^{q+1}}(\Bl(\Z))$ was due to
Pardon~\cite{Par76,Par77}, Ranicki~\cite[\S7.9]{Ran81} and
Smith~\cite{Smi81} for $\mu=1$ and to
Duval~\cite{Duv86} for $\mu\geq2$. More general results of
Vogel~\cite{Vog80,Vog82} on homology surgery and universal
localization are stated elsewhere in this
volume~\cite[\S1.4]{Ran03'}. An outline 
of their application to $C_n(F_\mu)$ is given
in~\cite[Ch1,\S4.4,5.3]{She03mem}.

\subsection{Overview}
Universal localization plays two roles in this paper. Firstly the
 ``augmentation localization'' of the group ring $\Z[F_\mu]$ of the free group
 appears in the definition of the
 Blanchfield-Duval form, our main object of study. 
Secondly, we prove that the category
 $\Bl(\Z)$ of $F_\mu$-link modules is (equivalent to) a universal
 localization of the category $\Sei(\Z)$ of Seifert modules.

Our first aim is to use~(\ref{identify_F_link_cobordism_with_W(Flink)}) 
to distinguish the elements of~$C_{2q-1}(F_\mu)$. We define complete invariants
(and secondary invariants if certain primary invariants vanish) by
analyzing the Witt groups $W^{(-1)^{q+1}}(\Bl(\Q))$.
We proceed in three steps, explained in more detail in
Section~\ref{section:Blanchfield_form_invariants},
which run parallel to steps 2, 3 and 4 in chapter 2 of~\cite{She03mem}:
\begin{enumerate}
\item Obtain a direct sum decomposition of
$W^{(-1)^{q+1}}(\Bl(\Q))$ by ``devissage''.
One must prove that $\Bl(\Q)$ is an abelian category in which each module
has a finite composition series. 
\item Use hermitian Morita equivalence to show that each summand of
the group $W^{(-1)^{q+1}}(\Bl(\Q))$ is isomorphic to some group $W^1(E)$ where
$E$ is a division ring of finite dimension over $\Q$.
\item Recall from the literature invariants of each $W^1(E)$.
\end{enumerate}
In the knot theory case $\mu=1$ there is one summand of
$W^{(-1)^{q+1}}(\Bl(\Q))$ for each maximal ideal
$(p)\in\Q[z,z^{-1}]$ which is invariant under the involution $z\mapsto
z^{-1}$. The generator $p$ is often called an Alexander polynomial.
The division ring $E$ coincides with the quotient field
$\Q[z,z^{-1}]/(p)$ and $W^1(E)$ is the Witt group of hermitian forms
over $E$ (compare Milnor~\cite{Mil69}).

The following theorem and  
corollary are restated and proved in
Section~\ref{section:Blanchfield_form_invariants}; see
Theorem~\ref{BD_form_invariants_for_W_suffice:precise}
and Corollary~\ref{BD_form_invariants_for_Flinks_suffice:precise}.
\begin{theorem}\label{BD_form_invariants_for_W_suffice:vague}
The invariants (and secondary invariants) defined in
Section~\ref{section:Blanchfield_form_invariants} are sufficient to
distinguish the elements of the Witt groups $W^\pm(\Bl(\Q))$ of Blanchfield-Duval forms with coefficients in $\Q$.
\end{theorem}
\begin{corollary}\label{BD_form_invariants_for_Flinks_suffice:vague}
Let $q>1$ and suppose $\phi^0$ and $\phi^1$ are the
Blanchfield-Duval forms for the $(2q-1)$-dimensional $F_\mu$-links
$(L^0,\theta^0)$ and $(L^1,\theta^1)$ respectively. 
These two $F_\mu$-links are cobordant if and only if all the invariants
(and possible secondary invariants) of
\begin{equation*}
[\Q\otimes_\Z(\phi^0\oplus -\phi^1)] \in W^{(-1)^{q+1}}(\Bl(\Q))
\end{equation*}
defined in Section~\ref{section:Blanchfield_form_invariants} are trivial.
\end{corollary}
\noindent 
Corollary~\ref{BD_form_invariants_for_Flinks_suffice:vague} follows 
from~(\ref{identify_F_link_cobordism_with_W(Flink)}) and the fact that
the canonical map 
\begin{equation*}
W^{(-1)^{q+1}}(\Bl(\Z))\to
W^{(-1)^{q+1}}(\Bl(\Q))
\end{equation*}
is an injection, which we deduce from
Theorem~\ref{main_covering_construction_theorem} at the end of
Section~\ref{section:overview_of_BD_form_invariants}.
Corollary~\ref{BD_form_invariants_for_Flinks_suffice:vague} is also
a consequence of Theorem~\ref{theorem_relating_invariants:vague} and
Theorem~B of~\cite{She03mem}.

Our second aim is to understand the algebraic relationship between the Seifert
forms and the Blanchfield-Duval form of an $F_\mu$-link and prove that the
cobordism invariants defined in~\cite{She03mem} using Seifert forms are
equivalent to those defined in
Section~\ref{section:Blanchfield_form_invariants} via the  
Blanchfield-Duval form. Example~\ref{worked_example_Seifert} gives
a sample calculation of the Seifert form invariants in~\cite{She03mem}. 

In the knot theory case $\mu=1$, the relationship between Seifert and
Blanchfield forms has been investigated extensively by Kearton~\cite{Kea75},
Levine~\cite[\S14]{Lev77}, Farber~\cite[\S7.1]{Far83} and
Ranicki~(\cite[ch32]{Ran98},\cite{Ran03}). For
$\mu\geq1$ K.H.Ko~\cite{Ko89} used geometric arguments to obtain a formula 
for Cappell and Shaneson's homology surgery obstruction in terms of
the Seifert form. A formula for the
Blanchfield-Duval form in terms of the Seifert form, again based on geometric
arguments, can also be found in Cochran and
Orr~\cite[Thm4.2]{CocOrr94} in the slightly more 
general context of ``homology boundary links''. 

M.Farber related Seifert and Blanchfield-Duval forms of $F_\mu$-links
in a purely algebraic way~\cite{Far91,Far92}. Although the present
paper is logically independent of his work, we take up a number
of his ideas in Sections~\ref{section:seifert}
and~\ref{section:covering_construction}, providing
a systematic treatment in the language of 
hermitian categories. Whereas Farber takes coefficients in a field or
in $\Z$, in these sections we allow the coefficients to lie 
in an arbitrary associative ring $A$. 

The first step is to show that an $F_\mu$-link module admits a canonical
Seifert module structure (cf~\cite[p193]{Far91}). An $F_\mu$-link module
$M\in\Bl(A)$ is not in general finitely generated (or projective) as an
$A$-module so we introduce a larger category $\Sei_\infty(A)$
which contains $\Sei(A)$ as a full subcategory (see
Notation~\ref{Sei_notation}). We obtain a ``forgetful'' functor  
\begin{equation*}
U:\Bl(A)\to \Sei_\infty(A).
\end{equation*}
For example, in the case $\mu=1$ of knot theory, an object in
$\Bl(A)$ is a module $M$ over the ring
$A[z,z^{-1}]$ of Laurent polynomials with a presentation
\begin{equation*}
0\to (A[z,z^{-1}])^m\xrightarrow{\sigma} (A[z,z^{-1}])^m \to M\to 0
\end{equation*}
such that $1-z:M\to M$ is an isomorphism. The Seifert module $U(M)$ is the
$A$-module $M$ together with the endomorphism $(1-z)^{-1}$.

If $A=k$ is a field, Farber defined, for each $M\in\Bl(k)$,
the ``minimal lattice''~\cite[p194-199]{Far91} of $M$, a Seifert submodule
of $U(M)$ which is of finite $k$-dimension. We prefer to work directly
with $U(M)$ which is defined regardless of the coefficients and avoids
technicalities of Farber's definition. His
minimal lattice becomes isomorphic to $U(M)$ after one performs a
universal localization of categories which we describe a few
paragraphs below.

Given a Seifert surface for an $F_\mu$-link one can construct the free
cover by cutting the link complement along the Seifert surface and
gluing together infinitely many copies of the resulting manifold in
the pattern of the Cayley graph of $F_\mu$. 
Figure~1
illustrates the geometric construction in the case of
a $2$-component link.

The algebraic analogue of this geometric construction is a functor
\begin{equation*}
B:\Sei(A)\to\Bl(A)
\end{equation*}
from Seifert modules to $F_\mu$-link modules (see Definition~\ref{define_covering_construction}). Since
$U$ takes values in the larger category $\Sei_\infty(A)$ we
expand the domain of $B$ to $\Sei_\infty(A)$, by necessity replacing
$\Bl(A)$ by a larger category $\Bl_\infty(A)$. This process of
enlargement stops here for there are functors 
\begin{align*}
&U:\Bl_\infty(A)\to \Sei_\infty(A) \\
&B:\Sei_\infty(A)\to \Bl_\infty(A).
\end{align*}

We show in Section~\ref{section:adjunction} that $B$ is left adjoint to $U$.
Roughly speaking, this means that $B(V)$ is the ``free''
$F_\mu$-link module generated by the Seifert module $V$
(with respect to the functor $U$). In other words, $B$ is
universal (up to equivalence) among functors from Seifert modules to
$F_\mu$-link modules.

Returning our attention to the subcategories $\Sei(A)$ and $\Bl(A)$
whose definitions involve a ``finitely generated~projective''
condition we show that $B$ is compatible with the notions
of duality in $\Sei(A)$ and $\Bl(A)$, extending $B$ to a
``duality-preserving functor'' between ``hermitian categories'' 
\begin{equation}\label{introduce_duality_preserving_functor}
(B,\Phi,-1):\Sei(A)\to \Bl(A).
\end{equation}
(see definitions~\ref{define_hermitian_category}
and~\ref{define_dpres_funct} and
proposition~\ref{extend_B_to_duality_pres_functor}). 
The following theorem concerns the induced homomorphism of Witt groups:
\begin{equation}\label{Witt_map_induced_by_B}
B:W^{\pm}(\Sei(A))\to W^{\mp}(\Bl(A)).
\end{equation}
Recall that by Wedderburn's Theorem, a ring $A$ is semi-simple and Artinian
if and only if it is a product of matrix rings over division rings.
%
\tracingstats=2
\newcounter{rot}
\newcounter{rot2}
\newcounter{posx}
\newcounter{posy}
%
\newcounter{psz}
\newcounter{pcentx}
\newcounter{pcenty}
\newcounter{subpcentx}
\newcounter{subpcenty}
\newcounter{subsubpcentx}
\newcounter{subsubpcenty}
%

%
%
\newcounter{ovalthk}
\newcounter{ovalthn}
\newcounter{subovalthk}
\newcounter{subovalthn}
\newcounter{subsubovalthk}
\newcounter{subsubovalthn}
\setcounter{psz}{60}
\setcounter{subsubpcentx}{{\value{psz}*19/6}}
\setcounter{subsubpcenty}{0}
\setcounter{subpcentx}{{\value{psz}*21/10}}
\setcounter{subpcenty}{0}
\setcounter{pcentx}{0}
\setcounter{pcenty}{0}
\setcounter{subsubovalthk}{\value{psz}/8}
\setcounter{subsubovalthn}{\value{psz}/10}
\setcounter{subovalthk}{\value{psz}/6}
\setcounter{subovalthn}{\value{psz}/8}
\setcounter{ovalthk}{\value{psz}/4}
\setcounter{ovalthn}{\value{psz}/5}
\begin{equation*}\label{picture}
\Draw(0.8pt,0.8pt)
%
%
%
%
%
%
%
\def\smallestpic
{
\setcounter{posx}{\value{subsubpcentx}-\value{subsubovalthk}}
\setcounter{posy}{\value{subsubpcenty}+\value{psz}*9/20}
\MoveTo(\value{posx},\value{posy})
\MarkLoc(c+2-)
\setcounter{posx}{\value{subsubpcentx}-\value{subsubovalthk}}
\setcounter{posy}{\value{subsubpcenty}+\value{subsubovalthk}*5/4}
\MoveTo(\value{posx},\value{posy})
\MarkLoc(C-+)
\setcounter{posx}{\value{subpcentx}+\value{psz}*2/3}
\setcounter{posy}{\value{subpcenty}+\value{subovalthk}}
\MoveTo(\value{posx},\value{posy})
\MarkLoc(b+1+)
\Curve(b+1+,C-+,C-+,c+2-)
\setcounter{posx}{\value{subsubpcentx}-\value{subsubovalthk}}
\setcounter{posy}{\value{subsubpcenty}-\value{psz}*9/20}
\MoveTo(\value{posx},\value{posy})
\MarkLoc(c-2-)
\setcounter{posx}{\value{subsubpcentx}-\value{subsubovalthk}}
\setcounter{posy}{\value{subsubpcenty}-\value{subsubovalthk}*5/4}
\MoveTo(\value{posx},\value{posy})
\MarkLoc(C--)
\setcounter{posx}{\value{subpcentx}+\value{psz}*2/3}
\setcounter{posy}{\value{subpcenty}-\value{subovalthk}}
\MoveTo(\value{posx},\value{posy})
\MarkLoc(b+1-)
\Curve(b+1-,C--,C--,c-2-)
\setcounter{posx}{\value{subsubpcentx}+\value{subsubovalthk}}
\setcounter{posy}{\value{subsubpcenty}+\value{psz}*9/20}
\MoveTo(\value{posx},\value{posy})
\MarkLoc(c+2+)
\setcounter{posx}{\value{subsubpcentx}+\value{psz}*9/20}
\setcounter{posy}{\value{subsubpcenty}+\value{subsubovalthk}}
\MoveTo(\value{posx},\value{posy})
\MarkLoc(c+1+)
\setcounter{posx}{\value{subsubpcentx}+\value{subsubovalthk}}
\setcounter{posy}{\value{subsubpcenty}+\value{subsubovalthk}*5/4}
\MoveTo(\value{posx},\value{posy})
\MarkLoc(C++)
\Curve(c+2+,C++,C++,c+1+)
\setcounter{posx}{\value{subsubpcentx}+\value{subsubovalthk}}
\setcounter{posy}{\value{subsubpcenty}-\value{psz}*9/20}
\MoveTo(\value{posx},\value{posy})
\MarkLoc(c-2+)
\setcounter{posx}{\value{subsubpcentx}+\value{psz}*9/20}
\setcounter{posy}{\value{subsubpcenty}-\value{subsubovalthk}}
\MoveTo(\value{posx},\value{posy})
\MarkLoc(c+1-)
\setcounter{posx}{\value{subsubpcentx}+\value{subsubovalthk}}
\setcounter{posy}{\value{subsubpcenty}-\value{subsubovalthk}*5/4}
\MoveTo(\value{posx},\value{posy})
\MarkLoc(C+-)
\Curve(c-2+,C+-,C+-,c+1-)
%
%
\setcounter{posx}{\value{subsubpcentx}+\value{psz}*9/20}
\setcounter{posy}{\value{subsubpcenty}}
\MoveTo(\value{posx},\value{posy})
\DrawOval(\value{subsubovalthn},\value{subsubovalthk})
\setcounter{posx}{\value{subsubpcentx}}
\setcounter{posy}{\value{subsubpcenty}-\value{psz}*9/20}
\MoveTo(\value{posx},\value{posy})
\DrawOval(\value{subsubovalthk},\value{subsubovalthn})
\setcounter{posx}{\value{subsubpcentx}}
\setcounter{posy}{\value{subsubpcenty}+\value{psz}*9/20}
\MoveTo(\value{posx},\value{posy})
\DrawOval(\value{subsubovalthk},\value{subsubovalthn})
}
\def\quartpic{
%
\setcounter{posx}{\value{pcentx}-\value{psz}}
\setcounter{posy}{\value{pcenty}}
\MoveTo(\value{posx},\value{posy})
\DrawOval(\value{ovalthn},\value{ovalthk})
\setcounter{posx}{\value{pcentx}+\value{psz}}
\setcounter{posy}{\value{pcenty}+\value{ovalthk}}
\MoveTo(\value{posx},\value{posy})
\MarkLoc(a+1+)
\setcounter{posx}{\value{pcentx}+\value{ovalthk}}
\setcounter{posy}{\value{pcenty}+\value{psz}}
\MoveTo(\value{posx},\value{posy})
\MarkLoc(a+2+)
\setcounter{posx}{\value{pcentx}+\value{ovalthk}}
\setcounter{posy}{\value{pcenty}+\value{ovalthk}}
\MoveTo(\value{posx},\value{posy})
\MarkLoc(A++)
\Curve(a+1+,A++,A++,a+2+)
\setcounter{posx}{\value{pcentx}+\value{psz}}
\setcounter{posy}{\value{pcenty}-\value{ovalthk}}
\MoveTo(\value{posx},\value{posy})
\MarkLoc(a+1-)
\setcounter{posx}{\value{pcentx}+\value{ovalthk}}
\setcounter{posy}{\value{pcenty}-\value{psz}}
\MoveTo(\value{posx},\value{posy})
\MarkLoc(a-2+)
\setcounter{posx}{\value{pcentx}+\value{ovalthk}}
\setcounter{posy}{\value{pcenty}-\value{ovalthk}}
\MoveTo(\value{posx},\value{posy})
\MarkLoc(A+-)
\Curve(a+1-,A+-,A+-,a-2+)
%
%
%
%
\setcounter{posx}{\value{subpcentx}-\value{subovalthk}}
\setcounter{posy}{\value{subpcenty}+\value{psz}*2/3}
\MoveTo(\value{posx},\value{posy})
\MarkLoc(b+2-)
\setcounter{posx}{\value{subpcentx}-\value{subovalthk}}
\setcounter{posy}{\value{subpcenty}+\value{subovalthk}*5/4}
\MoveTo(\value{posx},\value{posy})
\MarkLoc(B-+)
\Curve(a+1+,B-+,B-+,b+2-)
\setcounter{posx}{\value{subpcentx}-\value{subovalthk}}
\setcounter{posy}{\value{subpcenty}-\value{psz}*2/3}
\MoveTo(\value{posx},\value{posy})
\MarkLoc(b-2-)
\setcounter{posx}{\value{subpcentx}-\value{subovalthk}}
\setcounter{posy}{\value{subpcenty}-\value{subovalthk}*5/4}
\MoveTo(\value{posx},\value{posy})
\MarkLoc(B--)
\Curve(a+1-,B--,B--,b-2-)
\setcounter{posx}{\value{subpcentx}+\value{subovalthk}}
\setcounter{posy}{\value{subpcenty}+\value{psz}*2/3}
\MoveTo(\value{posx},\value{posy})
\MarkLoc(b+2+)
\setcounter{posx}{\value{subpcentx}+\value{psz}*2/3}
\setcounter{posy}{\value{subpcenty}+\value{subovalthk}}
\MoveTo(\value{posx},\value{posy})
\MarkLoc(b+1+)
\setcounter{posx}{\value{subpcentx}+\value{subovalthk}}
\setcounter{posy}{\value{subpcenty}+\value{subovalthk}*5/4}
\MoveTo(\value{posx},\value{posy})
\MarkLoc(B++)
\Curve(b+2+,B++,B++,b+1+)
\setcounter{posx}{\value{subpcentx}+\value{subovalthk}}
\setcounter{posy}{\value{subpcenty}-\value{psz}*2/3}
\MoveTo(\value{posx},\value{posy})
\MarkLoc(b-2+)
\setcounter{posx}{\value{subpcentx}+\value{psz}*2/3}
\setcounter{posy}{\value{subpcenty}-\value{subovalthk}}
\MoveTo(\value{posx},\value{posy})
\MarkLoc(b+1-)
\setcounter{posx}{\value{subpcentx}+\value{subovalthk}}
\setcounter{posy}{\value{subpcenty}-\value{subovalthk}*5/4}
\MoveTo(\value{posx},\value{posy})
\MarkLoc(B+-)
\Curve(b-2+,B+-,B+-,b+1-)

%
%
\setcounter{posx}{\value{subpcentx}+\value{psz}*2/3}
\setcounter{posy}{\value{subpcenty}}
\MoveTo(\value{posx},\value{posy})
\DrawOval(\value{subovalthn},\value{subovalthk})
\setcounter{posx}{\value{subpcentx}}
\setcounter{posy}{\value{subpcenty}-\value{psz}*2/3}
\MoveTo(\value{posx},\value{posy})
\DrawOval(\value{subovalthk},\value{subovalthn})
\setcounter{posx}{\value{subpcentx}}
\setcounter{posy}{\value{subpcenty}+\value{psz}*2/3}
\MoveTo(\value{posx},\value{posy})
\DrawOval(\value{subovalthk},\value{subovalthn})
}
\setcounter{rot}{0}
\Do(1,4){
\MoveTo(\value{pcentx},\value{pcenty})
\setcounter{rot2}{\value{rot}+90}
\RotatedAxes(\value{rot},\value{rot2})
\quartpic
\smallestpic
\EndRotatedAxes
\addtocounter{rot}{90}
}
\MoveTo(\value{subpcentx},-\value{subpcentx})
\RotatedAxes(90,180)
\smallestpic
\EndRotatedAxes
\MoveTo(\value{subpcentx},\value{subpcentx})
\RotatedAxes(180,270)
\smallestpic
\EndRotatedAxes
\MoveTo(-\value{subpcentx},\value{subpcentx})
\RotatedAxes(270,0)
\smallestpic
\EndRotatedAxes
\MoveTo(-\value{subpcentx},-\value{subpcentx})
\RotatedAxes(0,90)
\smallestpic
\EndRotatedAxes
%
\MoveTo(-\value{subpcentx},\value{subpcentx})
\RotatedAxes(0,90)
\smallestpic
\EndRotatedAxes
\MoveTo(\value{subpcentx},\value{subpcentx})
\RotatedAxes(270,0)
\smallestpic
\EndRotatedAxes
\MoveTo(\value{subpcentx},-\value{subpcentx})
\RotatedAxes(180,270)
\smallestpic
\EndRotatedAxes
\MoveTo(-\value{subpcentx},-\value{subpcentx})
\RotatedAxes(90,180)
\smallestpic
\EndRotatedAxes
%
\MoveTo(0,-260)
\Text(--\mbox{Figure 1}--)
\EndDraw
\end{equation*}
\vspace*{2ex}
\begin{theorem}\label{main_covering_construction_theorem}
If $A$ is a semi-simple Artinian ring then (\ref{Witt_map_induced_by_B}) is an isomorphism.
\end{theorem}
\noindent The map (\ref{Witt_map_induced_by_B}) will be considered for
more general rings $A$ in subsequent work (joint with A.Ranicki). It
follows from the isomorphisms~(\ref{identify_F_link_cobordism_with_W(Seifert)})
and~(\ref{identify_F_link_cobordism_with_W(Flink)}) above that
(\ref{Witt_map_induced_by_B}) is an isomorphism when
$A=\Z$. 
\begin{theorem}\label{theorem_relating_invariants:vague}
The duality-preserving functor $(B,\Phi,-1)$ identifies the Seifert
form invariants of~\cite{She03mem} with the Blanchfield-Duval form invariants
of Section~\ref{section:Blanchfield_form_invariants}.
\end{theorem}
Theorem~\ref{theorem_relating_invariants:precise} 
below is a more precise statement of
Theorem~\ref{theorem_relating_invariants:vague}. The
invariants of~\cite{She03mem} are outlined in
Section~\ref{section:cobordism_invariants_seifert}.

Theorem~\ref{main_covering_construction_theorem} is proved in two
stages. The first stage is to establish that, for any ring $A$,
there is an equivalence between $(B,\Phi,-1)$ and 
a certain universal localization of hermitian categories. 
Taking up Farber's terminology we call a Seifert module
$V\in\Sei_\infty(A)$ primitive if $B(V)\cong 0$. We denote by
$\Prim_\infty(A)$ the category of primitive modules. One may write
\begin{equation*}
\Prim_\infty(A)=\Ker(\, B:\Sei_\infty(A) \to \Bl_\infty(A)\,).
\end{equation*} 
The category quotient
\begin{equation}\label{introduce_universal_cat_loc}
F:\Sei_\infty(A)\to \Sei_\infty(A)/\Prim_\infty(A)
\end{equation}
is universal among functors which make invertible morphisms whose kernel and
cokernel are primitive. In particular, primitive modules in
$\Sei_\infty(A)$ are made isomorphic to $0$ in
$\Sei_\infty(A)/\Prim_\infty(A)$. 

Since $B:\Sei_\infty(A)\to\Bl_\infty(A)$ is
left adjoint to $U$ it follows that 
$B$ exhibits the same universal property as $F$ although
only ``up to natural isomorphism''
(Proposition~\ref{B_is_universal_localization_up_to_equivalence}).
We conclude that $B$ is equivalent to $F$ and, with a little extra work,
establish that $(B,\Phi,-1):\Sei(A)\to\Bl(A)$ is equivalent
to a universal localization of hermitian categories
(Theorem~\ref{Bbar_is_equiv_on_fgproj_cats} and
Proposition~\ref{Bbar_is_eq_of_hermitian_cats}). 

In the knot theory case $\mu=1$ the category $\Sei_\infty(A)$
coincides with the category of (left) modules over the polynomial ring $A[s]$
in a central indeterminate $s$. Setting $t=s(1-s)$, the functor
(\ref{introduce_universal_cat_loc}) is the central localization
$A[s,t^{-1}]\otimes_{A[s]}\functor$ from the category of
$A[s]$-modules to the category of $A[s,t^{-1}]$-modules; see
Farber~\cite[Thm2.6]{Far83} and Ranicki~\cite{Ran03}.

Pere Ara recently gave an independent proof~\cite[Thm 6.2]{Ara03}
using Farber's minimal lattice that if $A=k$ is a field then
$\Bl(k)$ is equivalent to a localization of $\Sei(k)$ by a category of
primitive modules (for all $\mu\geq1$)\footnote{
His context differs
slightly in that the free algebra $k\langle X\rangle$ on a set
$X=\{x_1,\cdots,x_\mu\}$ takes the place of the group ring $k[F_\mu]$
in the present paper; the category denoted
$\cy{Z}$ in~\cite{Ara03} plays the role of $\Bl(k)$.
Consequently, there is only one kind of ``trivially primitive'' module
(denoted $M_0$ in~\cite{Ara03}) as compared with the two kinds
in~\cite{Far92} and Section~\ref{section:structure_of_primitives} below. 
Ara also related the modules in $\cy{Z}$ to modules over the Leavitt
algebra $L$. By definition, $L$ is the universal localization of
$k\langle X\rangle$ which makes invertible the map 
$(x_1\cdots x_\mu):k\langle X \rangle^{\mu} \to k\langle X\rangle$
(compare the Sato condition, Lemma~~\ref{Sato_condition} below). The
category $\cy{Z}$ turns out to be equivalent to the category of
finitely presented $L$-modules of finite length.}. \\[0.5ex]
\indent The second stage in the proof of
Theorem~\ref{main_covering_construction_theorem} involves the analysis
of primitive modules.
A Seifert module $V$ is called ``trivially primitive'' if the
endomorphism with which it is endowed is either zero or the identity.
In Proposition~\ref{identify_prim_infty} we show that every primitive
module $V\in\Prim_\infty(A)$ is composed of (possibly infinitely many)
trivially primitive modules.
Restricting attention to $\Sei(A)$ we show that if $A$ is semi-simple
and Artinian then every primitive module is composed of finitely many
trivially primitive modules in $\Sei(A)$ (cf Farber~\cite[\S3,\S7.10]{Far92}). 
The proof of Theorem~\ref{main_covering_construction_theorem} is completed in
Section~\ref{section:pf_of_Witt_isomorphism_thm} by establishing that
the Witt group of the subcategory $\Prim(A)\subset \Sei(A)$ of
primitive modules is trivial. See
Proposition~\ref{Witt_isomorphism_for_hermitian_quotient} 
and part 2.~of Lemma~\ref{Witt_group_of_primitives_is_zero}.

The definitions of Blanchfield-Duval form invariants in
Section~\ref{section:Blanchfield_form_invariants} parallel
the author's Seifert form invariants in~\cite[Ch2]{She03mem}. 
The three steps outlined above to analyze $W^\pm(\Bl(\Q))$ 
were applied to $W^\pm(\Sei(\Q))$ in~\cite{She03mem}.
Theorem~\ref{theorem_relating_invariants:vague} is proved in
Section~\ref{section:pf_of_eq_of_invariants} by checking that the
duality-preserving functor $(B,\Phi,-1)$ respects each of these three
steps.
\vspace*{1ex}

Let us summarize the contents of this paper.
Section~\ref{section:Blanchfield_form} discusses $F_\mu$-link modules
over an arbitrary ring $A$ and uses universal localization (cf
Vogel~\cite{Vog80,Vog82} and Duval~\cite{Duv86}) to describe hermitian
structure in the category $\Bl(A)$. We define the Witt groups
$W^{\pm}(\Bl(A))$ of Blanchfield-Duval forms.

In Section~\ref{section:Blanchfield_form_invariants} we set $A=\Q$ and 
define invariants of $W^{\pm}(\Bl(\Q))$, obtaining intrinsic cobordism
invariants for $F_\mu$-links. We discuss each of the steps
1-3.~listed above, reformulating and proving Theorem~\ref{BD_form_invariants_for_W_suffice:vague} and
Corollary~\ref{BD_form_invariants_for_Flinks_suffice:vague}.

In Section~\ref{section:seifert} we define Seifert modules and Seifert
forms with coefficients in an arbitrary ring $A$. We treat a worked
example of the invariants defined in~\cite{She03mem} and we define
a forgetful functor $U:\Bl_\infty(A)\to\Sei_\infty(A)$. 

Section~\ref{section:covering_construction} begins to study the functor
$B:\Sei_\infty(A)\to \Bl_\infty(A)$ from Seifert modules to $F_\mu$-link modules.
We prove that $B$ is left adjoint
to $U:\Bl_\infty(A)\to \Sei_\infty(A)$ and show that $B$ factors
through a category equivalence $\Sei_\infty(A)/\Prim_\infty(A) \to
\Bl_\infty(A)$. We describe the structure of the primitive modules -- those
which are sent to zero by $B$ -- and outline a construction of the localization
$\Sei_\infty(A)\to\Sei_\infty(A)/\Prim_\infty(A)$. We construct a 
duality-preserving functor
$(B,\Phi,-1):\Sei(A)\to\Bl(A)$ which is natural in $A$ and
factors through an equivalence $\Sei(A)/\Prim_\infty(A)\to\Bl(A)$ of
hermitian categories. If $A$ is a semi-simple Artinian ring we give a
simplified description of the primitive modules and the universal
localization of hermitian categories.

Section~\ref{section:Proof_that_B_induces_Witt_isomorphism} contains a
proof of Theorem~\ref{main_covering_construction_theorem} and 
a reformulation and proof of
Theorem~\ref{theorem_relating_invariants:vague}.
 
\paragraph{Acknowledgments:}
The invariants in Section~\ref{section:Blanchfield_form_invariants}
came into being during my thesis work at the University of Edinburgh
under the guidance of Andrew Ranicki. I am also indebted to
Andrew for several more recent conversations and e-mails
and for encouragement to complete this paper. I am grateful to
John Baez and James Dolan for helpful discussions in category theory, and to
Pere Ara and to the referee for their comments and corrections.

I thank the London Mathematical
Society for financial support to attend the ICMS workshop on
``Noncommutative Localization in Algebra and Topology'' in April 2002
and the Edinburgh Mathematical Society who
financed my visit to Edinburgh in August 2002.
\section{Blanchfield-Duval form}\label{section:Blanchfield_form}
\subsection{$F_\mu$-link Modules}
\label{section:blanchfield_modules}
Let $A$ be an associative ring with $1$. Modules will be left modules
except where otherwise stated. Let $A[F_\mu]$ denote
the group ring of the free group $F_\mu$; an element of $A[F_\mu]$ is
a formal sum of elements of $F_\mu$, with coefficients in $A$. Note
that elements of the group $F_\mu$ commute with elements of $A$ and
$A[F_\mu]\cong A\otimes_\Z \Z[F_\mu]$.

The symbols $\epsilon$ and $j$ will be used for three slightly distinct
purposes but the meaning will be clear from the context. 
Firstly, $j$ denotes the inclusion of $A$ in $A[F_\mu]$. Secondly $j$
denotes the functor $V\mapsto A[F_\mu]\otimes_A V$
from the category of $A$-modules to the category of
$A[F_\mu]$-modules. For brevity we write $V[F_\mu]$ in place of $A[F_\mu]\otimes_A
V$. Thirdly, we use $j$ to denote the inclusion of a module $V$ in
$V[F_\mu]$ given by
\begin{equation*}
V\cong A\otimes_A V\xrightarrow{j\otimes1} A[F_\mu]\otimes_A V = V[F_\mu].
\end{equation*}

In the opposite direction $\epsilon:A[F_\mu]\to A$ denotes the ring
morphism which sends every element of $F_\mu$ to $1\in A$ and is the
identity on $A$. We also write $\epsilon$ for the functor
$A\otimes_{A[F_\mu]}\functor$ from the category of $A[F_\mu]$-modules
to the category of $A$-modules. Thirdly, $\epsilon:V[F_\mu]\to V$ denotes
the morphism 
\begin{equation*}
V[F_\mu]=A[F_\mu]\otimes_A V \xrightarrow{\epsilon\otimes1} A\otimes_A V \cong V.
\end{equation*}
Note that the composite $\epsilon j$ of ring morphisms is the identity
$\id_A$ and the composite $\epsilon j$ of module  
morphisms is the identity on $V$. The composite
$\epsilon j$ of functors is naturally isomorphic to the identity functor on
$A$-modules and we sometimes suppress the natural isomorphism identifying
$A\otimes_{A[F_\mu]}(A[F_\mu]\otimes_A V)$ with $V$.
\begin{definition}\label{define_F_link_module}
An {\it $F_\mu$-link module} is an $A[F_\mu]$-module $M$ which lies in
an exact sequence:
\begin{equation}\label{Presentation_of_F-link_module}
0\to V[F_\mu]\xrightarrow{\sigma} V[F_\mu]\to M\to 0
\end{equation} 
such that $V$ is an $A$-module and $\epsilon(\sigma):V\to V$ is an isomorphism.
\end{definition}
As we remarked in the introduction, the examples of $F_\mu$-link modules 
in the literature are more often called ``link modules''.
Note that if $V$ is a finitely generated $A$-module then the $F_\mu$-link
module $M$ is finitely generated as an 
$A[F_\mu]$-module but usually not as an $A$-module
(see Lemma~\ref{Sato_condition} below).

It will be helpful to make the following observation about
the definition of $F_\mu$-link modules: The condition that
$\epsilon(\sigma)$ is an isomorphism implies that $\sigma$ is an
injection (see Lemma~\ref{sigma_is_injective} below).

\begin{notation}\label{Flk_notation}
Let $\Bl_\infty(A)$ denote the category of $F_\mu$-link modules and
$A[F_\mu]$-module homomorphisms. Thus $\Bl_\infty(A)$ is a full
subcategory of the category of $A[F_\mu]$-modules.

Let $\Bl(A)\subset \Bl_\infty(A)$ denote the category of
modules with a presentation~(\ref{Presentation_of_F-link_module})
such that $V$ is a finitely generated projective $A$-module and
$\epsilon(\sigma)$ is an isomorphism. The morphisms in $\Bl(A)$ are,
as usual,
the $A[F_\mu]$-module morphisms so $\Bl(A)$ is a full subcategory of
$\Bl_\infty(A)$.
\end{notation}
\noindent We show in
Lemma~\ref{Bl(functor)_is_functor} below that 
$\Bl_\infty(\functor)$ and $\Bl(\functor)$ are functorial in~$A$.
The following lemma gives an alternative characterization of
$F_\mu$-link modules. Let $z_1,\cdots,z_\mu$ denote generators for $F_\mu$.
\begin{lemma}(Sato~\cite{Sat81})\label{Sato_condition}
Suppose $M$ is an $A[F_\mu]$-module which has a presentation
\begin{equation}\label{presentation}
0\to V[F_\mu]\xrightarrow{\sigma} V'[F_\mu] \to M\to0.
\end{equation}
for some $A$-modules $V$ and $V'$.
The augmentation $\epsilon(\sigma):V\to V'$ is an isomorphism if and
only if the $A$-module homomorphism
\begin{equation}\label{Sato_map}
\begin{aligned}
\gamma:M^{\oplus \mu}
&\to M \\
(m_1,\cdots,m_\mu)&\mapsto \sum_{i=1}^\mu (1-z_i)m_i
\end{aligned}
\end{equation}
is an isomorphism.
\end{lemma}
\begin{proof}
There is an exact sequence
\begin{equation}\label{tree_presentation}
0\to (A[F_\mu])^{\oplus\mu}\xrightarrow{\gamma} A[F_\mu]
\xrightarrow{\epsilon} A\to0
\end{equation}
where $\gamma(l_1,\cdots,l_\mu)=\sum_{i=1}^\mu (1-z_i)l_i$ for all
$l_1,\cdots,l_\mu\in A[F_\mu]$. Now (\ref{tree_presentation}) is split
(by $j$) when regarded as a sequence of right $A$-modules so the functors
$\functor\otimes_A V$ and $\functor\otimes_A V'$ lead to a commutative diagram
\begin{equation*}
\xymatrix@R=3ex{
  & 0 \ar[d]        & 0      \ar[d] & 0\ar[d] & \\
0\ar[r] & (V[F_\mu])^{\oplus\mu}
  \ar[r]^{\sigma^{\oplus\mu}}\ar[d]_{1\otimes\gamma} &
  ({V'}[F_\mu])^{\oplus\mu}
  \ar[r]\ar[d]_{1\otimes\gamma} & M^{\oplus\mu}\ar[d]^\gamma \ar[r] & 0 \\ 
0\ar[r] & V[F_\mu] \ar[r]^\sigma\ar[d]_\epsilon & V'[F_\mu] \ar[r] \ar[d]_\epsilon & M\ar[d]\ar[r] & 0 \\
0\ar[r] & V \ar[r]^{\epsilon(\sigma)}\ar[d] & V' \ar[r]\ar[d] & 0 \\
 & 0 & 0 &}
\end{equation*}
in which the first two rows and the first two columns are exact. A standard
diagram chase (e.g.~\cite[p49]{Mac95}) shows that the third row is exact if
and only if the third column is exact and the Lemma follows.
\end{proof}

We discuss next completions of $A[F_\mu]$-modules, which we shall need
both to refine Definition~\ref{define_F_link_module} (in
Lemma~\ref{sigma_is_injective}) and later to study the universal
augmentation localization of $A[F_\mu]$ (see
Lemma~\ref{i_Sigma_injective_for_AF}).

Let $I=\Ker(\epsilon:A[F_\mu]\to A)$. If $N$ is an $A[F_\mu]$-module one
defines
\begin{equation*}
\widehat{N}= \varprojlim_n \frac{N}{I^nN}.
\end{equation*}
An $A[F_\mu]$-module morphism $N\to N'$ maps $I^nN$ to $I^nN'$ for each
$n$ and therefore induces a homomorphism $\widehat{N}\to\widehat{N'}$.
\begin{caveat}
The natural isomorphisms $(A[F_\mu]/I^n)\otimes_{A[F_\mu]} N \to
N/(I^nN)$ induce a map $\widehat{A[F_\mu]}\otimes N  \to \widehat{N}$
but the latter is not in general an isomorphism.
\end{caveat}
In the examples with which we are most concerned,
$N=V[F_\mu]=A[F_\mu]\otimes_A V$ for some $A$-module $V$. One can
describe $\widehat{V[F_\mu]}$ as a module of power
series as follows. Let $X=\{x_1,\cdots,x_\mu\}$ and let $\Z\langle
X\rangle$ denote the free ring generated by $X$ (in other words the ring of
``polynomials'' in non-commuting indeterminates
$x_1,\cdots,x_\mu$). Let $A\langle X\rangle= A\otimes_\Z \Z\langle
X\rangle$ so that the elements of $A$ are formal sums of words in the
alphabet $X$ with coefficients in $A$. Let $V\langle X\rangle$ be the
$A\langle X\rangle$-module $A\langle X\rangle \otimes_A V$ and denote
by $X^n V\langle X\rangle \subset V\langle X\rangle$ the submodule
whose elements are formal sums of words of length at least $n$ with
coefficients in $V$. We may now define the $X$-adic completion
\begin{equation*}
V\langle\langle X\rangle\rangle=\varprojlim_n \frac{V\langle
X\rangle}{X^nV\langle X\rangle}. 
\end{equation*}
in which an element is a formal power
series in non-commuting indeterminates $x_1,\cdots,x_\mu$ with
coefficients in $V$.

\begin{lemma}\label{identify_free_completions}
There is a natural isomorphism $V\langle\langle X\rangle\rangle \cong
\widehat{V[F_\mu]}$.
\end{lemma}
\begin{proof}
The required isomorphism is induced by the isomorphisms
\begin{align*}
\frac{V[F_\mu]}{I^nV[F_\mu]} &\cong \frac{V\langle X\rangle}{X^nV\langle X\rangle} \\
z_i &\mapsto 1+x_i \\
z_i^{-1} &\mapsto 1-x_i + x_i^2 \cdots + (-1)^{n-1}x_i^{n-1}\qedhere
\end{align*}
\end{proof}
\begin{lemma}\label{completion_is_injection}
If $V$ is an $A$-module, the canonical map $V[F_\mu]\to \widehat{V[F_\mu]}$ is
an injection.
\end{lemma}
\begin{proof}
The argument of Fox~\cite[Corollary 4.4]{Fox53} implies that
$\displaystyle{\bigcap_{n=0}^\infty I^nV[F_\mu]}=0$ and the Lemma follows.
\end{proof}
\begin{lemma}\label{one_plus_tau_invertible}
If $V$ is an $A$-module and $\tau:V[F_\mu]\to V[F_\mu]$ is an
$A[F_\mu]$-module morphism such that $\epsilon(\tau)=0:V\to V$ then
$\tau(V[F_\mu])\subset IV[F_\mu]$ and the map
$1+\tau:\widehat{V[F_\mu]}\to \widehat{V[F_\mu]}$ is invertible. 
\end{lemma}
\begin{proof}
The commutative diagram
\begin{equation*}
\xymatrix@R=3ex{V[F_\mu] \ar[d]_{\epsilon} \ar[r]^{\tau} & V[F_\mu] \ar[d]^{\epsilon} \\
V \ar[r]_{\epsilon(\tau)} & V
}
\end{equation*}
implies that if $\epsilon(\tau)=0$ then 
\begin{equation*}
\tau(V[F_\mu])\subseteq \Ker(\epsilon:V[F_\mu]\to V) =
\Ker(\epsilon\otimes 1:A[F_\mu]\otimes_A V \to A\otimes_A V)
\end{equation*}
Since the surjection $\epsilon:A[F_\mu]\to A$ is split by $j$ we
obtain
\begin{equation*}
\tau(V[F_\mu])\subseteq I\otimes_A V = I(A[F_\mu]\otimes_A V)=IV[F_\mu].
\end{equation*}
By induction, $\tau^n(V[F_\mu])\subset I^nV[F_\mu]$ for all
$n$ so $1+\tau:\widehat{V[F_\mu]}\to \widehat{V[F_\mu]}$ has inverse
\begin{equation*}
(1+\tau)^{-1}=1-\tau+\tau^2-\tau^3 +\cdots~.\qedhere
\end{equation*}
\end{proof}
\begin{lemma}\label{sigma_is_injective}
If $V$ and $V'$ are $A$-modules and $\sigma:V[F_\mu]\to V'[F_\mu]$ is an $A[F_\mu]$-module
homomorphism such that $\epsilon(\sigma):V\to V'$ is an isomorphism
then $\sigma$ is an injection and the induced map
$\sigma:\widehat{V[F_\mu]}\to\widehat{V'[F_\mu]}$ is an isomorphism.
\end{lemma}
\begin{proof}
Let $\sigma'=(j\epsilon(\sigma))^{-1}\sigma:V[F_\mu]\to V[F_\mu]$. Now
$\epsilon(\sigma')=1_{V[F_\mu]}$ so we may write $\sigma'=1+\tau$ where
$\epsilon(\tau)=0$. Now $\sigma':\widehat{V[F_\mu]} \to
\widehat{V[F_\mu]}$ is an isomorphism by
Lemma~\ref{one_plus_tau_invertible} so
$\sigma:\widehat{V[F_\mu]}\to\widehat{V'[F_\mu]}$ is an isomorphism. 
The commutative diagram
\begin{equation*}
\xymatrix@=3ex{V[F_\mu] \ar@{>->}[d] \ar[r]^{\sigma'} & V[F_\mu]\ar@{>->}[d] \\
{\widehat{V[F_\mu]}} \ar[r]_{\sigma'} & {\widehat{V[F_\mu]}}}
\end{equation*}
implies that $\sigma':V[F_\mu]\to V[F_\mu]$ and $\sigma:V[F_\mu]\to
V'[F_\mu]$ are injections. 
\end{proof}
\begin{lemma}\label{Bl(functor)_is_functor}
A homomorphism $A\to A'$ of rings induces functors
\begin{align*}
A'[F_\mu]\otimes_{A[F_\mu]}\functor 
&:\Bl_\infty(A)\to\Bl_\infty(A') \quad \mbox{and}\\
A'[F_\mu]\otimes_{A[F_\mu]}\functor &:\Bl(A)\to\Bl(A')
\end{align*}
\end{lemma}
\begin{proof}
If $M\in\Bl_\infty(A)$ then $M$ has a
presentation~(\ref{Presentation_of_F-link_module}) such that
$\epsilon(\sigma)$ is invertible. Applying
$A'[F_\mu]\otimes_{A[F_\mu]}\functor$ one obtains an exact sequence
\begin{equation*}
A'[F_\mu]\otimes_{A[F_\mu]}V[F_\mu]\xrightarrow{1\otimes\sigma}
A'[F_\mu]\otimes_{A[F_\mu]}V[F_\mu] \to
A'[F_\mu]\otimes_{A[F_\mu]} M \to 0.
\end{equation*}
The naturality of the identifications
$A'\otimes_A(A\otimes_{A[F_\mu]} V[F_\mu]) \cong A'\otimes_A V \cong
A'\otimes_{A'[F_\mu]}(A'[F_\mu]\otimes_{A[F_\mu]} V[F_\mu])$
implies that
$\epsilon(1\otimes\sigma)=1\otimes\epsilon(\sigma):A'\otimes_AV\to
A'\otimes_AV$. But $1\otimes\epsilon(\sigma)$ is an isomorphism so
$1\otimes\sigma$ is injective by Lemma~\ref{sigma_is_injective}. Thus
$A'[F_\mu]\otimes_{A[F_\mu]} M\in\Bl_\infty(A')$.
The argument for $\Bl(\functor)$ is
similar, for if $V$ is a finitely generated projective
$A$-module then $A'\otimes_A V$ is a finitely generated projective $A'$-module.
\end{proof}

\subsection{Hermitian Categories}
Recall that an involution on a ring $A$ is a map $A\to A; a\mapsto \overline{a}$ such that
$\overline{\overline{a}}=a$, $\overline{a+b}=\overline{a}+\overline{b}$ and
$\overline{(ab)}=\overline{b}\overline{a}$ for all $a,b\in A$.
If $A$ is a ring with involution then the category $\Bl(A)$ can be
endowed with a notion of duality. But let us begin with simpler examples:
\begin{example}\label{duality_functor_on_A_Proj}
Suppose $A$ is a ring with an involution.
Let $A\proj$ denote the category of finitely generated projective (left)
$A$-modules. There is a duality functor defined on modules by
$V\mapsto V^*=\Hom(V,A)$ for $V\in A\proj$ and on morphisms by
$f\mapsto f^*=\functor\circ f :(V')^*\to V^*$ for
$f\in\Hom_A(V,V')$. In short, 
\begin{equation*}
\functor^*=\Hom_A(\functor,A):A\proj\to
A\proj.
\end{equation*}
Note that $\Hom_A(V,A)$ is a left $A$-module with
\begin{equation*}
(a.\xi)(x) = \xi(x)\overline{a}
\end{equation*}
for all $a\in A$, $\xi\in\Hom_A(V,A)$ and $x\in V$.
\end{example}
\begin{example}\label{finite_abelian_group_duality}
The category of finite abelian groups admits the duality functor
\begin{equation}\label{duality_functor_on_finabgp}
\functor^\wedge=\Ext_\Z(\functor,\Z).
\end{equation} 
\end{example}
A finite abelian group $M$ bears similarity to an $F$-link
module in that there exists a presentation
\begin{equation}\label{present_fin_ab_gp}
0\to \Z^n \xrightarrow{\sigma} \Z^n \to M \to 0
\end{equation}
and $1\otimes\sigma: \Q\otimes_\Z \Z^n \to \Q\otimes_\Z \Z^n$ is an
isomorphism. A more explicit description of the duality functor in
Example~\ref{finite_abelian_group_duality} is the following: 
\begin{lemma}\label{make_finabgp_duality_explicit}
There is a natural isomorphism 
\begin{equation*}
\Ext_\Z(\functor,\Z)\cong \Hom_\Z\left(\functor,\frac{\Q}{\Z}\right).
\end{equation*}
\end{lemma}
\begin{proof}
Suppose $M$ is a finite abelian group. The short exact sequence 
\begin{equation*}
0\to \Z\to \Q\to \Q/\Z\to 0
\end{equation*}
gives rise to a long exact sequence
\begin{equation*}
0=\Hom(M,\Q)\to \Hom(M,\Q/\Z)\to \Ext(M,\Z)\to \Ext(M,\Q)\to\cdots
\end{equation*}
which is natural in $M$. The presentation~(\ref{present_fin_ab_gp})
implies that 
\begin{equation*}
\Ext(M,\Q)=\Coker(\Hom(\Z^n,\Q)\xrightarrow{\sigma^*}\Hom(\Z^n,\Q))=0
\end{equation*}
so $\Hom(M,\Q/\Z)\to \Ext(M,\Z)$ is an isomorphism.
\end{proof}

The following general definition subsumes
Examples~\ref{duality_functor_on_A_Proj}
and~\ref{finite_abelian_group_duality} and the
category $\Bl(A)$ which we wish to study:
\begin{definition}\label{define_hermitian_category}
A hermitian category is a triple $(\cy{C},\functor^*,i)$ where
\begin{itemize}
\item $\cy{C}$ is an additive category,
\item $\functor^*:\cy{C}\rightarrow \cy{C}$ is an (additive) contravariant
functor and
\item $(i_V)_{V\in\cy{C}}:id\rightarrow (\functor^*)^*=\functor^{**}$ is a natural isomorphism such that 
$i^*_Vi_{V^*}=id_{V^*}$ for all $V\in\cy{C}$. 
\end{itemize}
\end{definition}
The functor $\functor^*$ is called a duality functor.
We usually abbreviate $(\cy{C},\functor^*,i)$ to $\cy{C}$ and identify $V$
with~$V^{**}$ via $i_V$. It follows from
Definition~\ref{define_hermitian_category} that if $\cy{C}$ is an abelian
hermitian category then $\functor^*$ is an equivalence of categories
and hence respects exact sequences.

If $A$ is a ring with involution then there is a unique involution on
$A[F_\mu]$ such that $\overline{g}=g^{-1}$ for each $g\in F_\mu$ and
such that the inclusion of $A$ in $A[F_\mu]$ respects the
involutions. The category of finitely generated projective
$A[F_\mu]$-modules therefore admits a duality functor as in
Example~\ref{duality_functor_on_A_Proj}.

Returning to $\Bl(A)$, duality is defined in a manner analogous to
Example~\ref{finite_abelian_group_duality}.
\begin{definition}\label{define_duality_on_Flink_modules}
Define $\functor^\wedge = \Ext_{A[F_\mu]}(\functor,A[F_\mu]):
\Bl(A)\to \Bl(A)$. 
\end{definition}
Note that if $M$ has
presentation~(\ref{Presentation_of_F-link_module}) then
$M^\wedge=\Ext_{A[F_\mu]}(M,A[F_\mu])$ has presentation
\begin{equation*}
0\to (V[F_\mu])^*\xrightarrow{\sigma^*} (V[F_\mu])^*\to M^\wedge\to 0
\end{equation*}
where $(V[F_\mu])^*=\Hom_{A[F_\mu]}(V[F_\mu],A[F_\mu])\cong \Hom_A(V,A)[F_\mu]= V^*[F_\mu]$;
see Lemma~\ref{lemma:duality} and
Remark~\ref{explicit_formula_for_Mwedge_pres} below.
The natural isomorphism 
\begin{equation*}
V[F_\mu]\to (V[F_\mu])^{**}
\end{equation*}
induces a natural isomorphism
$i_M:M\to M^{\wedge\wedge}$ with
$i_M^\wedge i_{M^\wedge}=\id_{M^\wedge}$.

There is also a more explicit description of $M^\wedge$ which is
analogous to Lemma~\ref{make_finabgp_duality_explicit}; see
Lemma~\ref{lemma:duality} below. Unlike $\Z$, the ring $A[F_\mu]$ is
in general highly non-commutative so universal localization will be
required. 

\subsection{Universal
Localization}\label{section:universal_localization_of_rings}
Let $R$ be a ring (associative with unit) and let $\Sigma$ be a set of
 (isomorphism classes of) triples $(P_1\ ,\ P_0\ ,\ \sigma:P_1\to
 P_0)$ where $P_0$ and $P_1$ are finitely generated projective $R$-modules.
In our application $R$ will be $A[F_\mu]$ and $\Sigma$ will
contain the endomorphisms $\sigma:V[F_\mu]\to V[F_\mu]$ such that $V$ is finitely
generated and projective as an $A$-module and $\epsilon(\sigma)$ is
 an automorphism of $V$. 

A homomorphism $\nu:R\rightarrow S$ is said to be {\it
$\Sigma$-inverting} if 
\begin{equation*}
1\otimes\sigma:S\otimes P_1\rightarrow S\otimes P_0
\end{equation*}
is invertible for each morphism $\sigma\in\Sigma$.  
There exists a universal $\Sigma$-inverting
homomorphism\footnote{The ring $R_\Sigma$
was denoted $\Sigma^{-1}R$ in~\cite{She03mem}} which, for consistency
with the other papers in the volume, will be 
denoted $i_\Sigma:R\rightarrow R_\Sigma$. The universal property is that
every $\Sigma$-inverting homomorphism $\nu:R\to S$ may be written uniquely as a
composite 
\begin{equation*}
R\xrightarrow{i_\Sigma}R_\Sigma\xrightarrow{\overline \nu}
S;
\end{equation*}
If $R$ is commutative and each $\sigma\in \Sigma$ is an endomorphism
then the localization is the ring of fractions 
\begin{equation*}
R_\Sigma= R_S = \{p/q\ \mid\ p\in R, q\in S\}
\end{equation*}
whose denominators lie in the multiplicative set $S\subseteq R$
 generated by the determinants of the morphisms in $\Sigma$:
\begin{equation*}
S= \left\{\prod_{i=1}^r \det(\sigma_i) \mid r\in\Z, r\geq0,
\sigma_i\in\Sigma\right\}
\end{equation*}
More general constructions of $i_\Sigma$ may be found
in~\cite[Ch4]{Scho85}, \cite[p255]{Coh71} or~\cite{Coh03}.

If $R=A[F_\mu]$ and $\Sigma$ is defined as above, the inclusion of
$A[F_\mu]$ in $\widehat{A[F_\mu]}$ is $\Sigma$-inverting; see
Lemmas~\ref{completion_is_injection} and~\ref{sigma_is_injective}. 
By the universal property of $A[F_\mu]_\Sigma$ there is therefore a
commutative diagram
\begin{equation}\label{completion_and_localization}
\xymatrix{
A[F_\mu]\ar@/_1.5pc/[rr] \ar[r]^{i_\Sigma} & A[F_\mu]_\Sigma \ar[r]^(0.35)\gamma &
\widehat{A[F_\mu]}\cong A\langle\langle X\rangle\rangle
}.
\end{equation}
where the natural isomorphism $\widehat{A[F_\mu]}\cong A\langle\langle
X\rangle\rangle$ is defined as in Lemma~\ref{identify_free_completions}.

The image of $\gamma$ is the ring $A_\rat\langle\langle
X\rangle\rangle$ of rational power series (see~\cite[\S4]{She01}). If $A$
is a field or a principal ideal domain then $\gamma$ is known to be
injective so $A[F_\mu]_\Sigma$ can be identified with $A_\rat\langle\langle
X\rangle\rangle$
(Cohn and Dicks~\cite[p416]{CohDic76}, Dicks and Sontag~\cite[Thm 24]{DicSon78}, Farber and
Vogel~\cite{FarVog92}). In the knot theory case $\mu=1$, if $A$ is
commutative then the localization is a ring of fractions 
\begin{equation*}
A[F_1]_{\Sigma} ~\cong~\{p/q\ \mid\ p,q\in A[z,z^{-1}], q(1)\
\mbox{is invertible}\}~\cong~A_\rat[[x]].
\end{equation*}
and $\gamma$ is injective. However, there exist non-commutative rings $A$ such
that $\gamma$ is not injective~\cite[Prop 1.2]{She01}.

Diagram~(\ref{completion_and_localization}) and
Lemma~\ref{completion_is_injection} imply:
\begin{lemma}\label{i_Sigma_injective_for_AF}
The localization $i_\Sigma:A[F_\mu]\to A[F_\mu]_\Sigma$ is
injective.\hspace*{\fill}\qed
\end{lemma}
The following is a generalization of
Lemma~\ref{sigma_is_injective}:
\begin{lemma}\label{sigma_injective}
If $i_\Sigma:R\to R_\Sigma$ is injective then each
$\sigma\in\Sigma$ is injective. 
\end{lemma}
\begin{proof}
Suppose $\sigma:P_1\to P_0$ is in $\Sigma$. There is a commutative diagram
\begin{equation*}
\xymatrix{
R\otimes_R P_1 \ar[r]^{i_\Sigma\otimes 1} \ar[d]_{1\otimes\sigma}& R_\Sigma\otimes_R P_1 \ar[d]^{1\otimes\sigma}\\
R\otimes_R P_0 \ar[r]_{i_\Sigma\otimes 1} & R_\Sigma\otimes_R P_0}
\end{equation*}
and $1\otimes\sigma:R_\Sigma\otimes_R P_1\to R_\Sigma\otimes_R P_0$ is
an isomorphism. If $i_\Sigma$ is injective then 
$1\otimes\sigma:R\otimes_R P_1\to R\otimes_R P_0$ is also injective so
$\sigma$ is injective.
\end{proof}
%
%
%
\begin{lemma}
\label{lemma:duality}
Suppose $i_\Sigma:R\rightarrow
R_\Sigma$ is an injection and $M=\Coker(\sigma)$ with $\sigma\in\Sigma$.
\begin{enumerate}
\item The (right) $R$-module $M^\wedge=\Ext_R(M,R)$ is isomorphic to
$\Coker(\sigma^*)$.
\item There is a natural isomorphism $\Ext_R(M,R)\cong \Hom_R(M,R_\Sigma/R)$. 
\end{enumerate}
\end{lemma}
\noindent If $R$ has involution one can regard these right modules as left
modules.
\begin{proof}[Proof of Lemma~\ref{lemma:duality}] (Compare Example~\ref{finite_abelian_group_duality}). \\
\noindent 
1. By Lemma~\ref{sigma_injective}, the map $\sigma$ is injective so
   $M$ has presentation
\begin{equation}\label{present_M}
0\to P_1\xrightarrow{\sigma}P_0\to M\to 0.
\end{equation}
There is therefore an exact sequence
\begin{equation*}
P_0^*\xrightarrow{\sigma^*} P_1^* 
\to \Ext_R(M,R)\to \Ext_R(P_0,R)=0 
\end{equation*}
\smallskip \\ \noindent
2.~The short exact sequence of $(R,R)$-bimodules
\begin{equation*}
0\rightarrow R\rightarrow R_\Sigma\rightarrow
R_\Sigma/R\rightarrow 0
\end{equation*}
induces a long exact sequence of right $R$-modules
\begin{multline*}
\cdots \to \Hom_R(M,R_\Sigma)\to\Hom_R(M,{R_\Sigma}/R)  \\ 
\to \Ext_R(M,R) \to \Ext_R(M,R_\Sigma) \to \cdots.
\end{multline*}
which is natural in $M$.
It remains to prove that $\Hom_R(M,R_\Sigma)=\Ext_R(M,R_\Sigma)=0$.
The presentation~(\ref{present_M}) gives rise to the long exact sequence
\begin{multline*}
0 \to \Hom_R(M,R_\Sigma)\to \Hom_R(P_0,R_\Sigma)
\xrightarrow{\sigma^*} \Hom_R(P_1,R_\Sigma) \\
\to \Ext_R(M,R_\Sigma) \to 0.
\end{multline*}
There is a natural isomorphism 
\begin{equation*}
\Hom_{R_\Sigma}(R_\Sigma\otimes_R \functor\,,\, R_\Sigma)\to
\Hom_R(R\otimes_R\functor\,,\,R_\Sigma)
\end{equation*}
induced by $i_\Sigma:R\to R_\Sigma$ and, in particular, a commutative diagram
\begin{equation}\label{hom(functor,R_Sigma)}
\begin{gathered}
\xymatrix{
{\Hom_{R_\Sigma}(R_\Sigma\otimes_R P_0,R_\Sigma)}\ar[r]^{(\id\otimes\sigma)^*}\ar[d]_\cong &
{\Hom_{R_\Sigma}(R_\Sigma\otimes_R P_1,R_\Sigma)} \ar[d]^\cong \\
{\Hom_R(P_0,R_\Sigma)}\ar[r]_{\sigma^*} & {\Hom_R(P_1,R_\Sigma)}
}
\end{gathered}
\end{equation}
The upper horizontal arrow is an isomorphism since $\id\otimes\sigma$
is an isomorphism so the lower horizontal arrow is also an isomorphism.
Thus
\begin{equation*}
\Hom_R(M,R_\Sigma)=\Ext_R(M,R_\Sigma)=0
\end{equation*}
as required.
\end{proof}
\begin{remark}\label{explicit_formula_for_Mwedge_pres}
If $R\to R_\Sigma$ is injective and $\sigma\in\Sigma$ then 
the sequence $0\to P_1\xrightarrow{\sigma} P_0 \xrightarrow{q} M\to 0$
is exact by Lemma~\ref{sigma_injective}. Now Lemma~\ref{lemma:duality}
gives an exact sequence 
\begin{equation*}
P_0^*\xrightarrow{\sigma^*} P_1^* \xrightarrow{q'}
\Hom_R(M,R_\Sigma/R)\to 0.
\end{equation*}
Let us give an explicit formula for $q'$. There is a short exact sequence of
right $R$-module chain complexes
\begin{equation*}
\xymatrix{
0\ar[r] & \Hom_R(P_0,R) \ar[r]\ar[d]_{\sigma^*} & \Hom_R(P_0,R_\Sigma)\ar[r]\ar[d]_{\sigma^*}^\cong
&  \Hom_R(P_0,R_\Sigma/R) \ar[r]\ar[d]_{\sigma^*} & 0 \\ 
0\ar[r] & \Hom_R(P_1,R) \ar[r] & \Hom_R(P_1,R_\Sigma)\ar[r] &
\Hom_R(P_1,R_\Sigma/R)\ar[r] & 0
}
\end{equation*}
and the natural isomorphism
\begin{equation*}
\Ext_R(M,R)=\frac{\Hom_R(P_1,R)}{\Image(\sigma*)} \longrightarrow \Hom_R(M,R_\Sigma/R)
\end{equation*}
is inverse to the boundary map in the induced long exact homology
sequence from the kernel of the right-most
$\sigma^*$ to the cokernel of the left-most $\sigma^*$. Written out at length,
the map $q':P_1^*\to \Hom_R(M,R_\Sigma/R)$ is the composite
\begin{multline*}
\Hom_R(P_1,R)\to\Hom_{R_\Sigma}(R_\Sigma\otimes P_1,R_\Sigma)
 \xleftarrow[(\id\otimes\sigma)^*]{\cong} 
\Hom_{R_\Sigma}(R_\Sigma\otimes P_0,R_\Sigma) \\
\xrightarrow{\cong}
\Hom_R(P_0,R_\Sigma)\to\Hom_R(P_0,R_\Sigma/R),
\end{multline*}
the image of which lies in the submodule $\Hom_R(M,R_\Sigma/R)\subset
\Hom(P_0,R_\Sigma/R)$.   
Suppose $f\in P_1^*=\Hom_R(P_1,R)$ and $m\in M$. Choose $x\in P_0$
such that $q(x)=m$ and write $1\otimes x\in R_\Sigma\otimes P_0$.
Now
\begin{equation*}
q'(f)(m)=(\id\otimes f)\left((\id\otimes\sigma)^{-1}(1\otimes x)\right)
\in R_\Sigma/R.
\end{equation*}
\end{remark}
Combining Lemmas~\ref{lemma:duality} and~\ref{i_Sigma_injective_for_AF} we have
\begin{proposition}\label{identify_Mwedge_hom_to_loc}
There is a natural isomorphism of contravariant functors
\begin{equation*}
\Ext_{A[F_\mu]}(\functor,A[F_\mu])\cong\Hom_{A[F_\mu]}(\functor,A[F_\mu]_\Sigma/A[F_\mu]):\Bl(A)\to\Bl(A).
\qedhere
\end{equation*}
\end{proposition}
\subsection{Hermitian forms and the Witt group}\label{section:define_hermitian_forms_and_witt_gp}
As we noted in the introduction (equations~(\ref{identify_F_link_cobordism_with_W(Seifert)})
and~(\ref{identify_F_link_cobordism_with_W(Flink)})), the cobordism
group $C_n(F_\mu)$ can be identified with a Witt group of Seifert or
Blanchfield-Duval forms. Let us recall the 
definition of a hermitian form in a hermitian category and the
appropriate definition of Witt group.
\begin{definition}
Let $\zeta=1$ or $-1$. A {\it $\zeta$-hermitian form} in a hermitian
category $(\cy{C},*,i)$ is a pair $(V,\phi)$ where $\phi:V\to V^*$ and
$\phi^*i_V=\zeta \phi$. If $\phi$ is an isomorphism then $\phi$ is
called {\it non-singular}.
\end{definition}
For example, in the category $\Bl(A)$ a $\zeta$-hermitian form is
an $A[F_\mu]$-module isomorphism $\phi:M\to M^\wedge=\Hom(M,A[F_\mu]_\Sigma/A[F_\mu])$  
such that $\phi^\wedge=\zeta\phi$ (we suppress the natural
isomorphism $i_M$).
\begin{definition} An object $V$ in a hermitian category
$(\cy{C},*,i)$ is called {\it self-dual} if $V$ is isomorphic to
$V^*$. If there exists a non-singular $\zeta$-hermitian form
$(V,\phi)$ then $V$ is called {\it $\zeta$-self-dual}.
\end{definition}

When one has a suitable notion of exact sequences in a hermitian category one
can define the Witt group of the category. For simplicity, suppose $\cy{C}$ 
is a full subcategory of an abelian category $\cy{A}$, so that every
morphism in $\cy{C}$ has kernel, image and cokernel in
$\cy{A}$. Suppose further that $\cy{C}$ is admissible in $\cy{A}$ in
the following sense: If there is an exact sequence $0\to V\to V'\to
V''\to 0$ and the modules $V$ and $V''$ lie in $\cy{C}$ then $V'$ lies
in $\cy{C}$. In Section~\ref{section:seifert} below we consider
a Serre subcategory of an abelian category which is defined by $V'\in
\cy{C}$ if and only if $V\in\cy{C}$ and $V''\in\cy{C}$. For the
present we maintain greater generality; in particular an admissible
subcategory $\cy{C}$ is not required to be an abelian category.

\begin{definition}\label{define_metabolic}
Let $\zeta=1$ or $-1$. A non-singular $\zeta$-hermitian form
$(V,\phi)$ is called
metabolic if there is a submodule $L\subset V$ such that i) $L$ and
$V/L$ are in 
$\cy{C}$ and ii) $L=L^\perp$.
By definition $L^\perp=\Ker(j^*\phi:V\to L^*)$ where $j:L\to V$ is the
inclusion.
\end{definition}
\begin{definition}\label{define_Witt_group}
The Witt group $W^\zeta(\cy{C})$ is the abelian group
 with one generator $[V,\phi]$ for each isomorphism class of non-singular
 $\zeta$-hermitian forms $(V,\phi)\in \cy{C}$ subject
to relations 
\begin{equation*}
\begin{cases}
{[V',\phi']=[V,\phi]+[V'',\phi'']}, 
&\mbox{if $(V',\phi')\cong
(V,\phi)\oplus (V'',\phi'')$} \\
{[V,\phi]=0}, &\mbox{if $(V,\phi)$ is metabolic}.
\end{cases}
\end{equation*}
Two forms represent the same Witt class $[V,\phi]=[V',\phi']$ if
 and only if there exist metabolic forms  $(H,\eta)$ and $(H',\eta')$
such that 
\begin{equation*}
(V\oplus H,\phi\oplus\eta)\cong (V'\oplus H', \phi'\oplus\eta').
\end{equation*}  
\end{definition}
For example, a non-singular $\zeta$-hermitian form
$\phi:\Z^{2n}\to \Z^{2n}$ in the category $\Z$-proj (see
Example~\ref{duality_functor_on_A_Proj}) is metabolic if there exists a summand
$L\cong \Z^n$ with 
$\phi(L)(L)=0$.

In the category $\Bl(A)$ a metabolizer $L$ for a form $\phi:M\to
M^\wedge$ need not be a summand but $L$ and $M/L$ must lie in $\Bl(A)$
and one must have $\phi(L)(L)=0$ and $\phi(x)(L)\neq0$ if $x\notin L$.
Now that the notation is defined, we repeat equation~(\ref{identify_F_link_cobordism_with_W(Flink)}):
\begin{equation}\label{Duval_identification}
C_{2q-1}(F_\mu)\cong W^{(-1)^{q+1}}(\Bl(\Z)).\quad (q\geq3)
\end{equation}

\subsection{Duality-preserving
  functors}\label{section:duality_preserving_functors} 
A functor between hermitian categories which respects their structure
is called duality-preserving. Our first
examples will be functors induced by a morphism of rings with
involution. Duality-preserving functors will also play an essential
role in later sections (see
Theorem~\ref{main_covering_construction_theorem} above and 
Theorem~\ref{hermitian_Morita_equivalence} below).
\begin{definition}\label{define_dpres_funct}
A duality-preserving functor from $(\cy{C},\functor^*,i)$ to
$(\cy{D},\functor^*,i)$ is a triple $(G,\Psi,\eta)$ where
\begin{itemize}
\item $G:\cy{C}\to\cy{D}$ is a functor,
\item $\Psi=(\Psi_V)_{V\in\cy{C}}:G(\functor^*) \to G(\functor)^*$ is a natural isomorphism,
\item $\eta=1$ or $-1$
\end{itemize}
such that
\begin{equation}
\label{duality_functor}
\Psi_V^*i_{G(V)}= \eta \Psi_{V^*}G(i_V)
: G(V)\to G(V^*)^*
\end{equation}
for all $V\in \cy{C}$.
\end{definition}
We sometimes abbreviate $(G,\Psi,\eta)$ to $G$.

\begin{definition}
The composite of duality-preserving functors is defined by
\begin{equation}\label{how2_compose_dpres_functors}
(G,\Psi,\eta)\circ (G',\Psi',\eta')= (GG',\Psi G(\Psi') ,\eta\eta').
\end{equation}
\end{definition}
\begin{example}\label{simplest_ring_change_example}
A homomorphism $\nu:A\to A'$ of rings with involution induces a duality-preserving functor $(A'\otimes_A\functor, \Pi,1)$ from the category
$A\proj$ of finitely generated projective $A$-modules to the category
$A'\proj$ of finitely generated projective $A'$-modules. Explicitly
\begin{equation}\label{nat_iso_in_coeff_change_dpres_funct}
\begin{aligned}
{\Pi_V}:A'\otimes_A(V^*) &\xrightarrow{\cong} (A'\otimes_A V)^* \\
     a'_1\otimes \theta &\mapsto (a'_2\otimes x \mapsto
     a'_2\nu(\theta(x))\overline{a'_1}).
\end{aligned}
\end{equation}
for all $a_1',a_2'\in A'$, $\theta\in V^*$ and $x\in V$.
\end{example}
We are particularly concerned with the category $\Bl(A)$:
\begin{lemma}\label{Change_of_rings_Bl}
A homomorphism $\nu:A\to A'$ of rings with involution induces a canonical
duality-preserving functor
\begin{equation*}
(A'[F_\mu]\otimes_{A[F_\mu]}\functor,\Upsilon,1):\Bl(A)\to\Bl(A').
\end{equation*}
\end{lemma}
\noindent The natural
isomorphism $\Upsilon:A'[F_\mu]\otimes_{A[F_\mu]} \functor^\wedge\to
(A'[F_\mu]\otimes_{A[F_\mu]} \functor)^\wedge$ will be defined in the course
of the proof.
\begin{proof}
We saw in Lemma~\ref{Bl(functor)_is_functor} that there is a
functor $A'[F_\mu]\otimes_{A[F_\mu]}\functor$ from $\Bl(A)$ to
$\Bl(A')$. The dual $M^\wedge$ of a module $M\in\Bl(A)$ has presentation
\begin{equation*}
0\to (V[F_\mu])^*\xrightarrow{\sigma^*}(V[F_\mu])^* \to M^\wedge \to 0
\end{equation*}
(see Remark~\ref{explicit_formula_for_Mwedge_pres}).
Applying Example~\ref{simplest_ring_change_example} to $\nu:A[F_\mu]\to
A'[F_\mu]$ one obtains a natural isomorphism
\begin{equation*}
\Pi_{V[F_\mu]}:A'[F_\mu]\otimes_{A[F_\mu]}(V[F_\mu])^* \cong
 (A'[F_\mu]\otimes_{A[F_\mu]}V[F_\mu])^*
\end{equation*}
and hence a commutative diagram 
\begin{equation*}
\xymatrix@R=3ex@C=2ex{
0\ar[r] & A'[F_\mu]\otimes_{A[F_\mu]} (V[F_\mu])^* \ar[r]\ar[d]^{\cong} &
A'[F_\mu]\otimes_{A[F_\mu]} (V[F_\mu])^* \ar[d]^{\cong} \ar[r] &
A'[F_\mu]\otimes_{A[F_\mu]} M^\wedge \ar[r]\ar@{.>}[d]^{\Upsilon_M} & 0 \\
0\ar[r] & (A'[F_\mu]\otimes_{A[F_\mu]} V[F_\mu])^* \ar[r] & (A'[F_\mu]\otimes_{A[F_\mu]} V[F_\mu])^*
\ar[r] & (A'[F_\mu]\otimes_{A[F_\mu]} M)^\wedge \ar[r] & 0.}
\end{equation*}
 
We must check that the induced isomorphism 
\begin{equation*}
\Upsilon_M:A'[F_\mu]\otimes_{A[F_\mu]}
M^\wedge\to (A'[F_\mu]\otimes_{A[F_\mu]} M)^\wedge
\end{equation*}
is independent of the choice of presentation $\sigma$ and that $\Upsilon$ is
natural with respect to $M$. If $C$ denotes the chain complex
$V[F_\mu]\xrightarrow{\sigma}V[F_\mu]$ and
$C'$ denotes an alternative choice 
of resolution for $M$, say 
$C'=(V'[F_\mu]\xrightarrow{\sigma'}V'[F_\mu])$  then 
the identity map $\id:M\to M$ lifts to a chain equivalence $C\to C'$. 
The naturality of $\Pi$ in Example~\ref{simplest_ring_change_example}
implies that the diagram
\begin{equation*}
\xymatrix@=2ex{
A'[F_\mu]\otimes C^* \ar[r]\ar[d] & A'[F_\mu]\otimes(C')^*\ar[d] \\
(A'[F_\mu]\otimes C)^* \ar[r] & (A'[F_\mu]\otimes C')^*
}
\end{equation*}
commutes. The horizontal arrows induce the identity map on
$A'[F_\mu]\otimes M^\wedge$ and $(A'[F_\mu]\otimes M)^\wedge$ respectively, so the
vertical arrows induce the same map $\Upsilon_M$. The naturality of $\Upsilon$
follows similarly from the naturality of the transformation $\Pi$
in Example~\ref{simplest_ring_change_example}.
\end{proof}
\begin{remark}
Identifying $M^\wedge$ with $\Hom(M,A[F_\mu]_\Sigma/A[F_\mu])$ by
Proposition~\ref{identify_Mwedge_hom_to_loc}, an explicit formula for
$\Upsilon_M$ is  
\begin{align*}
\Upsilon_M:A'[F_\mu]\otimes_{A[F_\mu]}\Hom\left(M,\frac{A[F_\mu]_\Sigma}{A[F_\mu]}\right) &\to \Hom\left(A'[F_\mu]\otimes_{A[F_\mu]}
M\ , \frac{A'[F_\mu]_\Sigma}{A'[F_\mu]}\right) \\
a'_1\otimes \theta&\mapsto (a'_2\otimes m\mapsto a'_2\nu(\theta(m))\overline{a'_1})
\end{align*}
\end{remark}

We conclude this section by noting the effect of duality-preserving
functor on Witt groups:
\begin{lemma}\label{duality_preserving_functor_induces}
A duality-preserving functor $(G,\Psi,\eta):\cy{C}\to\cy{D}$ which
respects exact sequences induces a homomorphism of Witt groups
\begin{align*}
G:W^\zeta(\cy{C})&\to W^{\zeta\eta}(\cy{D}) \\
[V,\phi]&\mapsto[G(V),\Psi_VG(\phi)]
\end{align*}
\end{lemma}
\begin{proof}
See for example~\cite[p41-42]{She03mem}.
\end{proof}

\section{Intrinsic Invariants}\label{section:Blanchfield_form_invariants}
In~\cite{She03mem} the author defined invariants of the cobordism group
$C_{2q-1}(F_\mu)$ of $F_\mu$-links using the
identification (\ref{identify_F_link_cobordism_with_W(Seifert)}) due
to Ko~\cite{Ko87} and Mio~\cite{Mio87} of $C_{2q-1}(F_\mu)$ with a
Witt group of $\mu$-component Seifert forms, denoted
$W^{(-1)^q}(\Sei(\Z))$ below ($q\geq3$). To 
distinguish $F_\mu$-links one first chooses a Seifert surface for each and
then computes invariants of the associated Seifert forms.

In the present section we define $F_\mu$-link cobordism invariants via 
Duval's identification $C_{2q-1}(F_\mu)\cong W^{(-1)^q}(\Bl(\Z))$. The
definitions will parallel those in~\cite{She03mem} and we shall prove in
Section~\ref{section:Proof_that_B_induces_Witt_isomorphism} that the
invariants obtained are equivalent. Whereas Seifert forms
are convenient for computing the invariants in explicit examples, the 
Blanchfield-Duval form has the advantage that it is defined without
making a choice of Seifert surface. 
\subsection{Overview}\label{section:overview_of_BD_form_invariants}
Let $\zeta=1$ or $-1$. The inclusion $\Z\subset \Q$ induces a
duality-preserving functor
$\Bl(\Z)\to\Bl(\Q)$ (see Lemma~\ref{Change_of_rings_Bl} above) and
hence a homomorphism of Witt groups 
\begin{equation}\label{ZintoQ_W(BDform)_morphism}
W^\zeta(\Bl(\Z))\to
W^\zeta(\Bl(\Q)).
\end{equation}
It follows from Theorem~\ref{main_covering_construction_theorem}
that~(\ref{ZintoQ_W(BDform)_morphism}) is an injection; see the proof
of Corollary~\ref{BD_form_invariants_for_Flinks_suffice:precise} below.
We proceed to compute $W^\zeta(\Bl(\Q))$
in three steps, which were outlined in the introduction. We list them
again here in more detail:
\begin{enumerate}
\item {\em Devissage}. We prove that $\Bl(\Q)$ is an abelian category
with ascending and descending chain conditions. Recall that a module
$M\in\Bl(\Q)$ is called simple (or irreducible) if $M\neq0$ and there
are no submodules of $M$ in $\Bl(\Q)$ other than $0$ and $M$.
If $M$ is a simple module then $\Bl(\Q)|_M\subset
\Bl(\Q)$ denotes the full subcategory in which the objects are direct
sums of copies of $M$. Recall that $M$ is called $\zeta$-self-dual
if there is an isomorphism $b:M\to M^\wedge$ such that
$b^\wedge=\zeta b$. We obtain, by ``hermitian devissage'', the
decomposition  
\begin{equation}\label{devissage_applied_to_flk(Q)}
W^\zeta(\Bl(\Q))~\cong \bigoplus W^\zeta(\Bl(\Q)|_M) 
\end{equation}
with one summand for each isomorphism class of $\zeta$-self-dual
simple $F_\mu$-link modules $M$. Let $p_M$ denote the projection of
$W^\zeta(\Bl(\Q))$ onto $W^\zeta(\Bl(\Q)|_M)$.
\item {\em Morita equivalence}. For each $\zeta$-self-dual simple
module $M$ we choose a non-singular $\zeta$-hermitian form $b:M\to
M^\wedge$. We obtain by hermitian Morita equivalence an isomorphism
\begin{equation}\label{hermitian_morita_equivalence_on_Flk(Q)}
\Theta_{M,b}:W^\zeta(\Bl(\Q)|_M) \to W^1(E)
\end{equation}
where $E=\End_{\Q[F_\mu]}M$ is the endomorphism ring of $M$ and is endowed
with the involution $f\mapsto b^{-1}f^\wedge b$.
By Schur's Lemma $E$ is a division ring and, as we
discuss, $E$ turns out to be of finite dimension over $\Q$.
\item We recall from the literature invariants of each group $W^1(E)$. In
most cases some combination of dimension modulo $2$, signatures,
discriminant and Hasse-Witt invariant are sufficient to distinguish
the elements of $W^1(E)$ (see the
table~(\ref{table_of_Witt_invariants}) below). One class of division
algebra with involution requires a secondary invariant, such as the
Lewis $\theta$, which is defined only if all the other invariants vanish.
\end{enumerate}
Let us make two remarks about the modules $M$ which appear in item 1.
Firstly, every simple module $M\in\Bl(\Q)$ such that $M\cong M^\wedge$
is either $1$-self-dual or $(-1)$-self-dual (or both). Secondly, $M$ is both
$1$-self-dual and $(-1)$-self-dual if and only if the
involution $f\mapsto b^{-1}f^\wedge b$ induced on
$E=\End_{\Q[F_\mu]}(M)$ is not the identity map for some (and therefore
every) $\zeta$-self-dual form $b:M\to M^\wedge$. See lemmas~5.5
and~5.6 of~\cite{She03mem} for details.
\begin{example}
In the knot theory case $\mu=1$ we can add simplifying remarks to each
of the three steps:
\begin{enumerate}
\item A simple self-dual module $M\in\Bl(\Q)$ may be written
\begin{equation*}
M=\Q[z,z^{-1}]/(p)
\end{equation*}
where $p$ is an irreducible polynomial and
$(p(z^{-1}))=(p(z)) \vartriangleleft \Q[z,z^{-1}]$. 
\item The endomorphism ring $E=\End_{\Q[z,z^{-1}]}(M)$ may also be
written as a quotient $\Q[z,z^{-1}]/(p)$ and is an
algebraic number field of finite dimension over $\Q$. The involution
on $E$ is given by $z\mapsto z^{-1}$ and does not depend on 
the choice of form $b$.

Setting aside the case $M=\Q[z,z^{-1}]/(1+z)$, the involution on $E$ is not
 the identity so every self-dual $M$ is both $1$-self-dual and
 $(-1)$-self-dual.
The exceptional module $M=\Q[z,z^{-1}]/(1+z)$ is only $1$-self-dual
 but plays little role since $-1$ is not a root of any
 polynomial $p\in\Z[z,z^{-1}]$ such that $p(1)=\pm1$. In other words, 
the projection of $W^1(\Bl(\Z))$ on this exceptional summand of
 $W^1(\Bl(\Q))$ is zero.
\item As discussed in 2., one need only consider the Witt groups
$W^1(E)$ of number fields with non-trivial involution, or in other
words, hermitian forms over number fields. The dimension modulo $2$,
signatures and discriminant are sufficient to distinguish the elements
of $W^1(E)$. 
\end{enumerate}
Equation~(\ref{iso_class_of_knot_cobordism_gp}) can be derived as a
consequence of this analysis. 
\end{example}

Returning to the general case $\mu\geq1$ and putting together steps
1-3.~we obtain the following restatement of
Theorem~\ref{BD_form_invariants_for_W_suffice:vague}.
\begin{theorem}\label{BD_form_invariants_for_W_suffice:precise}
Let $\zeta=1$ or $-1$. An element $\alpha\in W^\zeta(\Bl(\Q))$ is
equal to zero if and only if for each finite-dimensional
$\zeta$-self-dual simple $F_\mu$-link module $M$ and non-singular
$\zeta$-hermitian form $b:M\to M^\wedge$, the dimension
modulo~$2$, the signatures, the discriminant, the Hasse-Witt invariant
and the Lewis $\theta$-invariant of 
\begin{equation*}
\Theta_{M,b}\,p_M(\alpha) \in W^1(\End_{\Q[F_\mu]}M)
\end{equation*}
are trivial (if defined).\hspace*{\fill}\qed
\end{theorem}
Note that if the invariants corresponding to one form
$b:M\to M^\wedge$ are trivial then $p_M(\alpha)=0$ so the invariants
are trivial for any other choice $b':M\to M^\wedge$.
We now restate Corollary~\ref{BD_form_invariants_for_Flinks_suffice:vague}:
\begin{corollary}\label{BD_form_invariants_for_Flinks_suffice:precise}
Suppose $(L^0,\theta^0)$ and $(L^1,\theta^1)$ are $(2q-1)$-dimensional
$F_\mu$-links, where $q>1$. Let $\overline{X_i}$ denote the free cover
of the complement of $L_i$, let $N_i=H_q(\overline{X_i})/(\Z{\rm
-torsion})$ and let $\phi^i:N_i \to \Hom(N_i,\Z[F_\mu]_\Sigma/\Z[F_\mu])$ denote the
Blanchfield-Duval form for $(L^i,\theta^i)$.

The $F_\mu$-links $(L^0,\theta^0)$ and $(L^1,\theta^1)$ are cobordant
if and only if  
for each finite-dimensional $\zeta$-self-dual simple
$F_\mu$-link module $M$ and each non-singular $\zeta$-hermitian form 
$b:M\to M^\wedge$, the dimension modulo~$2$, the
signatures, the discriminant, the Hasse-Witt invariant and the Lewis
$\theta$-invariant of 
\begin{equation*}
\Theta_{M,b}\,p_M [\Q\otimes_\Z(N^0\oplus
N^1,\phi^0\oplus-\phi^1)] \in W^1(\End_{\Q[F_\mu]}M) 
\end{equation*}
are trivial (if defined).\hspace*{\fill}\qed 
\end{corollary}
\begin{proof}
We deduce Corollary~\ref{BD_form_invariants_for_Flinks_suffice:precise}
from Theorem~\ref{BD_form_invariants_for_W_suffice:vague}
(=Theorem~\ref{BD_form_invariants_for_W_suffice:precise}) 
and Theorem~\ref{main_covering_construction_theorem}. As we remarked in the
introduction,
Corollary~\ref{BD_form_invariants_for_Flinks_suffice:precise} also 
follows from Theorem~\ref{theorem_relating_invariants:vague} and
Theorem~B of~\cite{She03mem}.

Proposition~\ref{B_respects_coefficient_change} below says that 
the duality-preserving functor 
\begin{equation*}
(B,\Phi,-1):\Sei(A)\to\Bl(A)
\end{equation*}
of Section~\ref{section:covering_construction}
respects coefficient change so there is a commutative diagram
\begin{equation}\label{B_respects_Z_into_Q}
\begin{gathered}
\xymatrix@R=2ex@C=7ex{W^{(-1)^q}(\Sei(\Z))\ar[d]\ar[r]^B &W^{(-1)^{q+1}}(\Bl(\Z))\ar[d] \\
W^{(-1)^{q}}(\Sei(\Q))\ar[r]_B & W^{(-1)^{q+1}}(\Bl(\Q)).}
\end{gathered}
\end{equation}
The category $\Sei(A)$ is defined in Section~\ref{section:seifert_modules}.
The lower horizontal map in~(\ref{B_respects_Z_into_Q}) is an isomorphism by
Theorem~\ref{main_covering_construction_theorem} and the
upper horizontal map is an isomorphism by (\ref{identify_F_link_cobordism_with_W(Seifert)}) and (\ref{identify_F_link_cobordism_with_W(Flink)}) above
(see also Remark~\ref{remark_that_B(Seifert_form)=BDform}). It is easy
to prove that the left hand vertical map is an injection (see for
example Lemma~11.1 of~\cite{She03mem}). Thus the right-hand vertical map
is also an
injection. Corollary~\ref{BD_form_invariants_for_Flinks_suffice:precise}
therefore follows from Theorem~\ref{BD_form_invariants_for_W_suffice:precise}
\end{proof}
\subsection{Step 1: Devissage}
Let us briefly recall some definitions.
If $\cy{\cy{A}}$ is an additive category and $M$ is an object in $\cy{A}$
the symbol $\cy{A}|_M$ denotes the full subcategory such that 
$N\in\cy{A}$ if and only if $N$ is a summand of some finite direct sum
of copies of $M$. 
If $\cy{A}$ is a hermitian category and $M$
is self-dual then $\cy{A}|_M$ is a hermitian subcategory. 

Suppose now that $\cy{A}$ is an abelian category. 
A non-zero module $M$ in $\cy{A}$
is called simple (or irreducible) if there are no submodules of $M$ in
$\cy{A}$ other than $0$ and $M$. 
The category $\cy{A}$ has both ascending and descending chain conditions
if and only if every module $M$ in $\cy{A}$ has a finite composition
series
\begin{equation}\label{composition_series}
0=M_0\subset M_1\subset M_2\subset M_3\subset \cdots \subset M_s=M
\end{equation}
where $M_i/M_{i-1}$ is simple for $i=1,\cdots,s$.
If $M$ is simple then every module in
$\cy{A}|_M$ is a direct sum of copies of $M$.
Let $\zeta=1$ or $-1$. If $(\cy{A},*,i)$ is a hermitian category then
a module $M$ is called $\zeta$-self-dual if there is an isomorphism
$\phi:M\to M^*$ such that $\phi^*i_M=\zeta\phi$. 

The general decomposition theorem we need is the following:
\begin{theorem}[Devissage]
\label{hermitian_devissage}
Suppose $\cy{A}$ is an abelian hermitian category with ascending and
descending chain conditions.  There is an isomorphism of Witt groups
\begin{equation*}
W^\zeta(\cy{C})~\cong \bigoplus~
W^\zeta(\cy{C}|_M)
\end{equation*}
with one summand for each isomorphism class of simple $\zeta$-self-dual
modules in $\cy{C}$.
\end{theorem}
\begin{proof}
See~\cite[Theorem 5.3]{She03mem} or~\cite{QSS79}.
\end{proof}
To prove equation~(\ref{devissage_applied_to_flk(Q)}) it therefore
suffices to show: 
\begin{proposition}\label{Blanchfield_category_is_abelian_with_chain_conditions}
If $k$ is a (commutative) field then the category $\Bl(k)$ is an abelian
category with ascending and descending chain conditions. 
\end{proposition}
We take coefficients in a field for simplicity. With little extra work
one can show that
Proposition~\ref{Blanchfield_category_is_abelian_with_chain_conditions} 
holds when $k$ is replaced by any semi-simple Artinian ring. See also
Remark~\ref{remark_on_nu_torsion} below.
\begin{caveat} When $\mu\geq2$ a simple module in the category
$\Bl(k)$ is not a simple module in the category of $k[F_\mu]$-modules.
\end{caveat}
We prove
Proposition~\ref{Blanchfield_category_is_abelian_with_chain_conditions} 
as follows: We first show that $\Bl(k)$ is a 
subcategory of the category $\cy{T}_{k[F_\mu]}$ of ``torsion'' modules
which P.M.Cohn introduced and proved to be an abelian
category with ascending and descending chain conditions (see
Proposition~\ref{torsion_category_is_abelian} below). After giving
details of Cohn's work, we conclude the proof of
Proposition~\ref{Blanchfield_category_is_abelian_with_chain_conditions} 
by checking that $\Bl(k)$ is closed under direct sums and that the
kernel and cokernel of every morphism in $\Bl(k)$ again lie in $\Bl(k)$.

\subsubsection{Firs and torsion modules}\label{section:firs_and_torsion_modules}
To describe Cohn's results 
we must state some properties of the group ring $k[F_\mu]$.
\begin{definition}
A ring $R$ has invariant basis number (IBN) if 
$R^n\cong R^m$ implies $n=m$. In other words $R$ has IBN if every
finitely generated free left $R$-module has unique rank.
\end{definition}
The existence of the augmentation $\epsilon:k[F_\mu]\to k$ implies that
$k[F_\mu]$ has IBN for if $k[F_\mu]^n\cong k[F_\mu]^m$ then
\begin{equation*}
k^n\cong k\otimes_{k[F_\mu]}k[F_\mu]^n\cong
 k\otimes_{k[F_\mu]}k[F_\mu]^m\cong k^m
\end{equation*}
so $m=n$.

Note that if $R$ has IBN then one can use the duality functor
$\Hom(\functor,R)$ to prove that finitely generated free right
$R$-modules also have unique rank.
\begin{definition}
An associative ring $R$ is called a {\it free ideal ring (fir)} if $R$
has IBN, every left ideal in $R$ is a free left $R$-module and every right
ideal is a free right $R$-module.
\end{definition}
Cohn showed in~\cite[Corollary 3]{Coh64} that if $k$ is a field then the
group ring $k[F_\mu]$ of the free group $F_\mu$ is a fir.

If $R$ is a fir then every submodule of a free $R$-module is
free (Cohn~\cite[p71]{Coh85}). Hence every $R$-module has a
presentation $0\to F_1\to F_0\to M\to0$ where $F_1$ and $F_0$ are
free.
\begin{definition}
If an $R$-module $M$ has a presentation
\begin{equation}\label{present_torsion_module}
0\to R^n\xrightarrow{\sigma} R^m\xrightarrow{p} M \to0
\end{equation}
the {\it Euler characteristic} of $M$ is $\chi(M)=m-n$.
\end{definition}
If a finite presentation exists the Euler
characteristic is independent of the choice of presentation by
Schanuel's Lemma. Note also that an exact sequence $0\to M\to M'\to
M''\to0$ of finitely presented modules implies the equation
$\chi(M')=\chi(M)+\chi(M'')$ (compare the diagram
(\ref{horseshoe_diagram}) below).

If $R$ is a fir then the category of finitely presented $R$-modules
(and $R$-module maps) is
an abelian category; in other words direct sums of finitely presented
modules are finitely presented and the cokernel and kernel of a map
between finitely presented modules are finitely presented.
In fact, if $R$ is any ring such that every finitely generated
one-sided ideal is finitely related then the finitely presented
$R$-modules form an abelian category (e.g.~Cohn~\cite[p554-556]{Coh85}).

\begin{definition}
A morphism $\sigma:R^n\to R^n$ between free left
$R$-modules is called {\it full} if every factorization
\begin{equation*}
\xymatrix{
R^n\ar@/^1pc/[rr]^\sigma
\ar[r]_{\sigma_2} & F \ar[r]_{\sigma_1} & R^n
}
\end{equation*}
where $F$ is a free module has $\Rank(F)\geq n$.
\end{definition}

\begin{lemma}\label{characterize_torsion}
Suppose $R$ is a fir and $M$ is a finitely presented $R$-module with
$\chi(M)=0$. The following are equivalent:
\begin{enumerate}
\item\label{all_full} In every finite
presentation~(\ref{present_torsion_module}) of $M$, the map $\sigma$ is full. 
\item\label{exists_full} There exists a
presentation~(\ref{present_torsion_module}) such that $\sigma$ is full.
\item\label{Euler_sub} $\chi(N)\geq0$ for all finitely
generated submodules $N$ of $M$.
\item\label{Euler_quot} $\chi(M/N)\leq0$ for all
finitely generated submodules $N$ of $M$.
\end{enumerate}
\end{lemma} 
\begin{proof}
The implication \ref{all_full}$\Rightarrow$~\ref{exists_full} is
immediate. To show~\ref{exists_full} implies~\ref{Euler_sub},
suppose we are given a presentation~(\ref{present_torsion_module})
such that $m=n$ and $\sigma$ is full. If $N$ is a finitely
generated submodule of $M$ then
$\sigma(R^n)\subset p^{-1}(N)\subset
R^n$. Now $p^{-1}(N)$ is a free module because
$R$ is a fir and $p^{-1}(N)$ has rank at least $n$ since $\sigma$ is full.
The exact sequence
\begin{equation*}
0\to R^n\to p^{-1}(N) \to N \to 0
\end{equation*}
implies that $\chi(N)\geq0$. This completes the proof that~\ref{exists_full}
implies~\ref{Euler_sub}.

The equation $\chi(M)=\chi(N)+\chi(M/N)$ implies that~\ref{Euler_sub}
and~\ref{Euler_quot} are equivalent so we can conclude the proof of
the Lemma by showing that~\ref{Euler_sub} implies~\ref{all_full}. The equation
$\chi(M)=0$ says that every finite presentation
(\ref{present_torsion_module}) has $m=n$. We must prove that $\sigma$ is
full. Suppose $\sigma$ can be written as a composite
$R^n\xrightarrow{\sigma_2}R^k\xrightarrow{\sigma_1}R^n$. We aim to show
$k\geq n$. Now
$\sigma(R^n)$ and $\sigma_1(R^k)$ are free modules since $R$ is a fir and 
\begin{equation*}
\frac{\sigma_1(R^k)}{\sigma(R^n)}\subseteq
\frac{R^n}{\sigma(R^n)}=M. 
\end{equation*}
By statement~\ref{Euler_sub}, $\chi(\sigma_1(R^k)/\sigma(R^n))\geq0$ so
$\Rank(\sigma_1(R^k))\geq\Rank(\sigma(R^n))=n$. Thus $k\geq n$
and hence $\sigma$ is full.
\end{proof}

A module which satisfies the equivalent conditions in
Lemma~\ref{characterize_torsion} is called a {\it torsion} module. The
symbol $\cy{T}_R$ denote the category of torsion modules and $R$-module maps.

\begin{lemma}(Cohn \cite[p166]{Coh85})\label{extension_of_torsion_modules}
If $0\to M\to M'\to M''\to 0$ is an exact sequence and $M$ and $M''$
are torsion modules then $M'$ is a torsion module. In particular a
direct sum of torsion modules is again a torsion module.
\end{lemma}
\begin{proof}
If $M$ and $M''$ are finitely presented then $M'$ is also finitely presented
(compare~(\ref{horseshoe_diagram}) below). Since $\chi(M)=\chi(M'')=0$
we have $\chi(M')=0$. Now if $N\leq M'$ is finitely generated it
suffices by Lemma~\ref{characterize_torsion} to show that
$\chi(N)\geq0$. Note first that $\chi(N)=\chi(N\cap 
M)+\chi(N/(N\cap M))$. Now $N\cap M \leq M$ so $\chi(N\cap
M)\geq0$ and $N/(N\cap M)\cong (N+M)/M \leq M''$ so $\chi(N/(N\cap
M))\geq0$. Thus $\chi(N)\geq0$.
\end{proof}
\begin{proposition}(Cohn
\cite[p167,234]{Coh85})\label{torsion_category_is_abelian} Suppose $R$ is
a fir.
\begin{enumerate}
\item\label{abelian_part} The category $\cy{T}_R$ of torsion modules is an
abelian category.
\item\label{ACC_DCC_part} Every module in $\cy{T}_R$ has a finite
composition series~(\ref{composition_series}) in $\cy{T}_R$.
\end{enumerate}
\end{proposition}
\begin{proof}
To establish statement~\ref{abelian_part}, it suffices to show that if
$M$ and $M'$ lie in $\cy{T}_R$ then $M\oplus M'\in\cy{T}_R$ and
the kernel and cokernel of every morphism $f:M\to M'$ are in
$\cy{T}_R$. Lemma~\ref{extension_of_torsion_modules} gives $M\oplus
M'\in\cy{T}_R$. Suppose then that $f:M\to M'$ is an $R$-module
morphism. Since the finitely presented $R$-modules are an abelian
category the kernel, image and cokernel of $f$ are finitely presented.
Now $f(M)$ is a submodule of $M'$ and a quotient module of $M$ so
$\chi(f(M))\geq0$ and $\chi(f(M))\leq0$ by
Lemma~\ref{characterize_torsion}. Thus $\chi(f(M))=0$ and it follows
that $\chi(\Ker(f))=0$ and $\chi(\Coker(f))=0$. 
Every finitely generated submodule $N\leq\Ker(f)$ is a submodule of $M$
so $\chi(N)\geq0$ and hence $\Ker(f)\in\cy{T}_R$. Similarly every quotient
$N'$ of $\Coker(f)$ is a quotient of $M'$ and hence has
$\chi(N')\leq0$. Thus $\Coker(f)\in\cy{T}_R$ also.  

To prove part 2.~of Proposition~\ref{torsion_category_is_abelian}
we note first that it is sufficient to check the ascending chain
condition.
Indeed, if $M\in\cy{T}_R$ then
$M^\wedge=\Ext_R(M,R)$ is a torsion right $R$-module and a descending chain 
\begin{equation*}
M=M_0\supseteq M_1\supseteq M_2 \supseteq \cdots
\end{equation*}
in $\cy{T}_R$ gives rise to an ascending chain
\begin{equation*}
\left(\frac{M}{M_0}\right)^\wedge \subseteq
\left(\frac{M}{M_1}\right)^\wedge \subseteq
\left(\frac{M}{M_2}\right)^\wedge \subseteq \cdots
\end{equation*}

In fact Cohn showed that a larger class of modules, the finitely
related bound modules, have the ascending chain condition.
\begin{definition}
An $R$-module $M$ is {\it bound} if $\Hom(M,R)=0$.
\end{definition}
\begin{lemma}\label{torsion_implies_bound}
Every torsion module over a fir is bound.
\end{lemma}
\begin{proof}
If $M$ is an $R$-module and $\theta:M\to R$ then $\theta(M)$ is a free
module so $M\cong \Ker(\theta)\oplus \theta(M)$. Now $\theta(M)$ is a
quotient module of $M$ so if $M$ is a torsion module then
$\chi(\theta(M))\leq0$. It follows that $\theta(M)=0$ so
$\theta=0$. Thus $\Hom_R(M,R)=0$.
\end{proof}
\begin{lemma}(Cohn~\cite[p231]{Coh85})
If $R$ is a fir and $M$ is a finitely related $R$-module then every
bound submodule of $M$ is finitely presented.
\end{lemma}
\begin{proof}
There is an exact sequence $0\to R^n\to F\xrightarrow{p} M\to0$ where $F$ is a
free $R$-module. A submodule $B\leq M$ has presentation
\begin{equation*}
0\to R^n\to p^{-1}(B)\to B\to0
\end{equation*}
and $p^{-1}(B)$ is a free module since $R$ is a fir. The image of
$R^n$ is contained in a finitely generated summand of $p^{-1}(B)$ so
$B$ is a direct sum of a free module and a finitely presented
module. If $B$ is bound then $B$ does not have any non-zero free
summand so $B$ itself is finitely presented. 
\end{proof}
Thus every torsion submodule of a torsion module $M$ is finitely
generated so we have the ascending chain condition on torsion
submodules. It follows that $\cy{T}_R$ has both ascending and
descending chain conditions and the proof of
Proposition~\ref{torsion_category_is_abelian} is complete.
\end{proof}
We are now in a position to
deduce Proposition~\ref{Blanchfield_category_is_abelian_with_chain_conditions}.
\begin{lemma}\label{Bl(k)_is_inside_T}
The category $\Bl(k)$ is a full subcategory of the category
$\cy{T}_{k[F_\mu]}$ of torsion modules.
\end{lemma}
\begin{proof}
It suffices to show that every module in $\Bl(k)$ is a torsion
 module. By definition, a module in $\Bl(k)$ has a
 presentation~(\ref{Presentation_of_F-link_module}) where
 $\epsilon(\sigma)$ is an isomorphism. Since $\epsilon(\sigma)$ is full we
 can deduce that $\sigma$ is full since each factorization of $\sigma$
 induces a corresponding factorization of $\epsilon(\sigma)$.
\end{proof}
\begin{proof}[Proof of Proposition~\ref{Blanchfield_category_is_abelian_with_chain_conditions}]
Since $\cy{T}_{k[F_\mu]}$ has ascending and descending chain conditions
and $\Bl(k)\subset \cy{T}_{k[F_\mu]}$ we
need only show that $\Bl(k)$ is an abelian category. 

Suppose $M$ and $M'$ are in $\Bl(k)$. It follows directly from the
definition of $\Bl(k)$ that $M\oplus M'\in\Bl(k)$. We must
show that the kernel and cokernel of each map $f:M\to M'$
lie in $\Bl(k)$. By Proposition~\ref{torsion_category_is_abelian} the
kernel, image and cokernel all lie in
$\cy{T}_{k[F_\mu]}$ so 
\begin{equation*}
\chi(\Ker(f))=\chi(f(M))=\chi(\Coker(f))=0.
\end{equation*}
 After choosing
presentations $\sigma$ and $\sigma'$ for $\Ker(f)$ and $f(M)$
respectively one can fill in the dotted arrows below to obtain a
commutative diagram with exact rows and exact columns:
\begin{equation}\label{horseshoe_diagram}
\begin{gathered}\xymatrix@C=3ex@R=2ex{
& 0 \ar@{.>}[d] & 0\ar@{.>}[d] & 0\ar[d] \\
0 \ar[r] & k[F_\mu]^n \ar@{.>}[d] \ar[r]^\sigma & k[F_\mu]^n\ar[r]\ar@{.>}[d] &
\Ker(f) \ar[r] \ar[d] & 0   \\
0 \ar@{.>}[r] & k[F_\mu]^{n+n'} \ar@{.>}[r]^{\sigma''}\ar@{.>}[d] & k[F_\mu]^{n+n'}
\ar@{.>}[r]\ar@{.>}[d] & M \ar[r]\ar[d] & 0 \\
0 \ar[r] & k[F_\mu]^{n'} \ar[r]^{\sigma'} \ar@{.>}[d] & k[F_\mu]^{n'} \ar[r]
\ar@{.>}[d] & f(M)\ar[d] \ar[r] & 0  \\ 
& 0 & 0 & 0}
\end{gathered}
\end{equation}
The map $\sigma''$ is given by $\left(\begin{matrix} \sigma & \tau \\ 0 &
\sigma'\end{matrix}\right)$ for some $\tau$. Since $M\in\Bl(k)$, the
augmentation $\epsilon(\sigma'')$ is an isomorphism by
Lemma~\ref{Sato_condition} above. It follows that $\epsilon(\sigma)$
and $\epsilon(\sigma'')$ are isomorphisms and hence $\Ker(f)$ and
$f(M)$ are in $\Bl(k)$. The same argument, applied to the exact sequence
$0\to  f(M)\to M'\to \Coker(f)\to 0$ shows that $\Coker(f)$ is also in
$\Bl(k)$. 
\end{proof}
\begin{remark}\label{remark_on_nu_torsion}
The arguments above can be adapted to generalize
 Propositions~\ref{torsion_category_is_abelian}
 and~\ref{Blanchfield_category_is_abelian_with_chain_conditions} as
 follows. Suppose every finitely generated one-sided ideal of a ring $R$ is a
 projective module (i.e.~$R$ is semi-hereditary). Suppose $S$ is a ring with
 the property that $S^n\oplus P\cong S^n$ implies $P=0$ (i.e.~$S$ is weakly
 finite) and $\nu:R\to S$ is a ring homomorphism. The Grothendieck group
 $K_0(S)$ admits a partial order in which
 $x\leq y$ if and only $y-x$ lies in the positive cone 
\begin{equation*}
\{[P] \mid P\in S\proj\}\subset K_0(S) .
\end{equation*}
Now the following category $\cy{T}_\nu$ is an abelian category: 
An $R$-module $M$ lies in $\cy{T}_\nu$ if $M$ has a finite presentation
 by projective modules $0\to P_1\xrightarrow{\sigma} P_0 \to M \to0$ such that
 $1\otimes\sigma:S\otimes_R P_1\to S\otimes P_0$ is full with respect
 to the partial order on $K_0(S)$. 
The subcategory of modules for which $1\otimes\sigma$
 is invertible is also an abelian category. Under the additional
 hypotheses that all the one-sided ideals of $R$ are projective
 modules (i.e.~$R$ is hereditary) and that the equation $S\otimes_R
 P=0$ implies $P=0$ for projective $R$-modules $P$, one can also
 conclude that these abelian categories have ascending and descending
 chain conditions. One recovers
 Proposition~\ref{torsion_category_is_abelian} when $R$ is a fir by
 setting $\nu=\id:R\to R$ and one recovers
 Proposition~\ref{Blanchfield_category_is_abelian_with_chain_conditions}
 by setting $\nu=\epsilon:k[F_\mu]\to k$.
\end{remark}
\subsection{Step 2: Morita Equivalence}
Having reduced $W^\zeta(\Bl(\Q))$ to a direct sum of Witt groups
$W^\zeta(\Bl(\Q)|_M)$ where $M$ is simple and $\zeta$-self-dual
(see equation~(\ref{devissage_applied_to_flk(Q)})) we pass next
from $W^\zeta(\Bl(\Q)|_M)$ to the Witt group
$W^1(\End_{\Q[F_\mu]}M)$ of the endomorphism ring of $M$.

Recall that if $\cy{C}$ is an additive category then an object in
$\cy{C}|_M$ is a summand of a direct sum of copies of $M$.
The general theorem we employ in this section is the following:
\begin{theorem}[Hermitian Morita Equivalence]
\label{hermitian_Morita_equivalence}
Let $\eta=+1$ or $-1$. Suppose that $b:M\to M^*$ is a non-singular
$\eta$-hermitian form in a hermitian category $\cy{C}$, and assume
further that every idempotent endomorphism in the hermitian
subcategory $\cy{C}|_M$ splits. Let $E=\End_\cy{C}M$
be endowed with the involution $f\mapsto \overline f=b^{-1}f^*b$. 
Then there is an equivalence of hermitian categories
\begin{equation*}
\Theta_{M,b}=(\Hom(M,\functor), \Omega^b, \eta) : \cy{C}|_M\to E\proj
\end{equation*}
where for $N\in\cy{C}|_M$ the map
$\Omega^b_N(\gamma)=(\alpha\mapsto\eta b^{-1}\alpha^*\gamma)$ is the
composite of natural isomorphisms
\begin{equation}
\label{natural_isomorphism_Phi}
\begin{aligned}
\Hom_\cy{C}(M,N^*)&\to\Hom_\cy{C}(N,M)\to \Hom_E(\Hom_\cy{C}(M,N),E) \\
\gamma &\mapsto b^{-1}\gamma^*;\hspace{8mm}\delta\mapsto
(\alpha\mapsto \overline{\delta\alpha}).
\end{aligned}
\end{equation}
\end{theorem}
\begin{proof}
See~\cite[Theorem 4.7]{She03mem},~\cite[\S I.9,ch.II]{Knu91} or~\cite{QSS79}.
\end{proof}
The following is a corollary of
Theorem~\ref{hermitian_Morita_equivalence} and
Lemma~\ref{duality_preserving_functor_induces}.
\begin{corollary}\label{apply_morita_eq_to_Bl}
If $M$ is a simple module in
$\Bl(k)$ and $b:M\to M^\wedge$ is a non-singular $\zeta$-hermitian
form then the duality-preserving functor $\Theta_{M,b}$ of
Theorem~\ref{hermitian_Morita_equivalence} induces an isomorphism of
Witt groups
\begin{align*}
W^\zeta(\Bl(k)|_M)&\to W^1(\End_{k[F_\mu]}M) \\
[N,\phi] &\mapsto [\Hom(M,N),\Omega^b_N\phi_*]
\end{align*}
where $\phi_*:\Hom(M,N)\to\Hom(M,N^*)$ and 
\begin{equation*}
(\Omega^b_N\phi_*)(\alpha)(\beta)=\zeta
b^{-1}\beta^*\phi\alpha \in \End_{k[F_\mu]}M
\end{equation*}
for all $\alpha,\beta\in
\Hom(M,N)$. Equation~(\ref{devissage_applied_to_flk(Q)}) implies that
\begin{equation*}
W^\zeta(\Bl(k))~\cong \bigoplus_M W^1(\End_{k[F_\mu]}M)~.
\end{equation*}
with one summand for each isomorphism class of $\zeta$-self-dual
simple $F_\mu$-link modules $M$.
\end{corollary}
\begin{proof}
Every exact sequence in $\Bl(k)|_M$ splits, so $\Theta_{M,b}$ is
exact and hence induces a morphism of Witt groups
\begin{equation*}
W^\zeta(\Bl(k)|_M) \to W^1(\End_{k[F_\mu]}(M)).
\end{equation*}
Since $\Theta_{M,b}$ is an equivalence of hermitian categories it
follows by Lemma~\ref{nat_isomorphic_dpres_funct_give_same_Witt_homo}
of Appendix~\ref{section:naturality_of_B} that the induced map of Witt
groups is an isomorphism.
\end{proof}
Equation~(\ref{hermitian_morita_equivalence_on_Flk(Q)}) above is a
 special case of corollary~\ref{apply_morita_eq_to_Bl} so the following
 proposition completes step~2:
\begin{proposition}\label{Bl_endo_rings_fin_diml}
The endomorphism ring of every module in $\Bl(k)$
is of finite dimension over $k$.
\end{proposition}
Proposition~\ref{Bl_endo_rings_fin_diml} follows from
Theorem~\ref{Bbar_is_equiv_on_fgproj_cats} and
Lemma~\ref{fully_simplify_quotient} below.
Theorem~\ref{Bbar_is_equiv_on_fgproj_cats} can be considered analogous
to the geometric fact that one can choose a Seifert surface for an
$F_\mu$-link. Since the homology of a Seifert surface is 
finite-dimensional the endomorphism ring of the associated Seifert
module is finite-dimensional.  Using chapter~12
of~\cite{She03mem}, part 2.~of
Lemma~\ref{simple_goes_to_comes_from_simple} and part 2.~of 
Lemma~\ref{quotient_lemma_with_duality} we can also deduce that every division ring 
with involution which is finite-dimensional over $\Q$ arises as
\begin{equation*}
(\End_{\Q[F_\mu]}M\ ,\ f\mapsto\overline f=b^{-1}f^\wedge b)
\end{equation*}
for some pair $(M,b)$
where $M$ is a simple module in $\Bl(\Q)$. 

The proofs of the results cited in the previous paragraph 
do not use Proposition~\ref{Bl_endo_rings_fin_diml} (i.e.~the arguments
presented are not circular). However, the spirit of
this section is to define invariants of $F_\mu$-links by studying the category
$\Bl(\Q)$ directly so we desire a proof of
Theorem~\ref{Bl_endo_rings_fin_diml} which avoids any choice
of Seifert surface or Seifert module. One such proof is due to
Lewin~\cite{Lew69}. In a
subsequent paper we shall give a proof which applies
when $k[F_\mu]$ is replaced by a wider class of rings.

Before leaving the subject of Morita equivalence we pause to note the
following
``naturality'' statement which we will need in
Section~\ref{section:Proof_that_B_induces_Witt_isomorphism} to prove
Theorem~\ref{theorem_relating_invariants:vague}. The reader may refer
to equation~(\ref{how2_compose_dpres_functors}) for the definition of
composition for duality-preserving functors.
\begin{proposition}\label{naturality_of_h_morita_eq}
Suppose $(G,\Psi,\eta'):\cy{C}\to\cy{D}$ is a duality-preserving functor
between hermitian categories and $b:M\to M^*$ is an $\eta$-hermitian
form in~$\cy{C}$.
Let $E=\End_\cy{C}(M)$ and let $E'=\End_\cy{D}G(M)$. The following diagram of
duality-preserving functors commutes up to natural isomorphism:
\begin{equation*} 
\xymatrix@C=11ex@R=3ex{
{\cy{C}|_M} \ar[r]^{(G\,,\,\Psi\,,\,\eta')}
\ar[d]_{(\Hom(M,\functor)\, ,\,\Omega^b\, ,\,\eta)} & {\cy{D}|_{G(M)}}
\ar[d]^{(\Hom(G(M),\functor)\, ,\, \Omega^{\eta'\Psi G(b)}\, ,\,\eta\eta')} \\
{E\proj} \ar[r]_{(E'\otimes_E\functor\,,\, \Pi\,,\,1)} & {E'\proj}.
}
\end{equation*}
\end{proposition}
\begin{proof}
See Section~\ref{section:natuality_of_h_Morita_eq} of
Appendix~\ref{section:naturality_of_B}.
\end{proof}
\begin{corollary}\label{Witt_naturality_of_h_morita_eq}
If $G$ is exact then the following square also commutes:
\begin{equation*}
\xymatrix@C=8ex@R=3ex{
W^\zeta(\cy{C}|_M) \ar[r]^{G} \ar[d]^\cong_{\Theta_{M,b}} &
W^{\zeta\eta'}(\cy{D}|_{G(M)}) \ar[d]_\cong^{\Theta_{G(M),\eta'\Psi
G(b)}}\\
W^{\zeta\eta}(E) \ar[r]_G & W^{\zeta\eta}(E').
}
\end{equation*}
\end{corollary}
\subsection{Step 3: Invariants}
\label{table_of_invariants}
Equation~(\ref{devissage_applied_to_flk(Q)}) leads one to consider
invariants to distinguish Witt classes of forms over division algebras
$E$ of finite dimension over $\Q$. Such division algebras are
well understood (see Albert~\cite[p149,p161]{Alb39},
Scharlau~\cite{Scha85} or our earlier summary in~\cite[\S11.1,11.2]{She03mem}).

One considers five distinct classes of division algebras with
involution. Firstly a division algebra $E$ may be commutative or
non-commutative. Secondly, if $I$ is an involution on $E$ let
$\Fix(I)=\{a\in E\,|\,I(a)=a\}$. The involution is said to be ``of the
first kind'' if $\Fix(I)$ contains the center $K=Z(E)$ of $E$. 
Otherwise, the involution is ``of the second kind''. 
Finally, one of these four classes is further partitioned. A
non-commutative division algebra with involution of the first kind is
necessarily a quaternion algebra, with presentation
$K\langle i,j\,|\,
i^2=a,j^2=b,ij=-ji\rangle$ for some number field $K$ and some elements
$a,b\in K$. If $\Fix(I)=Z(E)=K$ then the involution is called
``standard''. On the other hand if $\Fix(I)$ strictly larger than
$K$ then the involution is called ``non-standard''.

Table~(\ref{table_of_Witt_invariants}) below lists sufficient
invariants to distinguish the Witt classes of forms over each class of
division algebras with involution. The symbol $m\; (2)$ denotes
dimension modulo $2$. The letter $\sigma$ signifies all signature invariants
(if any) each of which takes values in $\Z$. The discriminant $\Delta$
is the determinant with a possible sign adjustment and
takes values in the group of ``square classes''  
\begin{equation*}
\frac{\Fix(I)\cap K}{\{aI(a)\mid a\in K^\bullet\}}
\end{equation*}
where $K$ is the center of $E$ and $K^\bullet=K\setminus0$.
\begin{equation}\label{table_of_Witt_invariants}
\begin{array}{|c|l|l|l|}
\hline
{\it Kind} & {\it Commutative?} & {\it Involution} & {\it Invariants} \\
\hline
1st & \mbox{Yes} & \mbox{(Trivial)} & \mbox{$m\; (2)$,
$\sigma$, $\Delta$, $c$} \\
\hline
1st & \mbox{No} & \mbox{Standard} & \mbox{$m\; (2)$,
$\sigma$}  \\
\hline
1st & \mbox{No} & \mbox{Non-standard} & 
\mbox{$m\; (2)$, $\sigma$, $\Delta$, $\theta$} \\
\hline
2nd & \mbox{Yes} & \mbox{(Non-trivial)} & \mbox{$m\; 
(2)$, $\sigma$, $\Delta$} \\
\hline
2nd & \mbox{No} & \mbox{(Non-trivial)} & \mbox{$m\; (2)$, 
$\sigma$, $\Delta$} \\
\hline
\end{array}
\end{equation}
Two symbols in the table have not yet been mentioned. The Hasse-Witt
invariant, $c$, which appears in the first row takes values in a
direct sum of copies of $\{1,-1\}$, one copy for each prime of the number
field $K$. Finally, if $E$ is a quaternion algebra with non-standard
involution of the first kind then the local-global principle fails and one
requires a secondary invariant such as the Lewis $\theta$ which is defined if
all the other invariants vanish. 
The value group for $\theta$ is the quotient $\{1,-1\}^S\slash\sim$ where $S$
is the set of primes $\mathfrak{p}$ of $K$ such that the completion
$E_\mathfrak{p}$ is a division algebra and the relation $\sim$ identifies
each element $\{\epsilon_\mathfrak{p}\}_{\mathfrak{p}\in S}$ with its antipode  $\{-\epsilon_\mathfrak{p}\}_{\mathfrak{p}\in S}$.
\section{Seifert forms}\label{section:seifert}
In this section we describe algebraic structures arising in the study
of a Seifert surface of an $F_\mu$-link. We define in
Section~\ref{section:seifert_modules} a category 
$\Sei_\infty(A)$ of ``Seifert modules'' and a full subcategory
$\Sei(A)\subset \Sei_\infty(A)$ in which the objects are
finitely generated and projective as $A$-modules. The category 
$\Sei(A)$ was denoted $(P_\mu\mbox{--}A)\proj$ in~\cite{She03mem}.

In Section~\ref{blanchfield_module_has_seifert_structure} we prove that
every $F_\mu$-link module is also a Seifert module in a canonical 
way (cf Farber~\cite{Far91}) and obtain a
``forgetful'' functor 
\begin{equation*}
U:\Bl_\infty(A)\to\Sei_\infty(A).
\end{equation*}
The image of $\Bl(A)$ is usually not
contained in $\Sei(A)$ which explains our motivation for introducing
$\Sei_\infty(A)$.
We construct later
(Section~\ref{section:covering_construction}) a
functor $B$ from $\Sei(A)$ to $\Bl(A)$ which extends in an obvious way
to a functor $B:\Sei_\infty(A)\to \Bl_\infty(A)$ and is left adjoint to~$U$.

In Section~\ref{section:cobordism_invariants_seifert} we put
hermitian structure on $\Sei(A)$.
We will see in
Section~\ref{section:Proof_that_B_induces_Witt_isomorphism} below that 
the functor $B$ induces an isomorphism of Witt groups 
$W^\zeta(\Sei(A))\to W^{-\zeta}(\Bl(A))$ when $A$ is semi-simple
and Artinian and, in particular, when $A=\Q$. 
\subsection{Seifert modules}\label{section:seifert_modules}
Suppose $V$ is a finite-dimensional vector space over a field $k$ and
$\alpha:V\to V$ is an endomorphism. A time-honoured technique in
linear algebra regards the pair $(V,\alpha)$ as a module over a
polynomial ring $k[s]$ in
which the action of~$s$ on~$V$ is given by~$\alpha$. Equivalently,
$(V,\alpha)$ is a representation of $\Z[s]$ in the category of
finite-dimensional vector spaces over $k$. We shall use the words
``module'' and ``representation'' interchangeably.

Given a Seifert surface $U^{n+1}\subset S^{n+2}$ for a knot
$S^n\subset S^{n+2}$, small translations in the
directions normal to $U$ induce homomorphisms
\begin{equation*}
f^+,f^-:H_i(U)\to H_i(S^{n+2}\setminus U).
\end{equation*}
Using Alexander duality
one finds that $f^+-f^-$ is an isomorphism for $i\neq0,n+1$, so
$H_i(U)$ is endowed with an endomorphism
$(f^+-f^-)^{-1}f^+$ and may therefore be regarded as a representation
of a polynomial ring $\Z[s]$. The homology of a Seifert surface for a
$\mu$-component boundary link has, in addition to the endomorphism
$(f^+-f^-)^{-1}f^+$, a system of $\mu$ orthogonal idempotents
which express the component structure of the Seifert surface. Following
Farber~\cite{Far92} (see also~\cite{She03mem}) we regard $H_i(U)$
as a representation of the ring
\begin{align*}
P_\mu~ &=~\Z\left\langle s,\pi_1,\cdots,\pi_\mu
\biggm| \pi_i^2=\pi_i;\; \pi_i\pi_j=0 \text{\
for \mbox{$i\neq j$}};\sum_{i=1}^\mu
\pi_i=1 \right\rangle. \\ 
&\cong~\Z[s] *_\Z \left(\prod_\mu \Z\right).
\end{align*}
$P_\mu$ is also the path ring of a quiver which we illustrate in the
case $\mu=2$:
\begin{equation*}
\xymatrix{
\bullet \ar@/^/[r] \ar@(ul,dl)[] & \bullet \ar@/^/[l] \ar@(ur,dr)[]
}
\end{equation*}
\begin{notation}\label{Sei_notation}
Let $\Sei_\infty(A)$ denote the category of representations of $P_\mu$ by
$A$-modules. A module in $\Sei_\infty(A)$ is a pair
$(V,\rho)$ where $V$ is an $A$-module and $\rho:P_\mu\to \End_A(V)$ is
a ring homomorphism. 

Let $\Sei(A)$ denote the category of representations of $P_\mu$ by
finitely generated projective $A$-modules. In other words, $\Sei(A)$
is the full subcategory of pairs $(V,\rho)$ such that $V$ is finitely
generated and projective.
\end{notation}
We sometimes omit $\rho$, confusing an element of $P_\mu$ with its
image in $\End_A(V)$.
\subsection{Seifert structure on $F_\mu$-link
modules}\label{blanchfield_module_has_seifert_structure} 
Every $F_\mu$-link module has a canonical Seifert module structure which we
describe next. In fact $\Bl_\infty(A)$ can be regarded as a
full subcategory of $\Sei_\infty(A)$. Note however that a module in
$\Bl(A)$ is in general neither finitely generated nor projective as an
$A$-module (e.g.~see Lemma~\ref{Sato_condition} above) and therefore
does not lie in $\Sei(A)$.

If $M\in\Bl_\infty(A)$ then Lemma~\ref{Sato_condition} implies that
the $A$-module map 
\begin{align*}
\gamma:M^{\oplus \mu}&\to M \\
(m_1,\cdots,m_\mu) &\mapsto \sum_{i=1}^\mu (1-z_i)m_i
\end{align*}
is an isomorphism.
Let $p_i$ denote the
projection of $M^{\oplus \mu}$ onto its $i$th component and let
\begin{align*}
\omega:M^{\oplus\mu}&\to M \\
(m_1,\cdots,m_\mu)&\mapsto \sum_{i=1}^\mu m_i
\end{align*} 
denote addition. Define $\rho:P_\mu\to\End_kM$ by  
\begin{equation}
\label{Seifert_structure_on_F-link_module}
{\begin{aligned}
&\rho(\pi_i)=\gamma p_i\gamma^{-1} \\
&\rho(s)= \omega\gamma^{-1}.
\end{aligned}}
\end{equation}
We denote by $U(M)$ the $A$-module $M$ with the Seifert module structure
$\rho$. As remarked in the introduction, in the case $\mu=1$ the
Seifert module structure $\rho:P_\mu\to \End_A(M)$ can be described
more simply by the equation 
\begin{equation*}
\rho(s)=(1-z)^{-1}.
\end{equation*}
The following lemma says that $U:\Bl_\infty(A)\to \Sei_\infty(A)$ is a
full and faithful functor, so $\Bl_\infty(A)$ can be regarded as a
full subcategory of $\Sei_\infty(A)$.
\begin{lemma}\label{Bl_infty_is_full_subcat_of_Sei_infty}
An $A$-module morphism $f:M\to M'$ between $F_\mu$-link modules $M$ and
$M'$ is an $A[F_\mu]$-module morphism if and only if
$f$ is a morphism of Seifert modules. In other words,
$f\in\Hom_{\Sei_\infty(A)}(U(M),U(M'))$ if and only if
$f\in\Hom_{\Bl_\infty(A)}(M,M')$. 
\end{lemma}
\begin{proof}
If $f:M\to M'$ is an $A[F_\mu]$-module morphism then the
 diagram
\begin{equation}\label{f_commute_with_gamma}
\xymatrix@C=3ex@R=2ex{
M^{\oplus\mu}\ar[r]^{f^{\oplus\mu}}\ar[d]_\gamma & {M'}^{\oplus\mu}\ar[d]^\gamma \\
M\ar[r]_f & M'}
\end{equation}
is commutative. Conversely, if the
 diagram~(\ref{f_commute_with_gamma}) commutes then the equation
 $f((1-z_i)x)=(1-z_i)f(x)$ holds for each $x\in M$ and $i=1,\cdots,\mu$, so
 $f(z_ix)=z_if(x)$ for each $i$ and hence $f$ is an $A[F_\mu]$-module
 morphism. 

It remains to show that~(\ref{f_commute_with_gamma}) commutes if and
only if $f$ is a Seifert morphism. If~(\ref{f_commute_with_gamma})
commutes then 
\begin{equation*}
f(\gamma p_i\gamma^{-1}x)=\gamma p_i\gamma^{-1}f(x)
\quad \mbox{and}\quad f(\omega\gamma^{-1}x)=\omega\gamma^{-1}f(x)
\end{equation*}
for each $x\in M$ so $f$ is a Seifert morphism. Conversely, suppose
$f$ is a Seifert morphism. Now 
\begin{equation*}
\gamma^{-1}=
\left(\begin{matrix}
\omega p_1\gamma^{-1} \\
\omega p_2\gamma^{-1} \\
\vdots \\
\omega p_\mu\gamma^{-1}
\end{matrix}\right)
=
\left(\begin{matrix}
(\omega\gamma^{-1})(\gamma p_1\gamma^{-1}) \\
(\omega\gamma^{-1})(\gamma p_2\gamma^{-1}) \\
\vdots \\
(\omega\gamma^{-1})(\gamma p_\mu\gamma^{-1})
\end{matrix}\right)
: M\to M^{\oplus\mu}
\end{equation*}
so the diagram~(\ref{f_commute_with_gamma}) commutes.
\end{proof}

\subsection{Seifert forms}\label{section:cobordism_invariants_seifert}
Let us make $\Sei(A)$ a hermitian category, assuming $A$ is a ring with
involution. Recall that if $V$ is a
finitely generated projective left $A$-module then $V^*=\Hom(V,A)$ is a left
$A$-module with $(a.\theta)(x)=\theta(x)\overline{a}$ for all $a\in
A$, $\theta\in V^*$ and $x\in V$.

If $(V,\rho)$ is an object
in $\Sei(A)$ define $(V,\rho)^*=(V^*,\rho^*)$ where
\begin{equation*}
\rho^*(\pi_i)=\rho(\pi_i)^*:V^*\to V^*
\quad\mbox{and}\quad \rho^*(s)=1-\rho(s)^*:V^*\to
V^*. 
\end{equation*}
It is easy to see that if $f:V\to V'$ is a morphism in $\Sei(A)$ then 
the dual $f^*:(V')^*\to V^*$ again lies in $\Sei(A)$. 
Equivalently, if one gives the ring $P_\mu$ the involution defined by
$\overline{s}=1-s$ and $\overline{\pi_i}=\pi_i$ for each $i$ then
$\rho^*:P_\mu\to \End(V^*)$ is given by
$\rho^*(r)(\theta)(x)=\theta(\rho(\overline{r})x)$ for all $x\in V$,
$\theta\in V^*$ and $r\in P_\mu$.

If $V=H_q(U^{2q})/\text{torsion}$ then the intersection form
$\phi:V\to V^*$ is a morphism in $\Sei(\Z)$. In other words, the
intersection form respects the projections $\pi_i$ and respects the
endomorphism $(f^+-f^-)^{-1}f^+$ (see section~\ref{section:seifert_modules})
in the sense that
\begin{equation*}
\phi((f^+-f^-)^{-1}f^+x)(y)=\phi(x)((1-(f^+-f^-)^{-1}f^+)y).
\end{equation*}
for all $x,y\in V$.
Furthermore, $\phi$ is an isomorphism by Poincar\'e duality and is
$(-1)^q$-hermitian. This form $\phi$ will be called the Seifert form
associated to a Seifert surface. 
Kervaire~\cite[p94]{Ker71} or Lemma~3.31 of~\cite{She03mem} shows that
by an elementary change of variables, this form $\phi$ is equivalent
to the Seifert matrix of linking numbers more commonly encountered in
knot theory.

By Ko~\cite{Ko87} and Mio~\cite{Mio87} (see also Lemma~3.31
of~\cite{She03mem}), the association of this form $\phi$ to a Seifert surface
induces the isomorphism 
\begin{equation*}
C_{2q-1}(F_\mu)~\cong~W^{(-1)^q}(\Sei(\Z))\quad\quad (q\geq3).
\end{equation*}
mentioned in the introduction. Every $(-1)^q$-hermitian Seifert form
is associated to some $2q$-dimensional Seifert surface. Although there
are many possible Seifert surfaces for a given $F_\mu$-link, all are
cobordant and the corresponding Seifert forms lie in the same Witt class. 

In~\cite{She03mem} the author applied to $\Sei(\Q)$ the steps 1-3.~described in
Section~\ref{section:overview_of_BD_form_invariants} obtaining
explicit invariants to distinguish $F_\mu$-link cobordism
classes. Although the Blanchfield-Duval form is more intrinsic, 
the advantage of the Seifert form is that it is easier to compute
the numerical invariants. For illustration, we treat a worked example:
\begin{example}\label{worked_example_Seifert}
Setting $\mu=2$, consider the Seifert module $V=\Z^6$ with the
endomorphism $s$ and $(-1)$-hermitian form $\phi$ given by
\begin{equation*}
s=\left(
\begin{array}{cccc|cc}
1 & 0 & \phantom{-}1 & \phantom{-}0 & \phantom{-}0 & \phantom{-}0 \\
0 & 1 & -1           &           -1 &           -1 & \phantom{-}0 \\ 
0 & 1 & \phantom{-}0 & \phantom{-}0 & \phantom{-}0 & -1 \\
0 & 0 & \phantom{-}0 & \phantom{-}0 & \phantom{-}0 & \phantom{-}0 \\ \hline
0 & 1 & \phantom{-}0 & \phantom{-}0 & \phantom{-}1 & -1 \\
0 & 0 & \phantom{-}1 & \phantom{-}0 & \phantom{-}1 & \phantom{-}0
\end{array}\right)\quad \mbox{and}\quad
\phi=
\left(
\begin{array}{cccc|cc}
\phantom{-}0 & 0 & \phantom{-}0 & 1 & 0 & \phantom{-}0 \\
\phantom{-}0 & 0 &           -1 & 0 & 0 & \phantom{-}0 \\
\phantom{-}0 & 1 & \phantom{-}0 & 0 & 0 & \phantom{-}0 \\
          -1 & 0 & \phantom{-}0 & 0 & 0 & \phantom{-}0 \\ \hline
\phantom{-}0 & 0 & \phantom{-}0 & 0 & 0 & -1 \\
\phantom{-}0 & 0 & \phantom{-}0 & 0 & 1 & \phantom{-}0
\end{array}
\right)
\end{equation*}
The horizontal and vertical lines indicate the component structure of
the Seifert form. In other words, $\pi_1$ projects onto the
span of the first four basis elements while $\pi_2$ projects onto the
span of the last two. 
The corresponding Seifert matrix of linking numbers is
\begin{equation*}
\phi s=\left(
\begin{array}{cccc|cc}
\phantom{-}0 & \phantom{-}0 & \phantom{-}0 & \phantom{-}0 & \phantom{-}0 & \phantom{-}0 \\
\phantom{-}0 &           -1 & \phantom{-}0 & \phantom{-}0 & \phantom{-}0 & \phantom{-}1 \\
\phantom{-}0 & \phantom{-}1 &           -1 &           -1 &           -1 & \phantom{-}0 \\
          -1 & \phantom{-}0 &           -1 & \phantom{-}0 &
\phantom{-}0 & \phantom{-}0 \\ \hline
\phantom{-}0 & \phantom{-}0 &           -1 & \phantom{-}0 &           -1 & \phantom{-}0 \\
\phantom{-}0 & \phantom{-}1 & \phantom{-}0 & \phantom{-}0 & \phantom{-}1 & -1
\end{array}
\right)
\end{equation*}
but we shall work directly with $s$ and $\phi$.

The first step is to pass from $\Z$ to $\Q$ so we regard the entries
in the matrices as rational numbers. Devissage is next; let
$e_1,\cdots,e_6$ denote the standard basis of $\Q^6$. 
Now $\Q e_1$ is $s$-invariant and
$\phi(e_1)(e_1)=0$ so our Seifert form is Witt-equivalent to the
induced form on $e_1^\perp/\Q e_1 \cong \Q\{e_2,e_3,e_5,e_6\}$. We
have reduced $s$ and $\phi$ to
\begin{equation*}
s'=\left(
\begin{array}{cc|cc}
1 & -1           &            -1 & \phantom{-}0 \\ 
1 & \phantom{-}0 & \phantom{-}0 & -1 \\ \hline
1 & \phantom{-}0 & \phantom{-}1 & -1 \\
0 & \phantom{-}1 & \phantom{-}1 & \phantom{-}0
\end{array}\right)\quad \mbox{and}\quad
\phi'=
\left(
\begin{array}{cc|cc}
0 &     -1      & 0 & \phantom{-}0 \\
1 & \phantom{-}0& 0 & \phantom{-}0 \\ \hline
0 & \phantom{-}0& 0 & -1 \\
0 & \phantom{-}0& 1 & \phantom{-}0
\end{array}
\right)
\end{equation*}
The two-dimensional representation (over $\Q$) of $\Z[s]$ given by the
matrix $r=\left(\begin{matrix} 
1 & -1 \\
1 & 0
\end{matrix}\right)$ 
is simple (=irreducible) since there do not exist eigenvalues in $\Q$. It
follows that the Seifert module $V'=\Q^4$ with the action $s'$ and
$\pi_1,\pi_2$ shown is simple. The devissage
process is therefore complete.

Turning to the Morita equivalence step, the
endomorphism ring of this module $V'$ has $\Q$-basis consisting of the
identity and  
$\left(\begin{array}{c|c}
r & 0 \\ \hline 
0 & r
\end{array}\right)$.
The minimum polynomial of $r$ is $x^2-x+1$ so $\End_{\Sei(\Q)}(V')$ is
isomorphic to $\Q(\sqrt{-3})$. 
We may choose $b=-\phi':V'\to (V')^*$. It is easy to verify that the
involution $f\mapsto b^{-1}f^* b$ is not the identity map so it must
send $\sqrt{-3}$ to $-\sqrt{-3}$. Morita equivalence sends the form
$\phi':V'\to {V'}^*$ to the composite
\begin{multline*}
\Hom_{\Sei(\Q)}(V',V')\xrightarrow{\phi'}
\Hom_{\Sei(\Q)}(V',{V'}^*) \\
\xrightarrow{\Omega^b_{V'}}\Hom_{\Q(\sqrt{-3})}(\Hom(V',{V'}^*),\Hom(V',V'))
\end{multline*}
 which is given by
\begin{equation*}
\Omega^b_{V'}\phi'_*(\alpha)(\beta)=-b^{-1}\beta^*\phi'\alpha=\overline{\beta}\alpha.  
\end{equation*}
for $\alpha,\beta\in\Hom(V',V')$. 
This form may be written $\langle 1\rangle$ as a form over $\Q(\sqrt{-3})$.

Reading the fourth line of the table~(\ref{table_of_Witt_invariants}),
an element of the Witt group $W(\Q(\sqrt{-3}))$ for non-trivial
involution is determined by signatures and discriminant (and rank
modulo $2$ if there are no signatures). Up to complex conjugation 
there is precisely one embedding of $\Q(\sqrt{-3})$ in $\C$ (with
the complex conjugate involution), so there is in fact one signature,
which takes value $1\in \Z$ with our choice of $b$. The
discriminant is 
\begin{equation*}
1\in \frac{\Q\setminus0}{\Q(\sqrt{-3})\overline{\Q(\sqrt{-3})}} 
= \frac{\Q\setminus0}{\{a^2+3b^2\mid a,b\in\Q\}}.
\end{equation*} 
\end{example}

\section{The Covering construction}\label{section:covering_construction}
In this section we introduce 
a functor $B:\Sei_\infty(A)\to\Bl_\infty(A)$
which is the algebraic analogue of the geometric construction of the
free cover of an $F_\mu$-link complement from a Seifert surface.
(illustrated in Figure~1 on page~\pageref{picture}). 
The
restriction of $B$ to $\Sei(A)$ takes values in $\Bl(A)$ and extends
to a duality-preserving functor
\begin{equation*}
(B,\Phi,-1):\Sei(A)\to\Bl(A)
\end{equation*}
which is natural in $A$ (see
Propositions~\ref{extend_B_to_duality_pres_functor}
and~\ref{B_respects_coefficient_change}).

We show that $B:\Sei_\infty(A)\to\Bl_\infty(A)$ is left
adjoint to the full and faithful functor
$U:\Bl_\infty(A)\to\Sei_\infty(A)$; in other words, among
functors $\Sei_\infty(A)\to\Bl_\infty(A)$, the geometrically motivated
functor $B$ satisfies a universal property with respect to $U$ (see Definition~\ref{define_left_adjoint}).
In particular there is a natural transformation $\theta_V: V\to UB(V)$ for
$V\in\Sei(A)$ and a natural isomorphism $\psi_M:BU(M)\to M$ for
$M\in\Bl_\infty(A)$. The reader is referred to chapter 3 of
Borceux~\cite{Bor94_vol1} or chapter IV of Mac Lane~\cite{Mac71} for a
detailed treatment of adjoint functors. 

We use the adjunction in Section~\ref{localization_of_cat} to show
that the covering construction $B:\Sei_\infty(A)\to\Bl_\infty(A)$ is
equivalent to a universal localization $\Sei_\infty(A)\to
\Sei_\infty(A)/\Prim_\infty(A)$ of categories. We describe the
structure of the ``primitive'' modules 
$V\in\Prim_\infty(A)$ in Section~\ref{section:structure_of_primitives}
and outline a construction of the quotient category in
Section~\ref{section:construction_of_quotient}.

We show in Sections~\ref{localization_of_cat}
and~\ref{section:duality_on_quotient_arbitrary_A} that
$B:\Sei(A)\to\Bl(A)$ is equivalent to a localization
$\Sei(A)\to\Sei(A)/\Prim_\infty(A)$ of hermitian 
categories. In the case where $A$ is a semi-simple Artinian ring we
simplify the descriptions of the quotient and 
primitive modules in Section~\ref{section:global_dim_zero}.
\subsection{Definition}
To simplify notation in this section and
Section~\ref{section:adjunction}, we suppress the symbol $\rho$ which
appears in the definition of a Seifert module $(V,\rho)$, identifying
an element $r\in P_\mu$ with $\rho(r)\in\End_A(V)$. We shall extend
Seifert structure from an $A$-module $V$ to the induced module
$A[F_\mu]\otimes_A V$ by $s(\alpha\otimes
v)=\alpha\otimes s(v)$ and $\pi_i(\alpha\otimes v)=\alpha\otimes
\pi_i(v)$ for $\alpha\in A[F_\mu]$.

Recall that $z_1,\cdots,z_\mu$ are distinguished generators of
$F_\mu$; let us now write $z=\sum z_i\pi_i$. 
\begin{definition}\label{define_covering_construction}
If $V$ is a module in $\Sei_\infty(A)$ let
\begin{equation*}\label{covering_constrution}
B(V)= \Coker\left((1 - s(1-z)): V[F_\mu]\to V[F_\mu]\right).
\end{equation*}
\end{definition}
Since $\epsilon(1-s(1-z))=1$, it is clear that $B(V)$ lies in
  $\Bl_\infty(A)$. To achieve more symmetric notation we write $z_i=y_i^2$ and
  deduce $z=y^2$ where $y=\sum y_i\pi_i$. We write $F_\mu=F_\mu(y^2)$ when we
  wish to indicate that elements of $F_\mu$ are to be written as words in the
  symbols $y_i^{\pm2}$. The $A[F_\mu]$-module $\displaystyle{\bigoplus_{i=1}^\mu
  \pi_iV[F_\mu(y^2)y_i]}$ is isomorphic to $V[F_\mu]$ and will be written
  $V[F_\mu(y^2)y]$ for brevity. Now we have
\begin{equation}
\begin{aligned}
B(V)&\cong \Coker\left((1 - s(1-y^2))y^{-1}: V[F_\mu(y^2)y]\to
    V[F_\mu(y^2)]\right) \\ 
&= \Coker\left((1-s)y^{-1} + sy: V[F_\mu(y^2)y]\to V[F_\mu(y^2)]\right).
\end{aligned}
\end{equation}
In detail, 
\begin{equation*}
(1-s)y^{-1} + sy=\sum_{i=1}^\mu (1-s)\pi_iy_i^{-1} + s\pi_iy_i
:vwy_i \mapsto (1-s)(v)w + s(v)wy_i^2
\end{equation*}
for $v\in \pi_i(V)$, $w\in F_\mu(y^2)$ and $i\in\{1,\cdots,\mu\}$. A
morphism $f:V\to V'$ induces a commutative diagram
\begin{equation}
\label{B_of_morphism}
\begin{gathered}
\xymatrix@R=3ex{ 
0 \ar[r] & V[F_\mu] \ar[r]^\sigma\ar[d]_f & V[F_\mu] \ar[r]^q\ar[d]_f & B(V)
\ar[r]\ar@{.>}[d]^{B(f)} & 0 \\  
0 \ar[r] & V'[F_\mu] \ar[r]^{\sigma'} & V'[F_\mu] \ar[r]^{q'} & B(V') \ar[r] & 0}
\end{gathered}
\end{equation}
where $\sigma=1-s(1-z)$ and hence induces an $A[F_\mu]$-module map
$B(f)$ as shown. 
\begin{lemma}\label{exactness_of_B}
The functor $B$ is exact. In other words, if $V\to V'\to V''$ is an
exact sequence in $\Sei_\infty(A)$ then the sequence $B(V)\to
B(V')\to B(V'')$ induced is also exact.
\end{lemma}
\begin{proof}
It suffices to show that $B$ preserves short exact sequences. Suppose
$0\to V\to V'\to V''\to 0$ is exact. There is a commutative diagram
\begin{equation*}
\begin{gathered}\xymatrix@=2ex{ 
& 0 \ar[d] & 0 \ar[d] & 0 \ar[d] & \\
0 \ar[r] & V[F_\mu] \ar[r]\ar[d] & V[F_\mu] \ar[r]\ar[d] & B(V)
\ar[r]\ar[d] & 0 \\  
0 \ar[r] & V'[F_\mu] \ar[r]\ar[d] & V'[F_\mu]\ar[d] \ar[r] & B(V')\ar[d] \ar[r] & 0 \\
0 \ar[r] & V''[F_\mu] \ar[r]\ar[d] & V''[F_\mu] \ar[r]\ar[d] & B(V'')
\ar[r]\ar[d] & 0 \\  
& 0 & 0 & 0 & 
}
\end{gathered}
\end{equation*}
in which the rows and the two left-most columns are exact. It follows
that the right-hand column is also exact.
\end{proof}
The category $\Sei_\infty(A)$ has limits and colimits. For example, the
coproduct of a system of modules is the direct sum. Since $B$
is exact and respects arbitrary direct sums $B$ respects all colimits:
\begin{lemma}\label{B_is_cocontinuous}
The functor $B$ is cocontinuous.\hspace*{\fill}\qed
\end{lemma}
In particular if $V$ is a direct limit $V=\varinjlim V_i$ then
$B(V)=\varinjlim B(V_i)$. On the other hand $B$ does not respect
infinite limits. For example one finds $B(\prod V_i)\ncong 
\prod(B(V_i))$ because $(\prod V_i)[F_\mu]\ncong \prod (V_i[F_\mu])$. However,
$B$ does respect finite limits.

The idea behind the proof of the following proposition is due to A.Ranicki.
\begin{proposition}~\label{extend_B_to_duality_pres_functor}
The functor $B$ extends to a duality-preserving functor
$(B,\Phi,-1):\Sei(A)\to\Bl(A)$.
\end{proposition}
\begin{proof}
For each finitely generated projective $A$-module $V$ there is
a natural isomorphism $\Pi_V:V^*[F_\mu]\to (V[F_\mu])^*$ by
Example~\ref{nat_iso_in_coeff_change_dpres_funct}.
Naturality asserts that for each morphism $\alpha:V\to W$ the diagram
\begin{equation*}
\xymatrix@C=4ex@R=3ex{
W^*[F_\mu]\ar[r]^{\alpha^*}\ar[d]_{\Pi_W} & V^*[F_\mu] \ar[d]^{\Pi_V} \\
(W[F_\mu])^*\ar[r]^{\alpha^*} & (V[F_\mu])^*
}
\end{equation*}
is commutative. Moreover one can check commutativity of
\begin{equation*}
\xymatrix@C=6ex@R=3ex{
V^*[F_\mu]\ar[r]^{z_i^{-1}}\ar[d]_{\Pi_V} & V^*[F_\mu] \ar[d]^{\Pi_V} \\
(V[F_\mu])^*\ar[r]^{z_i^*} & (V[F_\mu])^*
}
\end{equation*}
where, as usual, $z_i:V[F_\mu]\to V[F_\mu]$ and $z_i:V^*[F_\mu]\to
V^*[F_\mu]$ denote multiplication on the right by $z_i$.

Now if $V\in\Sei(A)$ then there is a commutative diagram:
\begin{equation}\label{diagram_defining_Phi}
\begin{gathered}
\xymatrix@C=6ex@R=3ex{
0 \ar[r] & V^*[F_\mu]\ar[d]_{-\Pi_V(1-z)}\ar[r]^{\sigma(V^*)} & V^*[F_\mu] \ar[r]\ar[d]^{\Pi_V(1-z^{-1})} & B(V^*)\ar@{.>}[d]^{\Phi_V}
\ar[r] & 0 \\
0 \ar[r] & (V[F_\mu])^*\ar[r]^{\sigma(V)^*} & (V[F_\mu])^* \ar[r] & B(V)^\wedge
\ar[r] & 0
}
\end{gathered}
\end{equation}
where $\sigma(V^*)=1-(1-s^*)(1-z)$ and $\sigma(V)^*=(1-s(1-z))^*$. By
definition, $\Phi_V:B(V^*)\to B(V)^\wedge$ is the induced
morphism. Plainly $\Phi_V$ is a natural transformation.

The duality-preserving functor $\Pi$ has the
property 
\begin{equation*}
\Pi_V^*i_{V[F_\mu]} = \Pi_{V^*}i_V:V[F_\mu]\to (V^*[F_\mu])^*
\end{equation*}
(indeed, this equation features in
Definition~\ref{define_dpres_funct}). The equations
\begin{align*}
(\Pi_V(1-z^{-1}))^* i_{V[F_\mu]} = \Pi_{V^*}(1-z) i_V \\
(\Pi_V(1-z))^* i_{V[F_\mu]} = \Pi_{V^*}(1-z^{-1}) i_V
\end{align*}
imply that
\begin{equation*}
\Phi_V^\wedge i_{B(V)}=-\Phi_{V^*}B(i_V).
\end{equation*}
To show that $(B,\Phi,-1)$ is a
duality-preserving functor it remains to check that $\Phi$ is an
isomorphism. There is a commutative diagram
\begin{equation*}
\xymatrix@C=6ex@R=3ex{
0 \ar[r] & (V[F_\mu])^*\ar[d]_{-(1-s^*)\Pi_V^{-1}}\ar[r]^{\sigma(V)^*} & (V[F_\mu])^* \ar[r]\ar[d]^{-(1-s^*)z\Pi_V^{-1}} & B(V)^\wedge\ar@{.>}[d]
\ar[r] & 0 \\
0 \ar[r] & V^*[F_\mu]\ar[r]^{\sigma(V^*)} & V^*[F_\mu] \ar[r] & B(V^*)
\ar[r] & 0
}
\end{equation*}
and the composite morphisms of chain complexes
\begin{align*}
&\xymatrix@C=4ex@R=3ex{
V^*[F_\mu]\ar[d]_{(1-s^*)(1-z)}\ar[r]^{\sigma(V^*)} & V^*[F_\mu]\ar@{.>}[dl]_1
\ar[d]^{(1-s^*)(1-z)} \\
V^*[F_\mu]\ar[r]_{\sigma(V^*)} & V^*[F_\mu]} \\
\intertext{and}
&\xymatrix@C=4ex@R=3ex{(V[F_\mu])^*\ar[d]_{(1-z^{-*})(1-s^*)}\ar[r]^{\sigma(V)^*} & (V[F_\mu])^*\ar@{.>}[dl]_{z^{-*}}
\ar[d]^{-(1-z^*)(1-s^*)z^{-*}} 
\\
(V[F_\mu])^*\ar[r]_{\sigma(V)^*} & (V[F_\mu])^*
}
\end{align*}
are chain homotopic to the identity by the indicated chain homotopies.
[Chain complexes are drawn horizontally and
morphisms of chain complexes are given by vertical arrows. The symbol
$z^{-*}$ is shorthand for $(z^{-1})^*$].
These composite chain maps therefore induce the identity on $B(V^*)$
and $B(V)^\wedge$ respectively so $\Phi_V$ is an isomorphism and
$(B,\Phi,-1)$ is a duality-preserving functor.
\end{proof}
Since $B$ is an exact functor we have:
\begin{corollary}
The functor $(B,\Phi,-1)$ induces a homomorphism of Witt groups
\begin{equation*}
B:W^\zeta(\Sei(A))\to W^{-\zeta}(\Bl(A)).
\end{equation*}
\end{corollary}
If $(V,\phi)$ is a $\zeta$-hermitian form in $\Sei(A)$
then the covering construction gives $B(V,\phi) = (B(V),
\Phi_VB(\phi))$ (Lemma~\ref{duality_preserving_functor_induces}) which
can be described explicitly as follows. 
The morphism $\phi:V\to V^*$ induces $\phi:V[F_\mu]\to
V^*[F_\mu]$. Let $\widetilde{\phi}:V[F_\mu]\to
\Hom_{A[F_\mu]_\Sigma}(V[F_\mu]_\Sigma,A[F_\mu]_\Sigma)$ denote
the composition of $\phi$ with 
\begin{equation*}
\Pi_V:V^*[F_\mu]\to
(V[F_\mu])^*=\Hom_{A[F_\mu]}(V[F_\mu],A[F_\mu])
\end{equation*}
(see Example~\ref{simplest_ring_change_example}) and the localization
\begin{equation*}
\Hom_{A[F_\mu]}(V[F_\mu],A[F_\mu])\to
\Hom_{A[F_\mu]_\Sigma}(V[F_\mu]_\Sigma,A[F_\mu]_\Sigma).
\end{equation*}
 If $m,m'\in B(V)$ we
may write $m=q(x)$, $m'=q(x')$ for some $x,x'\in V[F_\mu]$. It follows
from Remark~\ref{explicit_formula_for_Mwedge_pres} that in
$\displaystyle{\frac{A[F_\mu]_\Sigma}{A[F_\mu]}}$ we have
\begin{align}
\Phi_VB(\phi)(m)(m')&=\Phi_VB(\phi)(q(x))(q(x')) \notag\\
&= q'(\Pi_V(1-z^{-1})\phi(x))(q(x')) \quad\mbox{(using~(\ref{diagram_defining_Phi}))} \notag\\
&= (\id\otimes\Pi_V(1-z^{-1})\phi(x))(\id\otimes\sigma)^{-1}(1\otimes
x') \notag\\
&= \widetilde{\phi}(x)((1-z)(\id\otimes\sigma)^{-1}(1\otimes x')).
\label{formula_for_BD_form}
\end{align}
\begin{remark}\label{remark_that_B(Seifert_form)=BDform}
If $V\in\Sei(\Z)$ and
$\phi:V\to V^*$ is the Seifert form corresponding to a Seifert surface
for an $F_\mu$-link then by~(\ref{formula_for_BD_form})
$\Phi_{V}B(\phi):B(V)\to B(V)^\wedge$ is the corresponding
Blanchfield-Duval form for the $F_\mu$-link; compare
Kearton~\cite{Kea75}, Levine~\cite[Prop 14.3]{Lev77}, Cochran and
Orr~\cite[Thm4.2]{CocOrr94} and Ranicki~\cite[Defn32.7]{Ran98}. For
example, setting $r=\phi(x)$, $s=\phi(x')$, $\Gamma=z$, $\theta=\phi
s$, and $\epsilon=\zeta$ one obtains
from~(\ref{formula_for_BD_form}) the equations appearing immediately
prior to Theorem 4.2 in~\cite{CocOrr94}.
\end{remark} 

For the proof of Theorem~\ref{theorem_relating_invariants:vague} in
Section~\ref{section:Proof_that_B_induces_Witt_isomorphism} we will need
the observation that $(B,\Phi,-1)$ respects a change of coefficients
from $\Z$ to $\Q$. Let us make a more general statement. 
Recall from equation~(\ref{how2_compose_dpres_functors}) the
definition of composition for duality-preserving functors.
\begin{proposition}\label{B_respects_coefficient_change}
A ring homomorphism $A\to A'$ induces a diagram of duality-preserving functors
\begin{equation*}
\xymatrix@R=4ex@C=15ex{
{\Sei(A)}\ar[d]_{(B,\Phi,-1)} \ar[r]^{(A'\otimes_A\functor,\Pi,1)} & {\Sei(A')}\ar[d]^{(B,\Phi,-1)} \\
{\Bl(A)}\ar[r]_{(A'[F_\mu]\otimes_{A[F_\mu]}\functor,\Upsilon,1)} & {\Bl(A')}.
}
\end{equation*}
which commutes up to natural isomorphism. Consequently, there is a
commutative diagram of Witt groups
\begin{equation}\label{Witt_diagram:B_respects_coeff_change}
\begin{gathered}
\xymatrix@=2ex{
W^\zeta(\Sei(A))\ar[d]_B \ar[r] & W^\zeta(\Sei(A'))\ar[d]^B \\
W^{-\zeta}(\Bl(A))\ar[r] & W^{-\zeta}(\Bl(A')).
}
\end{gathered}
\end{equation}
\end{proposition}
\begin{proof}
See Appendix~\ref{section:naturality_of_B}.
\end{proof}
\subsection{Adjunction}\label{section:adjunction}
We leave duality structures
behind for the present and prove that the functor $B:\Sei_\infty(A)\to \Bl_\infty(A)$ is left adjoint to $U:\Bl_\infty(A)\to\Sei_\infty(A)$.

\begin{definition}\label{define_left_adjoint}
Suppose $F:\cy{C}\to\cy{D}$ is a functor. A functor
$G:\cy{D}\to\cy{C}$ is called {\it left adjoint} to $F$
if there exists a natural transformation
$\theta:\id_\cy{D}\to FG$ such that for every object $D\in\cy{D}$ the
morphism $\theta_D:D\to FG(D)$ has the following universal property:
For every morphism $d$ in $\cy{D}$ of the form $d:D\to F(C)$ there is a
unique morphism $c:G(D)\to C$ in $\cy{C}$ such that $d=F(c)\theta_D$.
\begin{equation}\label{universal_property_of_adjoint_functor}
\begin{gathered}
\xymatrix@=2ex{D\ar[rr]^d\ar[dr]_{\theta_D} & & F(C) \\
& FG(D) \ar@{.>}[ur]_{F(c)}&}
\end{gathered}
\end{equation}
\end{definition}
Let us recall a few examples: 1) The inclusion of the category of
abelian groups in the category of groups has left adjoint known as
``abelianization'' which sends a group $G$ to $G/[G,G]$. 
2) The inclusion of the category of compact
Hausdorff topological spaces in the category of (all) topological
spaces has a left adjoint known as ``Stone-\v{C}ech
compactification''. 3) Colimit constructions (e.g. direct limit or
coproduct) can be expressed via a left adjoint as follows. Suppose
$\cy{C}$ is a category, $J$ is a small category and $\cy{C}^J$
denotes the category of functors $J\to\cy{C}$. If there is a left adjoint
to the constant functor $\cy{C}\to\cy{C}^J$ then that left adjoint
sends each functor $J\to\cy{C}$ to its colimit in $\cy{C}$ (and the
colimit exists).
\begin{proposition}\label{B_is_left_adjoint_to_U}
The functor $B:\Sei_\infty(A)\to\Bl_\infty(A)$ is left adjoint to
$U:\Bl_\infty(A)\to\Sei_\infty(A)$.
\end{proposition}
The required map $\theta_V:V\to UB(V)$ is the restriction of
the map $q$ in the diagram (\ref{B_of_morphism}) above. In symbols
$\theta_V=q|: V\to UB(V)$. During 
the proof of Proposition~\ref{B_is_left_adjoint_to_U} below we show
that $\theta_V$ is a morphism of Seifert modules. It follows from the
diagram~(\ref{B_of_morphism}) that $\theta:\id\to UB$ is a natural
transformation.

Before proving Proposition~\ref{B_is_left_adjoint_to_U}, we note some
consequences:
\begin{corollary}\label{deductions_from_left_adjointness}
Let $V\in\Sei_\infty(A)$ and $M\in\Bl_\infty(A)$.
There is a natural isomorphism $\psi_M:BU(M)\to M$ and the composites 
\begin{align*}
&U(M)\xrightarrow{\theta_{U(M)}} UBU(M)\xrightarrow{U(\psi_M)} U(M) \\
&B(V)\xrightarrow{B(\theta_V)} BUB(V) \xrightarrow{\psi_{B(V)}} B(V)
\end{align*}
are identity morphisms. In particular $\theta_{U(M)}$ and
$B(\theta_V)$ are isomorphisms.
\end{corollary}
The existence of a natural transformation $\psi_M:BU(M)\to M$ follows
from Proposition~\ref{B_is_left_adjoint_to_U} alone. To prove that $\psi$
is an isomorphism one requires the additional information that $U$ is full
and faithful. We are not claiming that $\theta_V:V\to UB(V)$ is an
isomorphism. Indeed, $U$ and $B$ are not equivalences of
categories. 
\begin{proof}[Proof of Corollary~\ref{deductions_from_left_adjointness}]
Let $\psi_M:BU(M)\to M$ be the unique morphism such that
 $\id_{U(M)}=U(\psi_M)\theta_{U(M)}$. One can check
 that $\psi_M$ is a natural transformation and that
 $\psi_{B(V)}B(\theta_V)=\id_{B(V)}$ (see for example Theorem 3.1.5
 of~\cite{Bor94_vol1}). The functor $U:\Bl_\infty(A)\to\Sei_\infty(A)$ is
 full and faithful by  
 Lemma~\ref{Bl_infty_is_full_subcat_of_Sei_infty} so $\psi_M$ is an
 isomorphism (see Theorem 3.4.1 of~\cite{Bor94_vol1}). It follows that
 $\theta_{U(V)}$ and $B(\theta_V)$ are isomorphisms.
\end{proof}
\begin{proof}[Proof of Proposition~\ref{B_is_left_adjoint_to_U}]
By Definition~\ref{define_left_adjoint} there are two statements to prove:
\begin{enumerate}
\item The map $\theta_V:V\to UB(V)$ is a morphism of Seifert modules.
\item If $M\in\Bl_\infty(A)$ and $f:V\to U(M)$ is a morphism in
$\Sei_\infty(A)$ then there is a unique morphism $g:B(V)\to M$ such
that $f=U(g)\theta_V$.
\end{enumerate}
As we remarked above, it follows from the diagram~(\ref{B_of_morphism}) that
$\theta:\id\to UB$ is a natural transformation.
We shall need the following lemma which is proved shortly below:  
\begin{lemma}\label{relating_Seifert_structures}
Suppose $V\in\Sei_\infty(A)$,
$M\in\Bl_\infty(A)$ and $f:V\to M$ is an $A$-module morphism. Let
$\widetilde{f}:V[F_\mu]\to M$ denote the induced $A[F_\mu]$-module morphism.
The map $f:V\to U(M)$ is a morphism in $\Sei_\infty(A)$ if and
only if $f(x)=\widetilde{f}(s(1-z)x)$ for all $x\in V$.
\end{lemma}

Let us deduce statement 1.~above. By the definition of $B(V)$ there is
an exact sequence
\begin{equation*}
0\to V[F_\mu]\xrightarrow{1-s(1-z)}V[F_\mu]\xrightarrow{q} B(V)\to 0
\end{equation*}
so $q(x)=q(s(1-z)x)$ for all $x\in V$. By
Lemma~\ref{relating_Seifert_structures}, $\theta_V=q|_V$ is a morphism
of Seifert modules.

We turn now to statement 2. Since $V$ generates $B(V)$ as an
$A[F_\mu]$-module, an $A[F_\mu]$-module morphism $g:B(V)\to M$ satisfies
$f=U(g)\theta_V$ if and only if $g$ fits into the diagram 
\begin{equation*}
\xymatrix@R=3ex@C=7ex{
0 \ar[r] & V[F_\mu] \ar[r]^{1-s(1-z)} & V[F_\mu]
\ar[r]^q 
\ar[dr]_{\widetilde{f}} 
& B(V) 
\ar@{.>}[d]^g 
\ar[r] & 0 
\\
& & & M &
}
\end{equation*}
Since $f$ is a morphism of Seifert modules we have
$f(x)=\widetilde{f}(s(1-z)x)$ for
all $x\in V$ by Lemma~\ref{relating_Seifert_structures}. So
$\widetilde{f}\sigma=0$, and
therefore there exists unique $g:B(V)\to M$ such that
$gq=\widetilde{f}$. It follows that there exists unique $g$ such that
$f=U(g)\theta_V$. Thus we have established both 1.~and 2.~assuming
Lemma~\ref{relating_Seifert_structures}.
\begin{proof}[Proof of Lemma~\ref{relating_Seifert_structures}]
The Seifert module structure on $U(M)$ is given
by~(\ref{Seifert_structure_on_F-link_module}) above so $f$ is
a Seifert morphism if and only if \\ \noindent
a) $\omega\gamma^{-1}f(x)=f(sx)$ \ and \\ \noindent b) $\gamma
p_i\gamma^{-1}f(x)=f(\pi_ix)$ for each $x\in V$. 

To prove the `if' part of Lemma~\ref{relating_Seifert_structures},
suppose $f(x)=\widetilde{f}(s(1-z)x)$.\\
\noindent a) The equations $\displaystyle{\omega\gamma^{-1}\left(\sum_{i=1}^\mu
(1-z_i)x_i\right)=\sum_{i=1}^\mu x_i}$ and 
\begin{equation*}
f(x)=\widetilde{f}(s(1-z)x)=\sum_{i=1}^\mu \widetilde{f}(s(1-z_i)\pi_ix) =
\sum_{i=1}^\mu(1-z_i)f(s\pi_i x).
\end{equation*}
imply that
$\displaystyle{\omega\gamma^{-1}f(x)}=
\displaystyle{\omega\gamma^{-1}\sum_{i=1}^\mu(1-z_i)f(s\pi_ix)}=
\displaystyle{\sum_{i=1}^\mu f(s\pi_ix)}=f(sx)$. \\ \noindent
b) Observe that 
\begin{equation*}
f(\pi_ix)=\widetilde{f}(s(1-z)\pi_ix)=\sum_{j=1}^\mu
\widetilde{f}(s(1-z_j)\pi_j\pi_ix)=(1-z_i)f(s\pi_ix).
\end{equation*}
while
\begin{align*}
\gamma p_i\gamma^{-1}f(x)=\gamma p_i\gamma^{-1}\widetilde{f}(s(1-z)x)&= 
\gamma p_i\gamma^{-1}\sum_{j=1}^\mu(1-z_j)f(s\pi_jx) \\
&=(1-z_i)f(s\pi_ix).
\end{align*}
Thus $f(\pi_ix)=\gamma p_i\gamma^{-1}f(x)$.

To prove the ``only if'' part of
Lemma~\ref{relating_Seifert_structures}, suppose we have a) and b)
above. Now 
\begin{align*}
\widetilde{f}(s(1-z)x)&=\sum_{i=1}^\mu (1-z_i)f(s\pi_ix) \\
&= \sum_{i=1}^\mu (1-z_i)(\omega\gamma^{-1})(\gamma p_i\gamma^{-1})f(x) \\
&= \sum_{i=1}^\mu (1-z_i)\omega p_i\gamma^{-1}f(x) \\
&= f(x).\qedhere
\end{align*}
\end{proof}
This completes the proof of Proposition~\ref{B_is_left_adjoint_to_U}.
\end{proof}

\subsection{Localization}\label{localization_of_cat}
When one passes from Seifert modules to $F_\mu$-link modules, certain
Seifert modules disappear altogether; following Farber we shall call
such modules primitive. 
\begin{definition}\label{define_primitive}
Let $\Prim_\infty(A)$ denote the full subcategory of $\Sei_\infty(A)$
containing precisely the modules $V$ such that $B(V)=0$. Modules in
$\Prim_\infty(A)$ will be called primitive.
\end{definition}
For example, if $\rho(s)=0$ or
$\rho(s)=1$ then 
\begin{equation*}
(1-\rho(s))y^{-1} + \rho(s)y:V[F_\mu(y^2)y] \to V[F_\mu(y^2)]
\end{equation*} 
is an isomorphism and therefore has zero cokernel. A module in
$\Sei_\infty(A)$ with $\rho(s)=0$ or $1$ will be called {\it trivially
primitive}.
We show in Section~\ref{section:structure_of_primitives} that all the
primitive Seifert modules in $\Sei_\infty(A)$
can be ``built'' from trivially primitive modules. If $A$ is
semi-simple Artinian then a similar result applies when one restricts 
attention to the category $\Sei(A)$ of
representations of $P_\mu$ by finitely generated projective
$A$-modules: Every
primitive in $\Sei(A)$ can be ``built'' from a finite number of
trivially primitive modules in $\Sei(A)$ (see
Proposition~\ref{structure_of_artinian_finite_primitives}).
This statement is not true for all rings $A$; one must consider
primitives which exhibit a kind of nilpotence. Such primitives were
described by Bass, Heller and Swan when $\mu=1$ (see also
Ranicki~\cite{Ran03}). The general case $\mu\geq1$
will be analyzed in a subsequent paper (joint work with A.Ranicki). 

In the present section we construct an equivalence between
$\Bl_\infty(A)$ and a quotient category
$\Sei_\infty(A)/\Prim_\infty(A)$. This quotient is an example of
universal localization for categories; the objects in
$\Sei_\infty(A)/\Prim_\infty(A)$ are the same as the objects in
$\Sei_\infty(A)$ but the morphisms are different. The universal property is
that a morphism in $\Sei_\infty(A)$ whose kernel and cokernel are in
$\Prim_\infty(A)$ has an inverse in $\Sei_\infty(A)/\Prim_\infty(A)$. A more
detailed construction of the quotient appears in
Section~\ref{section:construction_of_quotient}.
We proceed to derive an equivalence between $\Bl(A)$ and a corresponding quotient
of $\Sei(A)$.
\begin{definition}\label{universal_localization_of_Sei_infty}
The functor $F:\Sei_\infty(A)\to\Sei_\infty(A)/\Prim_\infty(A)$ is
the universal functor which makes invertible all morphisms whose kernel
and cokernel are primitive. In other words, any functor which makes
these morphisms invertible factors uniquely through $F$.
\end{definition}
We outline in Section~\ref{section:construction_of_quotient} one construction of $F$ which will be convenient for
our purposes; see Gabriel~\cite{Gab62} or Srinivas~\cite[Appendix
B.3]{Sri96} for further details. A more general construction can be
found in Gabriel and Zisman~\cite{GabZis67} or
Borceux~\cite[Ch5]{Bor94_vol1}. It follows directly from the definition
that the localization $F$ is unique (up to unique isomorphism).

Applying Definition~\ref{universal_localization_of_Sei_infty} to the
functor $B:\Sei_\infty(A)\to\Bl_\infty(A)$,
there is a unique 
functor $\overline{B}$ such that $B=\overline{B}F$: 
\begin{equation*}
\xymatrix{ 
{\Sei_\infty(A)} \ar@/^2pc/[rr]^B \ar[r]_F & {\displaystyle{\frac{\Sei_\infty(A)}{\Prim_\infty(A)}}}
\ar[r]_{\overline{B}} & {\Bl_\infty(A)} \ar@/^2pc/[ll]^U
}
\end{equation*}

Proposition~\ref{B_is_left_adjoint_to_U} stated that
$B:\Sei_\infty(A)\to\Bl_\infty(A)$ is left adjoint to the forgetful
functor $U$. We deduce in the next proposition that $B$ satisfies the
same universal property as $F$, but only ``up to
natural isomorphism''.
If $f$ is a morphism in $\Sei_\infty(A)$ let us write $f\in\Xi$ if the
kernel and cokernel of $f$ both lie in $\Prim_\infty(A)$.
\begin{proposition}\label{B_is_universal_localization_up_to_equivalence}
If $G:\Sei_\infty(A)\to\cy{B}$ is a functor which sends
every morphism in $\Xi$ to an invertible morphism then there is a
functor 
\begin{equation*}
\widetilde{G}:\Bl_\infty(A)\to\cy{B}
\end{equation*}
such that
$\widetilde{G}B$ is naturally isomorphic to $G$. The functor
$\widetilde{G}$ is unique up to natural isomorphism.
\end{proposition}
\begin{proof}
We prove uniqueness first.
If there is a natural isomorphism $G\simeq\widetilde{G}B$ then
$GU\simeq\widetilde{G}BU\simeq \widetilde{G}$ by
Corollary~\ref{deductions_from_left_adjointness}. 

To prove existence we must show that if $\widetilde{G}=GU$ then
$\widetilde{G}B\simeq G$. Indeed,
by Corollary~\ref{deductions_from_left_adjointness}
$B(\theta_V):B(V)\to BUB(V)$ is an isomorphism for each
$V\in\Sei_\infty(A)$. Since $B$ respects exact sequences we have
$\theta_V\in\Xi$.
It follows that $G(\theta):G\to GUB=\widetilde{G}B$ is a natural isomorphism.
\end{proof}
The following is an immediate consequence of the fact that $F$ and $B$
have the same universal property (up to natural isomorphism):
\begin{corollary}\label{Bbar_is_equiv_on_infty_cats}
The functor $\overline{B}:\Sei_\infty(A)/\Prim_\infty(A)\to
\Bl_\infty(A)$ is an equivalence.\hspace*{\fill}\qed 
\end{corollary}
We turn now to the categories $\Sei(A)$ and $\Bl(A)$.
\begin{definition}
Let $\Sei(A)/\Prim_\infty(A)~\subset~\Sei_\infty(A)/\Prim_\infty(A)$
denote the full subcategory whose objects are precisely the
modules in $\Sei(A)$ (i.e.~the modules which are finitely generated
and projective as $A$-modules).
\end{definition}
There is a commutative diagram of functors
\begin{equation}\label{factor_B_through_quotient_diagram}
\begin{gathered}
\xymatrix@R=3ex{
{\Sei_\infty(A)} \ar@/^2pc/[rr]^B \ar[r]_F &
{\displaystyle{\frac{\Sei_\infty(A)}{\Prim_\infty(A)}}} \ar[r]_{\overline{B}} &
{\Bl_\infty(A)} \\ 
{\Sei(A)}  \ar[u] \ar@/_2pc/[rr]_B \ar[r]_F &
{\displaystyle{\frac{\Sei(A)}{\Prim_\infty(A)}}} \ar[u]\ar[r]_{\overline{B}} &
{\Bl(A)} \ar[u]
}
\end{gathered}
\end{equation}
in which all the vertical arrows are inclusions of full subcategories.

\begin{theorem}\label{Bbar_is_equiv_on_fgproj_cats}
The functor $\overline{B}:\Sei(A)/\Prim_\infty(A)\to \Bl(A)$ is an
equivalence of categories.
\end{theorem}
We will use the following general lemma in the proof of
Theorem~\ref{Bbar_is_equiv_on_fgproj_cats}. Recall that a functor
$G:\cy{C}\to\cy{D}$ is called full and faithful if it induces an
isomorphism $\Hom_\cy{C}(V,V')\to\Hom_\cy{D}(G(V),G(V'))$ for every
pair of objects $V,V'\in \cy{C}$.
\begin{lemma}\label{describe_cat_equivalence}
A functor $G:\cy{C}\to\cy{D}$ is
an equivalence of categories if and only if 
$G$ is full and faithful and every object
in $\cy{D}$ is isomorphic to $G(V)$ for some $V\in\cy{C}$.
\end{lemma}
\begin{proof}
See for example Borceux~\cite[Prop 3.4.3]{Bor94_vol1}.
\end{proof}
It follows from Corollary~\ref{Bbar_is_equiv_on_infty_cats} and
Lemma~\ref{describe_cat_equivalence} that
\begin{equation*}
\overline{B}:\Sei(A)/\Prim_\infty(A)\to\Bl(A)
\end{equation*}
is full and faithful. Theorem~\ref{Bbar_is_equiv_on_fgproj_cats} is therefore a
consequence of the following proposition:
\begin{proposition}\label{B_is_surjective}
Every module in $\Bl(A)$ is isomorphic to $B(V)$ for some $V\in\Sei(A)$.
\end{proposition}
\begin{proof}
By definition, every module $M\in\Bl(A)$ has presentation
\begin{equation*}
0\to V[F_\mu]\xrightarrow{\sigma} V[F_\mu]\to M \to 0
\end{equation*}
where $V$ is a finitely generated projective $A$-module and
$\epsilon(\sigma):V\to V$ is an isomorphism. Given any $A$-module $W$
there is a canonical isomorphism 
\begin{equation*}
\Hom_{A[F_\mu]}(V[F_\mu],W[F_\mu])\cong \Hom_A(V,W)[F_\mu]
\end{equation*}
and in particular $\sigma$ can be expressed uniquely as a sum $\sum_{w\in
F_\mu}\sigma_ww$ with each $\sigma_w\in\Hom_A(V,V)$.

\begin{lemma}\label{linearize}
Every $M\in\Bl(A)$ is isomorphic to the cokernel of an endomorphism
$\sigma:V[F_\mu]\to V[F_\mu]$ of the form
\begin{equation*}
\sigma=1+\sigma_1(1-z_1)+\cdots+\sigma_\mu(1-z_\mu)
\end{equation*} 
where $V$ is finitely generated and projective and
$\sigma_1,\cdots,\sigma_\mu\in\Hom_A(V,V)$.
\end{lemma}
\begin{proof}[Proof of Lemma]
By the definition of $\Bl(A)$, the module $M$ is isomorphic to
the cokernel of some map $\sigma:V[F_\mu]\to V[F_\mu]$ where $V$ is
finitely generated and projective. The idea of this proof is to reduce
the support of $\sigma$ to $\{1,z_1,\cdots,z_\mu\}\subset F_\mu$ at
the expense of replacing $V$ by a larger finitely generated projective module.
Note first that  
\begin{equation*}
\Coker(\sigma)\cong
\Coker\left(\begin{matrix}
\sigma & 0 \\
0 & 1
\end{matrix}\right):(V\oplus V')[F_\mu]\to (V\oplus V')[F_\mu]
\end{equation*}
where $V'$ is any $A$-module and $1$ denotes the identity
morphism. The equation
\begin{equation}\label{Gaussian_elimination}
\left(\begin{matrix}
1 & -b \\ 
0 & 1
\end{matrix}\right)
\left(\begin{matrix}
a+bc & 0 \\ 
0    & 1
\end{matrix}\right)
\left(\begin{matrix}
1 & 0 \\ 
c & 1
\end{matrix}\right)
=
\left(\begin{matrix}
a & -b \\ 
c & 1
\end{matrix}\right),
\end{equation}
therefore implies that $\Coker\left(\begin{matrix}
a & -b \\ 
c & 1
\end{matrix}\right)$ is isomorphic to $\Coker(a+bc)$.
Repeated application of equation~(\ref{Gaussian_elimination}) implies
that $M$ is isomorphic to the cokernel of an endomorphism
$\sigma=\sigma_0+\sum_{i=1}^\mu\sigma_i^+z_i +
\sum_{i=1}^\mu\sigma_i^-z_i^{-1}$ with $\sigma_0$, $\sigma_i^+$ and
$\sigma_i^-$ in $\Hom_A(V,V)$ for some $V$. For each of the indices
$i=1,\cdots,\mu$ in turn, one can apply the identity
$\Coker(\sigma)=\Coker(\sigma z_i)$ followed by further
equations~(\ref{Gaussian_elimination}).
One obtains an identity $M\cong \Coker(\beta)$ where
$\beta=\beta_0+\beta_1z_1+\cdots +\beta_\mu z_\mu$
and $\beta_i\in\Hom_A(V,V)$ for some finitely
generated projective module $V$ over $A$. By
Lemma~\ref{Sato_condition}, $\epsilon(\beta)$ 
is an isomorphism. Let $\sigma=\epsilon(\beta)^{-1}\beta$. Now
$\epsilon(\sigma)=1$ and so 
\begin{equation*}
\sigma=1+\sigma_1(z_1-1)+\cdots \sigma_\mu (z_\mu-1)
\end{equation*}
for some
$\sigma_1,\cdots,\sigma_\mu\in\Hom_A(V,V)$.
This completes the proof of Lemma~\ref{linearize}.
\end{proof}

We may now finish the proof of Proposition~\ref{B_is_surjective}.
If 
\begin{equation*}
\sigma=1+\sum_i\sigma_i(1-z_i)
\end{equation*}
then the equation 
\begin{align*}
&\left(\begin{smallmatrix}
1 & 0 & \cdots & 0 \\
1 & 1 & \cdots & 0 \\
\vdots & \vdots & \ddots & \vdots \\
1 & 0 & \cdots & 1
\end{smallmatrix}\right)
\left(\begin{smallmatrix}
1 & \sigma_2(1-z_2) & \cdots & \sigma_\mu(1-z_\mu) \\
0 & 1 & \cdots & 0 \\
\vdots & \vdots & \ddots & \vdots \\
0 & 0 & \cdots & 1
\end{smallmatrix}\right)
\left(\begin{smallmatrix}
\sigma & 0 & \cdots & 0 \\
0 & 1 & \cdots & 0 \\
\vdots & \vdots & \ddots & \vdots \\
0 & 0 & \cdots & 1
\end{smallmatrix}\right)
\left(\begin{smallmatrix}
1 & 0 & \cdots &  0 \\
-1 & 1 & \cdots & 0 \\
\vdots & \vdots & \ddots & \vdots \\
-1 & 0 & \cdots & 1
\end{smallmatrix}\right) \\
&= 
\left(\begin{smallmatrix}
1 +\sigma_1(1-z_1) & \sigma_2(1-z_2) & \cdots & \sigma_\mu(1-z_\mu) \\
\sigma_1(1-z_1) & 1+\sigma_2(1-z_2) & \cdots & \sigma_\mu(1-z_\mu) \\
\vdots & \vdots & \ddots & \vdots \\
\sigma_1(1-z_1) & \sigma_2(1-z_2) & \cdots & 1+\sigma_\mu(1-z_\mu)
\end{smallmatrix}\right)
\end{align*}
implies that 
\begin{equation*}
\Coker(\sigma)\cong \Coker\left(1-s(1-z) : V^{\oplus\mu}[F_\mu] \to
V^{\oplus\mu}[F_\mu]\right)
\end{equation*}
where $\pi_i$ acts as projection on the $i$th component of $V^{\oplus
\mu}$ and $s$ acts as
\begin{equation*}
\left(\begin{matrix}
\sigma_1 & \sigma_2 & \cdots & \sigma_\mu \\
\sigma_1 & \sigma_2 & \cdots & \sigma_\mu \\
\vdots & \vdots & \ddots & \vdots \\
\sigma_1 & \sigma_2 & \cdots & \sigma_\mu
\end{matrix}\right).
\end{equation*}
Thus $M\cong B(V^{\oplus\mu})$.
\end{proof}
\noindent This completes the proof of Theorem~\ref{Bbar_is_equiv_on_fgproj_cats}.\hspace*{\fill}\qed

\subsection{Duality in the quotient}\label{section:duality_on_quotient_arbitrary_A}
Having established that
$\overline{B}:\Sei(A)/\Prim_\infty(A)\to\Bl(A)$ is an equivalence,
we may use $\overline{B}$ to give duality structure to
$\Sei(A)/\Prim_\infty(A)$ and make the lower part
of~(\ref{factor_B_through_quotient_diagram}) a commutative diagram of
duality-preserving functors.
Since the objects in $\Sei(A)/\Prim_\infty(A)$ coincide with those in
$\Sei(A)$ we define $F(V)^*=F(V^*)$ and $i_{F(V)}=F(i_V):F(V)\to
F(V)^{**}$ where 
\begin{equation*}
F:\Sei(A)\to
\Sei(A)/\Prim_\infty(A)
\end{equation*}
is the canonical functor.
If $f:V\to V'$ is a morphism in $\Sei(A)/\Prim_\infty(A)$ let
\begin{equation}\label{duality_on_morphisms_in_Sei/Prim_infty}
f^*=\overline{B}^{-1}(\Phi_V^{-1}\overline{B}(f)^\wedge\Phi_{V'}):{V'}^*\to V^*.
\end{equation}
It is easy to see that $\functor^*$ is a contravariant functor and
that $i^*_Vi_{V^*}=id_{V^*}$ for all $V$ so $\Sei(A)/\Prim_\infty(A)$
is a hermitian category.

Recall that the composite of duality-preserving functors is defined by
\begin{equation*}
(G,\Psi,\eta)\circ (G',\Psi',\eta')= (GG',\Psi G(\Psi') ,\eta\eta').
\end{equation*}
\begin{proposition}\label{composite_of_d_pres_functors}
The duality-preserving functor 
\begin{equation*}
(B,\Phi,-1):\Sei(A)\to \Bl(A)
\end{equation*}
coincides with the composite $(\overline{B},\Phi,-1)\circ (F,\id,1)$.
\end{proposition}
\begin{proof}
It follows from equation~(\ref{duality_on_morphisms_in_Sei/Prim_infty})
and Proposition~\ref{extend_B_to_duality_pres_functor} that
$(\overline{B},\Phi,-1)$ is a duality-preserving functor. 

By definition $F(V)^*=F(V^*)$ and $i_{F(V)}= F(i_V)$; to show that
\begin{equation*}
(F,\id,1):\Sei(A)\to \Sei(A)/\Prim_\infty(A)
\end{equation*}
is a duality-preserving functor we must check 
that $F(f)^*=F(f^*)$ for each morphism $f:V\to V'$ in $\Sei(A)$. Indeed, 
\begin{align*}
F(f)^*&=\overline{B}^{-1}(\Phi_V^{-1}\overline{B}F(f)^\wedge\Phi_{V'}) \\
&= \overline{B}^{-1}(\Phi_V^{-1}B(f)^\wedge\Phi_{V'}) \\
&= \overline{B}^{-1}(B(f^*)) \quad\quad \mbox{(since $\Phi$ is natural)} \\
&= F(f^*).
\end{align*}
It is easy to verify that $(B,\Phi,-1)=(\overline{B},\Phi,-1)\circ (F,\id,1)$.
\end{proof}

\begin{proposition}\label{Bbar_is_eq_of_hermitian_cats}
Let $\zeta=1$ or $-1$. The duality-preserving functor 
\begin{equation*}
(\overline{B},\Phi,-1):\frac{\Sei(A)}{\Prim_\infty(A)}\to \Bl(A)
\end{equation*}
 is an equivalence of hermitian categories and induces an
isomorphism of Witt groups
\begin{equation}
\label{isomorphism_of_witt_groups_induced_by_Bbar}
\overline{B}:W^\zeta\left(\frac{\Sei(A)}{\Prim_\infty(A)}\right)
\to W^{-\zeta}(\Bl(A)).
\end{equation}
\end{proposition}
\begin{proof}
Since $\overline{B}$ is an equivalence of categories (by
Theorem~\ref{Bbar_is_equiv_on_fgproj_cats} above) it follows that
$(\overline{B},\Phi,-1)$ is an
equivalence of hermitian categories (see Proposition II.7
of~\cite{She03mem}). It also follows that $\overline{B}$
preserves limits and colimits so $\overline{B}$ preserves exact
sequences and hence induces a
homomorphism~(\ref{isomorphism_of_witt_groups_induced_by_Bbar}) of
Witt groups. By
Lemma~\ref{nat_isomorphic_dpres_funct_give_same_Witt_homo} of
Appendix~\ref{section:naturality_of_B} this homomorphism is an
isomorphism~(\ref{isomorphism_of_witt_groups_induced_by_Bbar}).
\end{proof}
\subsection{Structure of Primitives}\label{section:structure_of_primitives}
Recall that a module $(V,\rho)$ in $\Sei_\infty(A)$ is called
trivially primitive if $\rho(s)=0$ or $\rho(s)=1$. In this section we prove
that every primitive module in $\Sei_\infty(A)$ is composed of
trivially primitive modules. 
\begin{lemma}\label{s=0_submodule_lemma}
If $(V,\rho)\in\Sei_\infty(A)$ and there exists a non-zero element
$x\in V$ such that $\rho(s\pi_i)x=0$ for all $i$ then $V$ has a
non-zero submodule $(V',\rho')$ such that $\rho'(s)=0$. 
\end{lemma}
\begin{proof}
Note that $x=\sum\pi_i x$ and at least one of the terms $\pi_ix$ must
be non-zero.  
Choose non-zero $V'=A\pi_ix$ and define $\rho'$ by 
\begin{equation*}
\rho'(\pi_j)=\begin{cases}
1 \ \mbox{if $j= i$} \\
0 \ \mbox{if $j\neq i$} 
\end{cases},
\quad 
\rho'(s)=0.
\end{equation*}
Now $(V',\rho')$ is the required non-zero submodule of $(V,\rho)$.
\end{proof}
\begin{lemma}\label{primitive_has_s=1_or_s=0_submodule}
If $(V,\rho)$ is primitive and non-zero then there exists a non-zero
trivially primitive submodule $(V',\rho')$.
\end{lemma}
\begin{proof} (compare Lemma~7.10c in Farber~\cite{Far92})
Since $(V,\rho)$ is primitive, 
\begin{equation*}
\rho(1-s)y^{-1} + \rho(s)y : V[F_\mu(y^2)y]\to V[F_\mu(y^2)]
\end{equation*}
is an isomorphism with inverse $\alpha$ say. Now $\alpha$ can be written as
a finite sum $\sum_{w\in S}\alpha_w w$ where $S$ is a finite subset of
$\bigcup_{i=1}^\mu F_\mu(y^2)y_i$ and $\alpha_w:V\to V$ has non-zero image in
$\pi_iV$ for each $w\in S$. Choose an element $w\in S$ whose
expression in reduced form as a product of letters $y_i^{\pm}$ is of
maximal length. We consider two cases: \\[1ex]
\noindent {\bf Case 1: $w=w'y_i$} for some $w'\in F_\mu(y^2)$ and
some $i$. The equation 
\begin{equation}\label{linear_sigma_inverted}
((1-\rho(s))y^{-1} + \rho(s)y)\alpha=1
\end{equation}
implies that $\rho(s\pi_j)\alpha_w=0$ for each $j$. Any element $x$ in
the image of 
$\alpha_w$ satisfies the conditions of Lemma~\ref{s=0_submodule_lemma}
so $(V,\rho)$ has a non-zero submodule $(V',\rho')$ with
$\rho'(s)=0$. \\[1ex]
\noindent{\bf Case 2: $w=w'y_i^{-1}$} for some $w'\in F_\mu(y^2)$. The
equation~(\ref{linear_sigma_inverted}) implies that
$\rho((1-s)\pi_i)\alpha_w=0$ for each $i$. By
Lemma~\ref{s=0_submodule_lemma} there is a non-zero submodule
$(V',\rho')$ with $\rho'(1-s)=0$ or in other words $\rho'(s)=1$.
\end{proof}
Recall that a module $V$ is called simple if there are no
submodules other than $0$ and $V$. The following remark is a
consequence of Lemma~\ref{primitive_has_s=1_or_s=0_submodule}.
\begin{remark}\label{Every_simple_primitive_module_is_trivially_primitive}
Every simple primitive module is trivially primitive.\hspace*{\fill}\qed
\end{remark}
\begin{definition}
If $\cy{A}$ is an abelian
category then a non-empty full subcategory $\cy{E}\subset \cy{A}$ is called a
{\it Serre subcategory} if for every exact sequence $0\to V\to V'\to V''\to
0$ in $\cy{A}$ one has
\begin{equation*}
V'\in\cy{E}~\Leftrightarrow~V\in\cy{E}~\mbox{and}~V''\in\cy{E}~.
\end{equation*}
\end{definition}
Note that every Serre subcategory of an abelian category is again an
abelian category. Since $B$ preserves exact sequences and arbitrary
direct sums $\Prim_\infty(A)$ is a Serre subcategory of
$\Sei_\infty(A)$ and is closed under direct sums.
\begin{lemma}\label{maximal_primitive_submodule}
Suppose $\cy{E}\subset\cy{A}$ is a Serre subcategory of an abelian category
and $\cy{E}$ is closed under arbitrary direct sums. Every module
$V\in\cy{A}$ contains a unique maximal submodule $U\leq V$ such that
$U\in\cy{E}$. If $U'\leq V$ and $U'\in \cy{E}$ then $U'\leq U$.
\end{lemma}
\begin{proof}
Let $U$ be the sum in $V$ of all the submodules $U_i\leq
V$ with $U_i\in\cy{E}$. Since $U$ is a factor module of $\bigoplus
U_i$, one finds $U\in\cy{E}$.
\end{proof}
\begin{proposition}\label{identify_prim_infty}
The category $\Prim_\infty(A)$ is the smallest Serre subcategory of
$\Sei_\infty(A)$ which a) contains the trivially primitive modules and
b) is closed under arbitrary direct sums.
\end{proposition}
\begin{proof}
Let $\cy{P}_\infty(A)$ denote the smallest full subcategory of
$\Sei_\infty(A)$ satisfying the
conditions of the Proposition. Now $\Prim_\infty(A)$ satisfies these conditions
so $\cy{P}_\infty(A)\subset\Prim_\infty(A)$.

Conversely, we must show that $\Prim_\infty(A)\subset
\cy{P}_\infty(A)$. Suppose $B(V)=0$. Let $W\leq V$ be the maximal
submodule such that $W\in\cy{P}_\infty(A)$ (the module $W$ exists by
Lemma~\ref{maximal_primitive_submodule}). Now
$B(V/W)=0$ since $B$ respects exact sequences so
Lemma~\ref{primitive_has_s=1_or_s=0_submodule} implies that
either $V/W=0$ or there is a non-zero submodule $V'$ of $V/W$ which
lies in $\cy{P}_\infty(A)$. In the latter case, let $p:V\to
V/W$ denote the projection and note the exact sequence
\begin{equation*}
0\to W\to p^{-1}(V')\xrightarrow{p|} V' \to 0.
\end{equation*}
Since $W\in\cy{P}_\infty(A)$ and $V'\in\cy{P}_\infty(A)$ we have
$p^{-1}(V)\in\cy{P}_\infty(A)$ which contradicts the maximality of $W$. Thus
$V/W=0$ and hence $V=W$, so $V$ lies in $\cy{P}_\infty(A)$.
\end{proof}

\subsection{Construction of the
quotient}\label{section:construction_of_quotient}
We outline next a construction of $\Sei_\infty(A)/\Prim_\infty(A)$. We 
will use this construction in
Section~\ref{section:Proof_that_B_induces_Witt_isomorphism} to 
show that $B:W^\zeta(\Sei(A))\to W^{-\zeta}(\Bl(A))$ is an
isomorphism when $A$ is a semi-simple Artinian ring.
The notion of Serre subcategory was defined in the preceding
section. Let us note some basic properties:
\begin{lemma}\label{add_and_coadd_in_Serre_subcat}
Suppose $\cy{A}$ is an abelian category and $\cy{E}$ is
a Serre subcategory. Suppose $V\in\cy{A}$, $W\leq V$ and $W'\leq V$.
\begin{enumerate}
\item\label{add_in_Serre_subcat} If $W\in\cy{E}$ and $W'\in\cy{E}$ then $W+W'\in \cy{E}$.
\item\label{coadd_in_Serre_subcat} If $V/W\in\cy{E}$ and
$V/W'\in\cy{E}$ then $V/(W\cap W')\in\cy{E}$.
\end{enumerate}
\end{lemma}
\begin{proof}
\ref{add_in_Serre_subcat}. There is an exact sequence 
\begin{equation*}
0\to W\to W+W'\to (W+W')/W\to 0
\end{equation*}
and  $(W+W')/W$ is isomorphic to $W'/(W\cap W')\in\cy{E}$. Hence
$W+W'\in\cy{E}$. \\ \noindent
\ref{coadd_in_Serre_subcat}. There is an exact sequence 
\begin{equation*}
0\to W/(W\cap W')\to
V/(W\cap W')\to V/W\to 0.
\end{equation*}
Now $W/(W\cap W')$ is isomorphic to $(W+W')/W'$
which is contained in $V/W'$ and so $W/(W\cap W')\in\cy{E}$ and hence
$V/W\cap W'\in\cy{E}$.
\end{proof} 
We may now recall a construction for the quotient of an abelian
category by a Serre subcategory.
See Gabriel~\cite{Gab62} or Srinivas~\cite[Appendix B.3]{Sri96} for 
further details. 

Suppose $\cy{A}$ is an abelian category and $\cy{E}$
is a Serre subcategory. The symbol $\cy{A}/\cy{E}$ will denote a
category with the same objects as $\cy{A}$ but different groups of
morphisms. To define $\Hom_{\cy{A}/\cy{E}}(V,V')$, consider the pairs
$(W,U')$ where $W\leq V$, $U'\leq V'$, $V/W\in\cy{E}$ and
$U'\in\cy{E}$. One says that $(W_1,U'_1)\leq(W_2,U'_2)$ if $W_2\leq
W_1$ and $U'_1\leq U'_2$ [note the directions of
inclusion]. Lemma~\ref{add_and_coadd_in_Serre_subcat} 
above implies that these pairs are a directed set. Indeed, given
pairs $(W_1,U'_1)$ and $(W_2,U'_2)$ one finds $(W_1,U'_1)\leq (W_1\cap
W_2,U'_1+U'_2)$ and $(W_2,U'_2)\leq (W_1\cap W_2,U'_1+U'_2)$.
The following definition can now be made:
\begin{equation}\label{define_morphisms_in_quotient}
\Hom_{\cy{A}/\cy{E}}(V,V')=\varinjlim_{(W,U')} \Hom_\cy{A}(W,V'/U').
\end{equation}
We leave to the reader the definition of composition of morphisms and
the canonical functor $F:\cy{A}\to\cy{A}/\cy{E}$. 
Proofs of the following statements can be found in the references cited above:
{\renewcommand{\labelenumi}{{\normalfont (\alph{enumi})}}
\begin{enumerate}
\item The quotient category $\cy{A}/\cy{E}$ is an abelian category and
$F$ is an exact additive functor.
\item\label{morphisms_become_isomorphisms} If $f$ is a morphism in $\cy{A}$ then $F(f)$ is an isomorphism
if and only if $\Coker(f)\in\cy{E}$ and $\Ker(f)\in\cy{E}$.
\end{enumerate}}
\noindent In particular if $V$ is a module in $\cy{A}$ then $F(V)\cong0$ if and
only if $V\in \cy{E}$.

As we indicated in earlier sections, the functor $F:\cy{A}\to
 \cy{A}/\cy{E}$ is universal with respect to property~(b). In detail, if  
 $G:\cy{A}\to \cy{B}$ makes invertible every morphism whose kernel
 and cokernel lie in $\cy{E}$ then there is a unique functor
 $\widetilde{G}:\cy{A}/\cy{E}\to\cy{B}$ such that $\widetilde{G}F=G$. 
In particular the functor
\begin{equation*}
F:\Sei_\infty(A)\to\Sei_\infty(A)/\Prim_\infty(A)
\end{equation*}
satisfies Definition~\ref{universal_localization_of_Sei_infty}.
Let us be explicit about $\widetilde{G}$: 

If $V$ is an module in  $\cy{A}/\cy{E}$ then one writes
 $\widetilde{G}(V)=G(V)$. Every morphism
 $f\in\Hom_{\cy{A}/\cy{E}}(V,V')$ is represented by some 
 $\overline{f}\in\Hom_\cy{A}(W,V'/U')$ with $U'\in\cy{E}$ and $V/W\in\cy{E}$ If $i :W\to V$ and $p:V'\to V'/U'$ denote the
 canonical monomorphism and epimorphism respectively one must define
\begin{equation*}
\widetilde{G}(f)=G(p)^{-1}G(\overline{f})G(i)^{-1}:G(V)\to G(V').
\end{equation*} 

In our particular example Lemma~\ref{maximal_primitive_submodule}
provides one simplification in our description of the quotient category
$\Sei_\infty(A)/\Prim_\infty(A)$. If $V\in \Sei_\infty(A)$ let us call
a submodule $W\leq V$ {\it coprimitive} if $V/W\in\Prim_\infty(A)$. 
\begin{lemma}\label{half_simplify_quotient}
If $V,V'\in\Sei_\infty(A)$ and $U'$ denotes the maximal primitive
submodule of $V'$ then
\begin{equation*}
\Hom_{\Sei_\infty(A)/\Prim_\infty(A)}(V,V') = \varinjlim_W \Hom_{\Sei_\infty(A)}(W,V'/U')
\end{equation*}
where the direct limit is over coprimitive submodules $W$ of $V$. \hspace*{\fill}\qed
\end{lemma}
\noindent Note that there is not in general a minimal coprimitive in
$V$; the functor $B$ does not respect infinite limits and an infinite
intersection of coprimitives is not in general coprimitive (but see
Lemma~\ref{max_prim_min_coprim_assuming_finiteness} below).

\subsection{Global dimension zero}\label{section:global_dim_zero}
{\it In this section the ring $A$ will be assumed semi-simple and
Artinian or, in other words, a finite product of matrix rings over
division rings}. The basic theory of semi-simple
Artinian rings can be found in many algebra textbooks
(e.g.~Lam~\cite[\S1-4]{Lam91} or Lang~\cite[Ch.XVII]{Lan93}). In particular,
all $A$-modules are projective and $\Sei(A)$ is an abelian category with
ascending and descending chain conditions; these facts lead to simplifications
of results in Sections~\ref{section:structure_of_primitives}
and~\ref{section:construction_of_quotient} above.
We show that the primitive modules in
$\Sei(A)$ are composed of a finite number of simple trivially
primitive modules
(Proposition~\ref{structure_of_artinian_finite_primitives}) and give a
simplified description of the 
hermitian category $\Sei(A)/\Prim_\infty(A)$. We shall consider
semi-simple Artinian rings again in
Section~\ref{section:Proof_that_B_induces_Witt_isomorphism} but it is
not essential to read the present section before
Section~\ref{section:Proof_that_B_induces_Witt_isomorphism}.

The key lemma we will need is the following:
\begin{lemma}\label{max_prim_min_coprim_assuming_finiteness}
Suppose $\cy{A}$ is an abelian category with ascending and descending
chain conditions and $\cy{E}$ is a Serre subcategory.
\begin{enumerate}
\item Every module $V\in\cy{A}$ contains a unique maximal
  submodule in $\cy{E}$ which contains all others in $\cy{E}$.
\item
Every module $V\in\cy{A}$ contains a unique submodule $W\leq V$ which
is minimal with respect to the property $V/W\in\cy{E}$. If
$V/W'\in\cy{E}$ then $W\leq W'$.
\end{enumerate}
\end{lemma} 
\begin{proof}
1.~Since $\cy{A}$ has the ascending chain condition there is a
  submodule $U\leq V$ 
  which is maximal with respect to the property $U\in\cy{E}$. In other
  words, if
  $U\leq U'\leq V$ and $U'\in\cy{E}$ 
  then $U'=U$. If $U'$ is any other submodule in $\cy{E}$ then $U+U'\in\cy{E}$
  by Lemma~\ref{add_and_coadd_in_Serre_subcat} so
  $U+U'=U$ and hence $U'\leq U$. \\ 
\noindent 2.~Since $\cy{A}$ has the descending chain condition there is a
  submodule $W\leq V$ 
  which is minimal with respect to the property $V/W\in\cy{E}$ (i.e.~if
  $W'\leq W\leq V$ and $V/W'\in\cy{E}$ 
  then $W'=W$). If $V/W'\in\cy{E}$ then $V/(W\cap W')\in\cy{E}$
  by Lemma~\ref{add_and_coadd_in_Serre_subcat} so
  $W\cap W'=W$ and hence $W\leq W'$.
\end{proof}
Recall that if $V\in\Sei(A)$, a submodule $W\leq V$ is called coprimitive
if $V/W$ is primitive. Since $A$ is Artinian and Noetherian,
Lemma~\ref{max_prim_min_coprim_assuming_finiteness}
implies that there is a maximal primitive submodule $U\leq V$ and a minimal
 coprimitive submodule $W\leq V$ for each $V\in\Sei(A)$. 
\subsubsection{Structure of
 Primitives}\label{section:structure_of_primitives_artinian_case}
\begin{definition}
Let $\Prim(A)$ denote the intersection of $\Prim_\infty(A)$ and
$\Sei(A)$. In other words, $\Prim(A)\subset\Sei(A)$ is the full
subcategory containing those modules $V$ such that $B(V)=0$.
\end{definition}
Note that $\Prim(A)$ is both a Serre subcategory and a hermitian
subcategory of $\Sei(A)$. Moreover,
$\Sei(A)/\Prim(A)=\Sei(A)/\Prim_\infty(A)$.
\begin{proposition}\label{structure_of_artinian_finite_primitives}
The category $\Prim(A)$ is the smallest Serre subcategory of $\Sei(A)$
which contains the trivially primitive modules in $\Sei(A)$.
\end{proposition}
\begin{proof}
We proceed as in the proof of Proposition~\ref{identify_prim_infty},
using Lemma~\ref{max_prim_min_coprim_assuming_finiteness} in place of
Lemma~\ref{maximal_primitive_submodule}. Let
$\cy{P}(A)$ denote the smallest Serre subcategory of $\Sei(A)$ which
contains all the trivially primitive modules.
To show that $\cy{P}(A)\subset \Prim(A)$ it suffices to observe that
$\Prim(A)$ is a Serre subcategory which contains these modules.

Conversely, to show $\Prim(A)\subset \cy{P}(A)$ suppose $V\in\Prim(A)$.
There exists, by Lemma~\ref{max_prim_min_coprim_assuming_finiteness}, a
maximal submodule $U\leq V$ such that $U\in\cy{P}(A)$. Now $B(V/U)=0$
so either $V=U$ or by
Lemma~\ref{primitive_has_s=1_or_s=0_submodule} $V/U$ has a
non-zero trivially primitive submodule $U'$. If
$p:V\twoheadrightarrow V/U$ is the canonical map then the exact
sequence $0\to U\to p^{-1}(U')\to U'\to0$ implies that
$p^{-1}(U')\in\cy{P}(A)$ contradicting the maximality of $U$. Thus
$V=U$ and $V\in\cy{P}(A)$.
\end{proof}
\subsubsection{Construction of the
quotient}\label{section:structure_of_quotient_gd0}
With the benefit of
Lemma~\ref{max_prim_min_coprim_assuming_finiteness} we can give a simpler
description of the quotient category than Lemma~\ref{half_simplify_quotient}.
\begin{lemma}\label{fully_simplify_quotient}
The morphisms in $\Sei(A)/\Prim(A)$ are 
\begin{equation}\label{morphisms_in_Sei/Prim}
\Hom_{\Sei(A)/\Prim(A)}(V,V')=\Hom_{\Sei(A)}(W,V'/U').
\end{equation}
where $W\leq V$ is the minimal coprimitive and $U'\leq V'$ is the maximal
primitive.\hspace*{\fill}\qed
\end{lemma}

We simplify next the hermitian structure on $\Sei(A)/\Prim(A)$.
We have seen that the duality-preserving functor
$(B,\Phi,-1):\Sei(A)\to\Bl(A)$ factors through an equivalence of
hermitian categories $(\overline{B},\Phi,-1)$:
\begin{equation*}
\xymatrix@C=8ex{
\Sei(A) \ar@/^2pc/[rr]^{(B,\Phi,-1)} \ar[r]_{(F,\id,1)} &
{\frac{\Sei(A)}{\Prim(A)}} \ar[r]_{(\overline{B},\Phi,-1)} & {\Bl(A)}
}
\end{equation*}
(Theorem~\ref{Bbar_is_equiv_on_fgproj_cats} and
Proposition~\ref{composite_of_d_pres_functors} above). The duality
functor on the quotient $\Sei(A)/\Prim(A)$ was defined in
Section~\ref{section:duality_on_quotient_arbitrary_A} above by
$F(V)^*=F(V^*)$ and by equation~(\ref{duality_on_morphisms_in_Sei/Prim_infty}).  
Using the assumption that $A$ is Artinian we can
re-interpret equation~(\ref{duality_on_morphisms_in_Sei/Prim_infty}).
Suppose that $f\in\Hom_{\Sei(A)/\Prim(A)}(V,V')$. As usual, let $W$
denote the minimal coprimitive submodule of $V$ and let $U'$ denote
the maximal primitive submodule of $V'$. The morphism $f$ is identified 
with some $\overline{f}\in\Hom_{\Sei(A)}(W,V'/U')$. 
Since $\functor^*$ preserves exact sequences the
following are exact
\begin{align}\label{dualized_short_exact_sequences}
&0\to (V/W)^*\to V^* \to W^* \to 0      \\
&0\to (V'/U')^*\to (V')^*\to (U')^*\to 0.
\end{align}
Now $(B,\Phi,1)$ is a duality-preserving functor, so for each
$V\in\Sei(A)$ one has $B(V)=0$ if and only if $B(V^*)=0$. Thus
$\Prim(A)$ is a hermitian subcategory of $\Sei(A)$ and in particular
$(U')^*$ and $(V/W)^*$ are primitive. It follows that $(V/W)^*$ is the
maximal primitive
in $V^*$ and $(V'/U')^*$ is the minimal coprimitive in $(V')^*$. Since
$(F,\id,1)$ is
a duality-preserving functor, $f^*$ is represented by
\begin{equation}\label{duality_in_hermitian_quotient}
{\overline{f}}^*\in\Hom((V'/U')^*,W^*).
\end{equation}
\section{Equivalence of Invariants}
\label{section:Proof_that_B_induces_Witt_isomorphism}
Cobordism invariants of $F_\mu$-links have been defined in two
different ways in Sections~\ref{section:Blanchfield_form_invariants}
and~\cite{She03mem}.
In this section we use the duality-preserving functor $(B,\Phi,-1)$
which was studied in
Section~\ref{section:covering_construction} to relate the two
approaches, proving Theorems~\ref{main_covering_construction_theorem}
and~\ref{theorem_relating_invariants:vague}.
To prove Theorem~\ref{theorem_relating_invariants:vague} we show that
the functor $B$ respects each of the three steps laid out in
Sections~\ref{section:overview_of_BD_form_invariants} and \ref{section:cobordism_invariants_seifert}. A more detailed
version of Theorem~\ref{theorem_relating_invariants:vague} is set
out in Theorem~\ref{theorem_relating_invariants:precise} below.
\subsection{Proof of
Theorem~\ref{main_covering_construction_theorem}}\label{section:pf_of_Witt_isomorphism_thm} 
Suppose $\cy{A}$ is an abelian category with ascending and descending
chain conditions and $\cy{E}$ is a Serre
subcategory. Let $F:\cy{A}\to\cy{A}/\cy{E}$ denote the quotient functor.
Recall that a module $V$ in $\cy{A}$ is called simple if $V$ is not
isomorphic to $0$ and $V$ does not have submodules other than $0$ and $V$.
\begin{lemma}\label{simple_goes_to_comes_from_simple}
\begin{enumerate}
\item If $V\in\cy{A}$ is simple then either $V\in\cy{E}$ or $F(V)$ is simple.
\item Every simple module in $\cy{A}/\cy{E}$ is isomorphic to $F(V)$
  for some simple module $V\in\cy{A}$ which does not lie in $\cy{E}$.
\end{enumerate}
\end{lemma}
\begin{proof}
1.~Suppose $V\in\cy{A}$ is simple, $V\notin\cy{E}$ and $i:V'\to F(V)$ is the inclusion
  of a submodule in $\cy{A}/\cy{E}$. Now $i$ is represented by some morphism
  $\overline{i}:W'\to V$ where $W'\leq V'$ and
  $V'/W'\in\cy{E}$. Either $\overline{i}=0$ in which case $V'\cong 0$
  in $\cy{A}/\cy{E}$ or $\overline{i}$ is an epimorphism which implies
  that $V'=F(V)$ (recall that $F$ is exact). \\
2.~Every module in $\cy{A}/\cy{E}$ is $F(V)$ for some module
  $V\in\cy{A}$. Suppose $F(V)$ 
  is simple. Now $V$ has a finite filtration 
$0=V_0\leq V_1\leq\cdots \leq V_n=V$ where each quotient $V_i/V_{i-1}$
is a simple module. Since $F$ respects exact sequences
  $F(V_i/V_{i-1})=0$ for all $i\in\{1,\cdots,n\}$ except one, for
  which there is an isomorphism $F(V_i/V_{i-1})\cong F(V)$. This
  module $V_i/V_{i-1}$ does not lie in $\cy{E}$.
\end{proof}

Suppose now that $\cy{A}$ and $\cy{A}/\cy{E}$ are hermitian categories and the
quotient functor extends to a duality-preserving functor
\begin{equation*}
(F,\id,1):\cy{A}\to\cy{A}/\cy{E}.
\end{equation*}
\begin{lemma}\label{quotient_lemma_with_duality}
\begin{enumerate}
\item The Serre subcategory $\cy{E}$ is a hermitian subcategory.
\item Let $\zeta=1$ or $-1$. If $V\in\cy{A}$ is simple and $V\notin\cy{E}$ then $V$ is
$\zeta$-self-dual if and only if $F(V)$ is $\zeta$-self-dual.
\end{enumerate}
\end{lemma}
\begin{proof}
1.~If $V\in\cy{E}$ then $F(V^*)=F(V)^*\cong0\in\cy{A}/\cy{E}$ so
  $V^*\in\cy{E}$. \\
2.~To prove the ``only if'' part it suffices to recall that for $\phi:V\to
  V^*$ one has $F(\phi^*)=F(\phi)^*$. For the ``if'' part, note also
  that 
\begin{equation*}
F:\Hom_{\cy{A}}(V,V^*)\to\Hom_{\cy{A}/\cy{E}}(V,V^*)
\end{equation*} is an isomorphism.
\end{proof}
\begin{proposition}\label{Witt_isomorphism_for_hermitian_quotient}
Suppose $\cy{A}$ and $\cy{A}/\cy{E}$ are hermitian categories and
\begin{equation*}
(F,\id,1):\cy{A}\to\cy{A}/\cy{E}
\end{equation*}
is a duality-preserving functor. 
For each $\zeta$-self-dual simple module
$V\in\cy{A}$ such that $V\notin\cy{E}$ there is a canonical isomorphism
\begin{equation*}
W^\zeta(\cy{A}|_V)\cong W^\zeta((\cy{A}/\cy{E})|_V).
\end{equation*}
If $\cy{A}$ has ascending and descending chain conditions then
there is a canonical isomorphism
\begin{equation*}
W^\zeta(\cy{A})\cong W^\zeta(\cy{E})\oplus W^\zeta(\cy{A}/\cy{E}).
\end{equation*}
\end{proposition}
\begin{proof}
If $V\in\cy{A}$ is a simple module and $V\notin\cy{E}$ then $F(V)$ is
simple by part 1.~of
Lemma~\ref{simple_goes_to_comes_from_simple}. Every module 
in $\cy{A}|_V$ is a direct sum of copies of $V$ so by
equation~(\ref{define_morphisms_in_quotient})
the restriction $F:\cy{A}|_V\to(\cy{A}/\cy{E})|_{F(V)}$ is a
full and faithful functor and hence an equivalence of
categories. 

By part 2.~of
Lemma~\ref{quotient_lemma_with_duality}, $V$ is
$\zeta$-self-dual if and only if $F(V)$ is  
$\zeta$-self-dual, in which case
$(F,\id,1):F:\cy{A}|_V\to(\cy{A}/\cy{E})|_{F(V)}$ is an equivalence of
hermitian categories and induces an isomorphism
\begin{equation}\label{restricted_witt_homomorphism}
W^\zeta(\cy{A}|_V)\to W^\zeta\left((\cy{A}/\cy{E})|_{F(V)}\right).
\end{equation}

To prove the last sentence of the Lemma, note first that
by part 1.~of Lemma~\ref{quotient_lemma_with_duality}, $\cy{E}$ is a
hermitian subcategory of $\cy{A}$. Theorem~\ref{hermitian_devissage}
provides canonical decompositions  
\begin{align*}
&W^\zeta(\cy{A})\cong \bigoplus W^\zeta(\cy{A}|_V) \\ 
&W^\zeta(\cy{E})\cong \bigoplus W^\zeta(\cy{E}|_V) \\
&W^\zeta(\cy{A}/\cy{E})\cong \bigoplus
W^\zeta\left((\cy{A}/\cy{E})|_V\right)
\end{align*}
where the right hand side of each identity has one summand for each
isomorphism class of $\zeta$-self-dual simple modules $V$. 

By part 2.~of Lemma~\ref{simple_goes_to_comes_from_simple} and part
2.~of Lemma~\ref{quotient_lemma_with_duality} every
summand of $W^\zeta(\cy{A}/\cy{E})$ is the isomorphic image of
$W^\zeta(\cy{A}|_V)$ for some simple $\zeta$-self-dual module
$V$ in $\cy{A}$.

On the other hand, if $V\in\cy{E}$ is simple and $\zeta$-self-dual
then $(F,\id,1)$ sends $W^\zeta(\cy{A}|_V)$ to zero. The last
sentence of the Lemma follows.
\end{proof}
In our application, we set $\cy{A}=\Sei(A)$ and
$\cy{E}=\Prim(A)$ where $A$ is semi-simple Artinian. Recall that
$\Prim(A)=\Prim_\infty(A)\cap \Sei(A)$ is an abelian category with ascending
and descending chain conditions and $\Sei(A)/\Prim(A)=\Sei(A)/\Prim_\infty(A)$.
\begin{lemma}\label{Witt_group_of_primitives_is_zero}
\begin{enumerate}
\item None of the simple primitive modules in $\Prim(A)$ are self-dual.
\item $W^\zeta(\Prim(A))=0$.
\end{enumerate}
\end{lemma}
\begin{proof}
By Remark~\ref{Every_simple_primitive_module_is_trivially_primitive}
above, every simple primitive module is trivially primitive. If
$(V,\rho)\in\Prim(A)$ then $\rho(s)=0$ if and only if $\rho^*(s)=1$ so
none of the simple trivially primitive modules are self-dual. Thus
part 1.~is proved, and part 2.~follows immediately from
Theorem~\ref{hermitian_devissage}.
\end{proof}
\begin{proof}[Proof of Theorem~\ref{main_covering_construction_theorem}]
By Proposition~\ref{composite_of_d_pres_functors} the
duality-preserving functor $(B,\Phi,-1)$ is the composite
$(\overline{B},\Phi,-1)\circ (F,\id,1)$. 
Setting $\cy{A}=\Sei(A)$ and $\cy{E}=\Prim(A)$ in
Proposition~\ref{Witt_isomorphism_for_hermitian_quotient},
and invoking also
Lemma~\ref{Witt_group_of_primitives_is_zero}, we learn that
$(F,\id,1)$ induces an isomorphism
\begin{equation*}
W^\zeta(\Sei(A))\to W^\zeta\left(\frac{\Sei(A)}{\Prim(A)}\right).
\end{equation*}
By Theorem~\ref{Bbar_is_equiv_on_fgproj_cats}, $(\overline{B},\Phi,1)$ is an
equivalence and hence induces an isomorphism
\begin{equation*}
W^\zeta\left(\frac{\Sei(A)}{\Prim(A)}\right)\to W^{-\zeta}(\Bl(A)).
\end{equation*}
(see Lemma~\ref{nat_isomorphic_dpres_funct_give_same_Witt_homo} in
Appendix~\ref{section:naturality_of_B}). Thus $(B,\Phi,-1)$ induces an
isomorphism
\begin{equation*}
W^\zeta(\Sei(A))\to W^{-\zeta}(\Bl(A)).
\end{equation*}
This completes the proof of Theorem~\ref{main_covering_construction_theorem}.
\end{proof}
\subsection{Proof of Theorem~\ref{theorem_relating_invariants:vague}}\label{section:pf_of_eq_of_invariants}
We prove in this section that the functor $(B,\Phi,-1)$ identifies the
invariants defined in~\cite{She03mem} with those of
Section~\ref{section:Blanchfield_form_invariants}. More precisely, we
prove the following theorem:
\begin{theorem}{\text{\bf (Equivalence of
invariants)}}\label{theorem_relating_invariants:precise}
\begin{enumerate}
\item If $V\in\Sei(\Q)$ is simple and $(-1)^q$-self-dual then
$B(V)\in\Bl(\Q)$ is simple and $(-1)^{q+1}$-self-dual. 
\item Every simple $(-1)^{q+1}$-self-dual module $M\in\Bl(\Q)$ is
isomorphic to $B(V)$ for some simple $(-1)^q$-self-dual module $V\in\Sei(\Q)$.
\item If $V\in\Sei(\Q)$ and $B(V)\in\Bl(\Q)$ are simple then the
functor $B$ induces an isomorphism of rings
$\displaystyle{B:\End_{\Sei(\Q)}(V)\xrightarrow{\cong}\End_{\Bl(\Q)}(B(V))}$.
\item Suppose $V\in\Sei(\Q)$ is simple and $b:V\to V^*$ is a
$\zeta$-hermitian form. The ring isomorphism in part 3.~respects
involutions. Explicitly, if $f\in\End_{\Sei(\Q)}V$ then 
$\displaystyle{B(b^{-1}f^*b)=(\Phi_V B(b))^{-1}B(f)^\wedge\Phi_V B(b)}$.
\item Suppose $W\in\Sei(\Z)$ and $\phi:W\to W^*$ is a
$(-1)^q$-hermitian form. The dimension modulo $2$, signatures,
discriminant, Hasse-Witt invariant and Lewis $\theta$-invariant of
\begin{equation*}
\Theta_{V,b}p_{V}\,[\Q\otimes_\Z (W,\phi)]\in
W^1(\End_{\Sei(\Q)}(V))
\end{equation*}
coincide (if defined) with the corresponding invariants of
\begin{equation*}
\Theta_{(B(V),-\Phi_V B(b))}p_{B(V)}\,[\Q\otimes_\Z (B(W),\Phi_W B(\phi))]\in
W^1(\End_{\Bl(\Q)}(B(V))).
\end{equation*}
\end{enumerate}
\end{theorem}

Recall that $B=\overline{B}\circ F$ and
$\overline{B}:\Sei(\Q)/\Prim(\Q)\to\Bl(\Q)$ is an equivalence of
categories. In parts
1.~through 3.~of Theorem~\ref{theorem_relating_invariants:precise} it
therefore suffices to prove 
corresponding statements with the functor $F$ in place of $B$ and
$(-1)^q$ in place of $(-1)^{q+1}$: \\ \noindent
1. %
The statement follows from 
part 1.~of Lemma~\ref{simple_goes_to_comes_from_simple}, part 2.~of
Lemma~\ref{quotient_lemma_with_duality} and part 1.~of
Lemma~\ref{Witt_group_of_primitives_is_zero}. \\ \noindent
2. %
The statement follows from
part 2.~of Lemma~\ref{simple_goes_to_comes_from_simple} and part 2.~of
Lemma~\ref{quotient_lemma_with_duality}. \\ \noindent
3. %
This is a consequence of equation~(\ref{morphisms_in_Sei/Prim}). \\ \noindent
4. Since $B(b^{-1}f^*b)=B(b^{-1})B(f^*)B(b)$ it suffices to prove that
\begin{equation*}
B(f^*)=\Phi_V^{-1}B(f)^\wedge\Phi_V.
\end{equation*}
This equation is a consequence of the fact that $\Phi$ is a natural
isomorphism.

The proof of part 5.~of
Theorem~\ref{theorem_relating_invariants:precise} is slightly more involved. 
Recall from proposition~\ref{B_respects_coefficient_change} that $B$
respects changes of coefficients and, in particular, that the
inclusion of $\Z$ in $\Q$ induces the
commutative diagram (\ref{B_respects_Z_into_Q}).
One must check that $(B,\Phi,-1)$ respects each of the three steps
in the definitions of the $F_\mu$-link invariants (see
Section~\ref{section:overview_of_BD_form_invariants}).
\vspace*{0.5ex}

\noindent {\bf Devissage:}
Let $V\in\Sei(\Q)$ be a $\zeta$-self-dual
simple module. If $W$ is isomorphic to a
direct sum of copies of $V$ then $B(W)$ is isomorphic to a direct sum
of copies of $B(V)$. Hence the image of
$W^\zeta(\Sei(\Q)|_V)$ under $B$ lies in
$W^{-\zeta}(\Bl(\Q)|_{B(V)})$ and there is a commutative diagram of isomorphisms
\begin{equation}\label{B_respects_devissage}
\begin{gathered}
\xymatrix@R=2ex@C=9ex{
W^\zeta(\Sei(\Q)) \ar[r]^B\ar@{<->}[d] & W^{-\zeta}(\Bl(\Q))\ar@{<->}[d] \\
{\displaystyle{\bigoplus_V W^\zeta(\Sei(\Q)|_V)}}
\ar[r]_(0.45){\bigoplus B|}
& {\displaystyle{\bigoplus_{V}
W^{-\zeta}(\Bl(\Q)|_{B(V)})}}. 
}
\end{gathered}
\end{equation}
where the direct sums are indexed by the isomorphism classes of simple
$\zeta$-self-dual modules in $\Sei(\Q)$.
\vspace*{0.5ex}

\noindent {\bf Morita Equivalence:}
Suppose $V\in\Sei(\Q)$ is a simple module and
\begin{equation*}
b:V\to V^*
\end{equation*}
is a non-singular $\zeta$-hermitian form. Let us denote the
endomorphism rings $E=\End_{\Sei(\Q)}V$ and $E'=\End_{\Bl(\Q)}B(V)$.
By Corollary~\ref{Witt_naturality_of_h_morita_eq} above the
duality-preserving functor $(B,\Phi,-1)$ induces a
commutative diagram
\begin{equation}\label{B_respects_HME:Witt_level}
\begin{gathered}
\xymatrix@R=3ex@C=7ex{
W^\zeta(\Sei(\Q)|_V) \ar[d]_{\Theta_{V,b}}
\ar[r]^(.45)B & W^{-\zeta}(\Bl(\Q)|_{B(V)}) 
\ar[d]^{\Theta_{B(V),-\Phi_V B(b)}} \\
W^1(E) \ar[r]_B & W^1(E').
}
\end{gathered}
\end{equation}
\noindent {\bf Invariants:}
The isomorphism $E\to E'$ in part 4.~induces isomorphisms
between the target groups for the invariants in part 5. For example,
if $E$ and $E'$ are commutative
with trivial involution then the discriminant $\Delta$ of
$\Theta_{V,b}p_{V}\,[\Q\otimes_\Z (W,\phi)]$ lies
in $E/E^2$ and the functor $B$ induces an isomorphism $E/E^2\to
E'/(E')^2$. The word ``coincide'' in part 5.~is understood to mean that the 
image of $\Delta$ in $E'/(E')^2$ is equal to the discriminant of
$\Theta_{(B(V),-\Phi_V B(b))}p_{B(V)}\,[\Q\otimes_\Z (B(W),\Phi_W
B(\phi))]$. 

The isomorphism $B:E\to E'$ 
of rings with involution induces an isomorphism $W^1(E)\to W^1(E')$. 
We leave it to the reader to check that if $\alpha\in W^1(E)$ then all
the listed invariants of $\alpha$ coincide (in this sense) with the
corresponding invariants of $B(\alpha)\in W^1(E')$. Further details of
the invariants can be found in chapter 11 of~\cite{She03mem}. 

This completes the proof of part 5.~and hence of
theorems~\ref{theorem_relating_invariants:vague}
and~\ref{theorem_relating_invariants:precise}.\hspace*{\fill}\qed
\begin{appendix}
\section{Naturality of constructions}\label{section:naturality_of_B}
In this appendix we prove naturality theorems for 
the covering construction $(B,\Phi,-1)$ and for hermitian Morita
equivalence, proving Propositions~\ref{B_respects_coefficient_change}
and~\ref{naturality_of_h_morita_eq} above.

To compare duality-preserving functors one requires the following definition.
\begin{definition}
Suppose $(G,\Psi,\eta):\cy{C}\to\cy{D}$ and
$(G',\Psi',\eta):\cy{C}\to\cy{D}$ are duality-preserving functors
between hermitian categories $\cy{C}$ and $\cy{D}$. A natural
transformation $\alpha:(G,\Psi,\eta)\to (G',\Psi',\eta)$ is a natural
transformation between the underlying functors, $\alpha:G\to G'$, such that
\begin{equation}\label{condition_for_nat_trans}
\Psi_V=\alpha_V^*\Psi'_V\alpha_{V^*}
\end{equation}
for each object $V\in\cy{C}$.
\end{definition}
If $\alpha:G\to G'$
is a natural isomorphism between the underlying functors and $\alpha$
satisfies~(\ref{condition_for_nat_trans}) then $\alpha^{-1}:G'\to G$
also satisfies~(\ref{condition_for_nat_trans}) so $\alpha$ is in fact a
natural isomorphism of duality-preserving functors. 

We noted in Lemma~\ref{duality_preserving_functor_induces} that an exact
duality-preserving functor induces a homomorphism of Witt groups. The
following lemma says that naturally isomorphic duality-preserving
functors induce the same homomorphism on Witt groups. 
\begin{lemma}\label{nat_isomorphic_dpres_funct_give_same_Witt_homo}
Suppose $(G,\Psi,\eta),(G',\Psi',\eta):\cy{C}\to\cy{D}$ are
duality-preserving functors which respect exact sequences and
 $\alpha:(G,\Psi,\eta)\to(G',\Psi',\eta)$ is a natural isomorphism. If
$(V\ ,\ \phi:V\to V^*)$ is a hermitian form in $\cy{C}$ then there is a
natural isomorphism between the induced hermitian forms
$(G(V),\Psi G(\phi))\cong (G'(V),\Psi'G'(\phi))$. Let $\zeta=1$ or
$-1$. The duality-preserving
functors $(G,\Psi,\eta)$ and $(G',\Psi',\eta)$ induce the same homomorphism
of Witt groups $W^\zeta(\cy{C})\to W^{\zeta\eta}(\cy{D})$.
\end{lemma}
\begin{proof}
In the diagram
\begin{equation*}
\xymatrix@R=2ex@C=5ex{
G(V)\ar[d]_{\alpha_V} \ar[r]^{G(\phi)} & G(V^*) \ar[d]_{\alpha_{V^*}}\ar[r]^{\Psi_V} & G(V)^* \\
G'(V) \ar[r]_{G'(\phi)} & G'(V^*) \ar[r]_{\Psi'_V} & G'(V)^* \ar[u]_{\alpha_V^*}
}
\end{equation*}
the left-hand square commutes by the naturality of $\alpha$ while the
right-hand square commutes because $\alpha$ satisfies
equation~(\ref{condition_for_nat_trans}). The Lemma follows easily.
\end{proof}
It is a consequence of Lemma~\ref{nat_isomorphic_dpres_funct_give_same_Witt_homo} that an equivalence of hermitian
categories induces an isomorphism of Witt groups.
\subsection{The covering construction}
In this section we prove that the covering construction $B$ respects
changes to coefficients
(Proposition~\ref{B_respects_coefficient_change}). We need one more
observation which is straightforward to verify: 
\begin{lemma}\label{nat_isos_for_composed_ring_homos}
If $A\to A'\to A''$ are ring homomorphisms then the diagram
\begin{equation*}
\xymatrix@R=0ex@C=2ex{
A''\otimes_A V^* \ar[dd] \ar[r] & A''\otimes_{A'}(A'\otimes_A
V^*) \ar[dr] & \\
& &  A''\otimes_{A'}(A'\otimes_A V)^* \ar[dl]  \\
(A''\otimes_A V)^* & (A''\otimes_{A'}(A'\otimes_A V))^*\ar[l]}
\end{equation*}
 of natural isomorphisms is commutative.\hspace*{\fill}\qed
\end{lemma}
\begin{proof}[Proof of Proposition~\ref{B_respects_coefficient_change}]
Suppose $A\to A'$ is a ring homomorphism and $V$ is a module in $\Sei(A)$. The
natural isomorphism 
\begin{equation*}
(A'\otimes_A V)[F_\mu]\to A'[F_\mu]\otimes_{A[F_\mu]} (V[F_\mu])
\end{equation*}
induces a natural isomorphism (see Lemma~\ref{Bl(functor)_is_functor})
\begin{equation*}
\{\alpha_V\}_{V\in\Sei(A)}: B(A'\otimes V) \to A'[F_\mu]\otimes_{A[F_\mu]}
B(V) .
\end{equation*}

We aim to show that $\alpha$ is
a natural isomorphism between duality-preserving functors
\begin{equation*}
(B,\Phi,-1)\circ (A'\otimes_A\functor,\Pi,1) \to
(A'[F_\mu]\otimes_{A[F_\mu]}\functor,\Upsilon,1)\circ (B,\Phi,-1).
\end{equation*}
Applying Lemma~\ref{nat_isos_for_composed_ring_homos} to both
composites in the commutative square
\begin{equation*}
\xymatrix@=3ex{A_{\phantom{\mu}} \ar[r]\ar[d] & A[F_\mu] \ar[d] \\
A'_{\phantom{\mu}}\ar[r] & A'[F_\mu]}
\end{equation*}
of ring homomorphisms one obtains commutative diagrams
\begin{equation*}
\xymatrix@R=2ex@C=6ex{(A'\otimes_A V^*)[F_\mu] \ar[r]\ar[d] &
(A'\otimes_A V)^*[F_\mu] \ar[r]^{\gamma} &
(A'\otimes_A V[F_\mu])^* \\
A'[F_\mu]\otimes_{A[F_\mu]}(V^*[F_\mu]) \ar[r]_{\delta} &
A'[F_\mu]\otimes_{A[F_\mu]}(V[F_\mu])^*
\ar[r] &
(A'[F_\mu] \otimes_{A[F_\mu]} V[F_\mu])^* \ar[u]
}
\end{equation*}
where $\gamma=\pm\Pi_{A'\otimes V}(1-z^{\pm})$ and
$\delta=1\otimes\pm\Pi_V(1-z^{\pm})$ and hence the commutative diagram
\begin{equation*}
\xymatrix@R=3ex@C=6ex{
B(A'\otimes_A V^*) \ar[r]^{B(\Pi_V)}\ar[d]_{\alpha_{V^*}} & B((A'\otimes_A V)^*)
\ar[r]^{\Phi_{A'\otimes V}} & B(A'\otimes
V)^\wedge \\
A'[F_\mu]\otimes_{A[F_\mu]} B(V^*) \ar[r]_{1\otimes\Phi_V} & A'[F_\mu]\otimes_{A[F_\mu]} B(V)^\wedge
\ar[r]_{\Upsilon_{B(V)}} & (A'[F_\mu]\otimes_{A[F_\mu]} B(V))^\wedge \ar[u]_{(\alpha_V)^\wedge~.}
}
\end{equation*}
Thus $\alpha$ is a natural transformation between duality-preserving
functors as claimed. It follows by
Lemma~\ref{nat_isomorphic_dpres_funct_give_same_Witt_homo} that the
diagram~(\ref{Witt_diagram:B_respects_coeff_change}) of Witt group
homomorphisms commutes. The proof of
Proposition~\ref{B_respects_coefficient_change} is complete.
\end{proof}
\subsection{Hermitian Morita Equivalence}\label{section:natuality_of_h_Morita_eq}
In this section we prove that hermitian Morita equivalence respects
duality-preserving functors (Proposition~\ref{naturality_of_h_morita_eq}).
Let $(G,\Psi,\eta'):\cy{C}\to\cy{D}$ denote a duality-preserving functor and let $M\in\cy{C}$, $E=\End_\cy{C}M$ and
$E'=\End_{\cy{D}}G(M)$. 

We shall define a natural isomorphism between the composite
functors 
\begin{multline*}
\alpha:(E'\otimes_E\functor,\Pi,1)\circ(\Hom(M,\functor),\Omega^b,\eta) \\
\xrightarrow{\simeq}
(\Hom(G(M),\functor),\Omega^{\eta'\Psi_MG(b)},\eta\eta')\circ
(G,\Psi,\eta')
\end{multline*}
If $N\in\cy{C}|_M$ then $\Hom(M,N)$ is a left $E$-module for the
 action 
\begin{equation*}
f.\theta=\theta \overline f= \theta b^{-1}f^*b
\end{equation*}
where $f\in E$ and $\theta\in\Hom(M,N)$. The group $\Hom(G(M),G(N))$ is
 regarded as a left $E'$-module in the same way. Define
\begin{equation}\label{nat_iso_for_naturality_of_h_Morita_Eq}
\begin{aligned}
\alpha_N&:E'\otimes_E\Hom_{\cy{C}}(M,N)\to
\Hom_{\cy{D}}(G(M),G(N)) \\
&f\otimes \gamma \mapsto f.G(\gamma)=G(\gamma)\overline{f}=G(\gamma)(\eta'\Psi_MG(b))^{-1}f^*(\eta'\Psi_MG(b)).
\end{aligned}
\end{equation}
Since $\alpha_N$ is an isomorphism in the case $N=M$ it follows that
$\alpha_N$ is an isomorphism for all $N\in \cy{C}|_M$. It is easy to
see that $\{\alpha\}_{N\in\cy{C}|_M}$ is a natural transformation 
\begin{equation*}
(E'\otimes_E\functor)\circ\Hom(M,\functor)\to\Hom(G(M),\functor)\circ
G. 
\end{equation*}
One must check that 
$\alpha$ is a natural transformation of duality-preserving
functors. By equations~(\ref{how2_compose_dpres_functors})
and~(\ref{condition_for_nat_trans}) one must show that
\begin{equation*}
\alpha_N^*\Omega^{\eta'\Psi_MG(b)}_{G(N)}\Hom(G(M),\Psi_N)\alpha_{N^*}=\Pi_{\Hom(M,N)}(1\otimes\Omega^b_N)
\end{equation*}
This equation can be checked by direct calculation, substituting the
formulae~(\ref{nat_iso_for_naturality_of_h_Morita_Eq}),
(\ref{natural_isomorphism_Phi})
and~(\ref{nat_iso_in_coeff_change_dpres_funct}) for $\alpha$, $\Omega$
and $\Pi$ respectively and applying the naturality of $\Phi$ and the
equation~(\ref{duality_functor}). This completes the proof of
Proposition~\ref{naturality_of_h_morita_eq}.
\end{appendix}

Desmond Sheiham \\
International University Bremen \\
Campus Ring 1 \\
Bremen 28759 \\
Germany \\
des@sheiham.com
\end{document}